\documentclass{amsart}
\RequirePackage{amssymb}
\usepackage[all]{xy}

\theoremstyle{plain}
 
\newtheorem{theorem}{Theorem}[section]
\newtheorem{proposition}[theorem]{Proposition}
\newtheorem{lemma}[theorem]{Lemma}
\newtheorem{corollary}[theorem]{Corollary}
 
\theoremstyle{definition}
 
\newtheorem{definition}[theorem]{Definition}

\newtheorem{notations}[theorem]{Notations}
\newtheorem{properties}[theorem]{Properties}
\newtheorem{claim}[theorem]{Claim}
\newtheorem{remark}[theorem]{Remark}
\newtheorem{remarks}[theorem]{Remarks}
\newtheorem{example}[theorem]{Example}

\newtheorem{assumptions}[theorem]{Assumptions}
\newcommand{\C}{\mathbf{C}}
\newcommand{\R}{\mathbf{R}}
\newcommand{\N}{\mathbf{N}}
\newcommand{\Z}{\mathbf{Z}}
\newcommand{\q}{\mathcal{Q}}
\newcommand{\D}{\mathrm{D}}
\newcommand{\K}{\mathrm{K}}
\newcommand{\M}{\mathrm{M}}
\newcommand{\Cinf}{\mathcal{C}}
\newcommand{\F}{\mathcal{F}}

\renewcommand{\O}{\mathcal{O}}
\renewcommand{\S}{\mathcal{S}}

\newcommand{\fA}{\mathcal{A}}
\newcommand{\fB}{\mathcal{B}}
\newcommand{\fAintro}{\mathcal{A}}
\newcommand{\fBintro}{\mathcal{B}}
\newcommand{\fC}[2]{\mathcal{R}^{#1}_{#2}}
\newcommand{\fcun}[2]{\mathcal{RA}^{#1}_{#2}}
\newcommand{\fcde}[2]{\mathcal{RB}^{#1}_{#2}}
\newcommand{\fctr}[2]{\mathcal{RC}^{#1}_{#2}}
\newcommand{\fcqu}[2]{\mathcal{RD}^{#1}_{#2}}
\newcommand{\fcci}[2]{\mathcal{RE}^{#1}_{#2}}
\newcommand{\fcs}[2]{\mathcal{RS}^{#1}_{#2}}
\newcommand{\fcsun}[2]{\mathcal{RT}^{#1}_{#2}}
\newcommand{\fck}[2]{\mathcal{RK}^{#1}_{#2}}
\newcommand{\fckun}[2]{\mathcal{RL}^{#1}_{#2}}
\newcommand{\fckde}[2]{\mathcal{RM}^{#1}_{#2}}
\newcommand{\fD}{\mathcal{A}}
\newcommand{\E}{\mathcal{E}}
\renewcommand{\H}{\mathcal{H}}

\renewcommand{\to}{\rightarrow}
\newcommand{\from}{\leftarrow}

\newcommand{\xto}[2][]{\xrightarrow[#1]{#2}}
\newcommand{\xfrom}[2][]{\xleftarrow[#1]{#2}}
\newcommand{\isoto}[1][]{\xrightarrow[#1]{\raisebox{-2pt}[0pt][0pt]{$\sim$}}}
\newcommand{\isofrom}[1][]{\xleftarrow[#1]{\raisebox{-2pt}[0pt][0pt]{$\sim$}}}

\newcommand{\sect}{\Gamma}
\newcommand{\rsect}{R\Gamma}
\newcommand{\im}{\operatorname{im}}
\newcommand{\codim}{\operatorname{codim}}

\newcommand{\RHom}{\operatorname{RHom}}
\newcommand{\Hom}{\operatorname{Hom}}

\newcommand{\Ext}{\operatorname{Ext}}

\renewcommand{\hom}{\mathcal{H}om}
\newcommand{\Rhom}{R\mathcal{H}om}
\newcommand{\supp}{\operatorname{supp}}
\newcommand{\Int}{\operatorname{Int}}
\newcommand{\poinc}{\pi_1}

\newcommand{\aaa}{\alpha}
\newcommand{\bbb}{\beta}
\newcommand{\ccc}{\gamma}
\newcommand{\ddd}{\delta}
\newcommand{\epm}{A} 

\newcommand{\lieh}{\mathfrak{h}}

\newcommand{\liek}{\mathfrak{k}}

\newcommand{\prive}{\setminus}
\newcommand{\pxzx}{\pi}
\newcommand{\ptutd}{\tau}
\newcommand{\pxf}{\phi}
\newcommand{\pxef}{\psi}
\newcommand{\pcxf}{\varphi}
\newcommand{\sv}{X}
\newcommand{\uu}{U'}
\newcommand{\vv}{V}
\newcommand{\vp}{V^+}

\newcommand{\Ll}{\mathcal{L}}

\newcommand{\ovl}{\overline}
\newcommand{\udl}{\underline}
\newcommand{\ponctuation}[1]{\makebox[0pt]{\makebox[3mm][r]{#1}}}
\newcommand{\fintd}{\smash{\xto{+1}}}
\newcommand{\Split}{S}
\newcommand{\Splitbis}{S'}
\newcommand{\Ksharp}{K}
\newcommand{\lieKsharp}{\liek}
\newcommand{\Scomp}{S^c}
\newcommand{\Scompbis}{S'^{c}}
\newcommand{\Scompz}{S^{c0}}

\title{Equivariant derived category of a complete symmetric variety}
\author{St\'ephane Guillermou}

\address{Universit\'e de Grenoble I\\
D\'epartement de Math\'ematiques\\
Institut Fourier, UMR 5582 du CNRS\\
38402 Saint-Martin d'H\`eres Cedex, France}
\email{Stephane.Guillermou@ujf-grenoble.fr}

\begin{document}

\begin{abstract}
  Let $G$ be a complex algebraic semi-simple adjoint group and $\sv$ a
  smooth complete symmetric $G$-variety.  Let $L= \oplus_\alpha
  L_\alpha$ be the direct sum of all irreducible $G$-equivariant
  intersection cohomology complexes on $\sv$, and let $\E=
  \Ext^\cdot_{\D_G(X)}(L,L)$ be the extension algebra of $L$, computed
  in the $G$-equivariant derived category of $X$. We considered $\E$
  as a dg-algebra with differential $d_\E =0$, and the $\E_\alpha =
  \Ext^\cdot_{\D_G(X)}(L,L_\alpha)$ as $\E$-dg-modules.  We show that
  the bounded equivariant derived category of sheaves of $\C$-vector
  spaces on $\sv$ is equivalent to $\D_\E\langle \E_\alpha \rangle$,
  the subcategory of the derived category of $\E$-dg-modules generated
  by the $\E_\alpha$.
\end{abstract}
\maketitle

\section{Introduction}
The aim of this paper is to give a description of the equivariant
derived category of a smooth complete symmetric variety.  Let $G$ be a
complex algebraic semi-simple adjoint group, $\sigma$ an automorphism
of $G$ of order $2$ and $H=G^\sigma$.  Let $\sv$ be a complete
symmetric variety containing $G/H$, as defined in~\cite{DP85},
section~5: this is a smooth compactification of $G/H$ with a
$G$-action (extending the action on $G/H$) and with a $G$-equivariant
morphism to the canonical compactification described in~\cite{DP83}.
From~\cite{DP83} and~\cite{DP85} we have the following results on the
$G$-orbits of $X$: $X\setminus (G/H$) is the union of irreducible
smooth $G$-stable divisors with normal crossings, say $D_1,\ldots,
D_m$; any non-empty intersection $D_{i_1}\cap\ldots \cap D_{i_n}$ is
the closure of a single $G$-orbit, and, conversely, for any $G$-orbit
$\O$ in $X$, $\ovl{\O}$ is the intersection of the $D_i$ containing
$\O$.  We consider $\sv$ as a complex analytic variety, with its
transcendental topology.

We denote by $\D^b_{G}(\sv)$ the bounded equivariant derived category
of sheaves of $\C$-vector spaces on $\sv$; we let $\D^b_{G,c}(\sv)$ be
the subcategory formed by constructible objects (it is introduced
in~\cite{BL94} -- we recall the points we need in
section~\ref{eq_der_cat}).  Let $\S$ be the set of $G$-orbits of
$\sv$. For a $G$-orbit $\O$, let $\tau_\O$ be the group of components
of the stabiliser. A representation $\rho$ of $\tau_\O$ induces a
$G$-equivariant local system on $\O$ and we let $L_\O^\rho$ be the
corresponding intersection cohomology complex. Since $\ovl{\O}
\setminus \O$ consists of normal crossings divisors, $L_\O^\rho$ is in
fact a sheaf. It is known that, for any orbit $\O$, there exists $s$
with $\tau_\O \simeq (\Z/2\Z)^s$ (we recall this in
section~\ref{symmetricvarieties}).  In particular, for $\rho$
irreducible, the local system $L_\O^\rho|_{\O}$ is of rank $1$.  The
category $\D^b_{G,c}(\sv)$ is generated by the $L_\O^\rho$, for all
$\O$, $\rho$ as above.  We set $L=\bigoplus L_\O^\rho\in
\D^b_{G,c}(\sv)$, where $\O$ runs over $\S$ and $\rho$ runs over the
irreducible representations of $\tau_\O$. We consider the graded ring
$\E=\bigoplus_{i\in\N}\Ext^i(L,L)$ (here $\Ext^i(L,L)$ denotes
$\Hom_{\D^b_{G,c}(\sv)}(L,L[i])$). We view $\E$ as a differential
graded algebra (``dg-algebra''), with differential $d_\E=0$, and we
denote by $\D_\E$ the derived category of $\E$-dg-modules, introduced
in~\cite{BL94} (we recall its definition in section~\ref{eq_der_cat}).
We let $\D_{\E}\langle\E_\O^\rho \rangle$ be the subcategory generated
by the $\E$-modules $\E_\O^\rho =
\bigoplus_{i\in\Z}\Ext^i(L,L_{\O}^\rho)$, for all $\O$ in $\S$ and all
irreducible representations $\rho$ of $\tau_\O$.

On general grounds recalled below (a kind of derived version of the
Freyd-Mitchell embedding theorem) there exists a dg-algebra $R$ with
cohomology $H^\cdot(R) \simeq \E$ such that the category
$\D^b_{G,c}(\sv)$ is equivalent to a subcategory of $\D_{R}$.  We will
prove that we can actually take $R= \E$.
\begin{theorem}
  \label{thm}
  With the above notations, the categories $\D^b_{G,c}(\sv)$ and
  $\D_{\E}\langle\E_\O^\rho \rangle$ are equivalent.
\end{theorem}
The statement of this theorem is inspired by questions of Soergel
(see~\cite{So98}, \cite{So01} where it is asked whether it holds for a
Langlands parameter space, instead of a symmetric variety here).  The
case of general (possibly singular) toric varieties was done by Lunts
in~\cite{L95} and we follow the principle of his proof. Some
difficulties appear: unlike the toric case, for a $G$-orbit $\O$, the
smallest open $G$-stable set containing $\O$ is in general not
homotopically equivalent to $\O$; moreover we may have non-connected
isotropy groups.

Let us say a word about the algebra $\E$. For two $G$-orbits $\O$ and
$\O'$, $\ovl{\O} \cap \ovl{\O'}$ is a (smooth) orbit closure, say
$\ovl{\O} \cap \ovl{\O'} = \ovl{\O''}$.  Let $c$ be its complex
codimension in $\ovl{\O'}$. If $\rho$ and $\rho'$ are the trivial
representations of $\tau_\O$ and $\tau_{\O'}$, then
$\Ext^\cdot(L_{\O}^\rho,L_{\O'}^{\rho'}) \simeq H^{\cdot
  +2c}_G(\ovl{\O''})$ (in this paper for a group $G$ and a topological
space $Y$ with a $G$-action, we denote by $H^{\cdot}_G(Y) =
H^{\cdot}_G(Y;\C)$ the $G$-equivariant cohomology of $Y$, with
coefficient in $\C$). Hence, in the case where all $\tau_\O$ are
trivial, the computation of $\E$ reduces to the computation of
$H^{\cdot}_G(\ovl{\O})$, for all $G$-orbits $\O$ of $\sv$.  But
$\ovl{\O}$ is a ``regular embedding'' in the sense of definition~5
of~\cite{BDP90} (this means that $\ovl{\O}$ has a finite number of
$G$-orbits, that each $G$-orbit closure in $\ovl{\O}$ is the
transversal intersection of the codimension $1$ orbits containing it,
and that, for each $p\in \ovl{\O}$, $G_p$ has a dense orbit in the
normal space at $p$ to the $G$-orbit containing $p$).  In~\cite{BDP90}
there is a description of the equivariant cohomology of a regular
embedding, say $Y$, as a subalgebra of the product of the equivariant
cohomology algebras of the $G$-orbits of $Y$.  Note that, for a
$G$-orbit $G\cdot p$, $p\in Y$, $H^\cdot_G(G\cdot p) \simeq
H^\cdot_{G_p}(\{pt\})$ only depends on the stabiliser of $p$.  For a
symmetric variety, these stabilisers $G_p$ can be determined from the
action of the involution $\theta$ on the root system of $G$ with
respect to a suitable maximal torus. Hence we have a method to compute
the $H^{\cdot}_G(\ovl{\O})$, for all $G$-orbits $\O$ of $\sv$.  This
could lead to a combinatorial description of the algebra $\E$ and the
modules $\E_\O^\rho$, at least when the groups $\tau_\O$ are trivial.
We note that in the case of toric varieties the computation of the
equivariant intersection cohomology from combinatorial data has been
carried out in~\cite{BBFK02} and \cite{BrL03}.



\medskip

The plan of the proof mostly follows that of Lunts in the toric case
(see~\cite{L95}).  The principle is the following.  We first show that
$\D^b_{G,c}(\sv)$ is equivalent to $\D_\H \langle H_{\O}^\rho
\rangle$, where $\H$ is a sheaf of dg-algebras on a finite set $I$,
$\D_\H$ denotes the derived category of sheaves of $\H$-dg-modules,
and the $H_{\O}^\rho$ are $\H$-modules corresponding to the
$L_{\O}^\rho$ ($I$ and $\H$ are described below). On a finite set, the
category of sheaves has enough projectives, so that, for each
$H_{\O}^\rho$, we may choose a projective resolution $P_{\O}^\rho \to
H_{\O}^\rho$.  We set $P^\cdot =\bigoplus P_{\O}^\rho$, where the sum
is over all $G$-orbits $\O$ and irreducible representations $\rho$ of
$\tau_\O$.  We consider the dg-algebra $R^\cdot =
\Hom(P^\cdot,P^\cdot)$ and the left $R$-modules $R_{\O}^\rho =
\Hom(P^\cdot, P_{\O}^\rho)$.  Since $\D^b_{G,c}(\sv)$ and $\D_\H
\langle H_{\O}^\rho \rangle$ are equivalent, we have $H^\cdot(R)
\simeq \E$ and $H^\cdot(R_{\O}^\rho) \simeq \E_\O^\rho$.  We let
$\D_R$ be the derived category of left $R$-dg-modules. By general
arguments, the functor $\D_\H \langle H_{\O}^\rho \rangle \to \D_R
\langle R_{\O}^\rho \rangle$, $F \mapsto \Hom(P^\cdot,F)$ is an
equivalence of categories. Now, any quasi-isomorphism $R \to R'$
between two dg-algebras induces an equivalence of categories between
$\D_R$ and $\D_{R'}$, by restriction and extension of scalars.  Hence
the theorem is proved if we show that $R$ is quasi-isomorphic to its
cohomology (such dg-algebras are called ``formal'').

We give the details first assuming that all $\tau_\O$ are trivial.

\smallskip

1) We recall the following facts about $\D^b_{G,c}(\sv)$ (see
section~\ref{eq_der_cat}). Let $E$ be a universal bundle for $G$.  By
the construction of Bernstein-Lunts in~\cite{BL94}, the category
$\D^b_{G,c}(\sv)$ is the subcategory of $\D(E\times_G \sv)$ generated
by the sheaves induced by $G$-equivariant constructible sheaves on
$\sv$.  More precisely, in our situation it is generated by the
sheaves $E\times_G L_{\O}^\rho$, which are equal to $E\times_G
\C_{\ovl{\O}} = \C_{E\times_G \ovl{\O}}$, since we assume that the
$\tau_\O$ are trivial.  However, there is no slice theorem for a
$G$-action and we find it easier to work in the equivariant derived
category for the action of a compact group.  In fact, if $K$ is a
maximal compact subgroup of $G$, the restriction functor $\D^b_G(X)
\to \D^b_K(X)$ is fully faithful.  Since $\D^b_K(X)$ is itself a
subcategory of $\D(E\times_K \sv)$, this identifies $\D^b_{G,c}(\sv)$
with the subcategory of $\D(E\times_K \sv)$ generated by the
$\C_{E\times_K \ovl{\O}}$.

\smallskip

2) Following~\cite{L95}, we obtain a category equivalent to
$\D(E\times_G \sv)$ as follows.  From now on, we assume that the
maximal compact subgroup $K$ is compatible with $\sigma$ (i.e.
$\sigma$ commutes with the conjugation on $G$ induced by $K$).  We
decompose $\sv$ according to the $K$-orbit types, and take the
connected components: $\sv = \bigsqcup_{i\in I} \sv_i$.  This is a
stratification of $\sv$, which is precisely described in~\cite{BDP90}.
We let $\phi:\sv \to I$ be the map such that $\phi(\sv_i) = \{i\}$ and
endow $I$ with the quotient topology.  We also denote by
$\psi:E\times_K \sv \to I$ the induced map.  Since $G$ is linear, we
may take for $E$ an increasing union of $G$-manifolds (e.g. Stiefel
manifolds), $E=\bigcup_{k\in\N} E_k$. Then the sheaves of
$\Cinf^\infty$-forms of degree $i$ on $E_k\times_K \sv$,
$\Omega_{E_k\times_K \sv}^i$, form a projective system, and we define
$\Omega_{E\times_K \sv}^i$ as the projective limit of the
$\Omega_{E_k\times_K \sv}^i$.  The complex $\Omega_{E\times_K
  \sv}^\cdot$ has a natural structure of sheaf of dg-algebras and it
gives a soft resolution of $\C_{E\times_K \sv}$.

We set $\fAintro = \psi_*(\Omega_{E\times_K \sv})$. This is a sheaf of
dg-algebras on $I$. We consider the direct image functor
$\gamma:\D^b_{G,c}(\sv) \to \D_\fAintro$, $F^\cdot \mapsto
\psi_*(\Omega_{E\times_K \sv} \otimes F^\cdot)$.  We prove that
$\gamma(\C_{\ovl{\O}}) \simeq \fAintro_{\phi(\ovl{\O})}$, for any
$G$-orbit $\O$, and that $\gamma$ gives an equivalence between
$\D^b_{G,c}(\sv)$ and $\D_\fAintro \langle \fAintro_{\phi(\ovl{\O})}
\rangle$.  This point uses the following property of our
stratification. For $j\in I$, we let $V_j \subset \sv$ be the smallest
open subset of $\sv$ containing $\sv_j$ and constructible with respect
to the stratification $\sv = \bigsqcup_{i\in I} \sv_i$.  Then, there
exists a $K$-equivariant homotopy contracting $V_j$ to $K\cdot x_j$,
for some $x_j \in \sv_j$.

\smallskip

3) Let $\H$ be the cohomology sheaf of $\fAintro$, i.e. the sheaf on
$I$ associated to $U\mapsto H^\cdot (\fAintro(U))$. We consider $\H$
as a sheaf of dg-algebras with differential $0$. We prove that there
exists a sequence of quasi-isomorphisms $\fAintro \from \fAintro' \to
\fAintro'' \from\cdots \to \H$ (actually there are 5 steps in this
sequence).  This implies that the categories $\D_\fAintro \langle
\fAintro_{\phi(\ovl{\O})} \rangle$ and $\D_\H \langle
\H_{\phi(\ovl{\O})} \rangle$ are equivalent.

Let us remark that the sheaf $\fAintro$ is determined by the stalks
$\fAintro_i$, $i\in I$, since we are on a finite set. For a given
$i\in I$, the formality of the dg-algebra $\fAintro_i$ is easy: we
have $\fAintro_i = \fAintro(\phi(V_i))$, with $V_i$ as above, and
$H^\cdot (\fAintro_i)$ is the $K$-equivariant cohomology of $K\cdot
x_i$, since $V_i$ has a retraction to $K\cdot x_i$. Since
$H^\cdot_K(K\cdot x_i) = H^\cdot_{K_{x_i}}(\{pt\})$ is a polynomial
algebra (because $K_{x_i}$ is connected), any choice of
representatives for its generators gives a quasi-isomorphism
$H^\cdot_{K_{x_i}}(\{pt\}) \to \fAintro_i$.  But, of course, to obtain
the formality of the sheaf of dg-algebras $\fAintro$, we need
additionally that these quasi-isomorphisms commute with the
restriction maps $\fAintro_i \to \fAintro_j$, for $i\in \ovl{\{j\}}$.

The construction of the sequence of quasi-isomorphisms above makes use
of the description of the stabilisers given in~\cite{BDP90}. Let us
briefly recall it.  Let $\Split$ be a maximal split torus of $G$ and
denote by $x_0 \in \sv$ the class of $1_G$, $x_0 \in G/H \subset \sv$.
Then $\ovl{\Split \cdot x_0}$ is a smooth toric variety for the action
of $\Splitbis = \Split/(\Split \cap H) =\Split / \{t\in \Split; t^2
=1\}$ and contains a toric subvariety $Z$, with the following
properties. Taking intersection with $Z$ gives a bijection between the
set of $G$-orbits of $\sv$ and the set of $\Split$-orbits of $Z$.
Moreover the action of $K$ in $X$ has a fundamental domain $C_\sv
\subset Z$. We set $\Scomp = \Split \cap K$; this is a maximal compact
subgroup of $\Split$. For $x_i \in \sv_i \cap C_\sv$, we have, modulo
a finite group, the decomposition $K_{x_i} = \Scomp_i \times
\Ksharp_i$, where $\Scomp_i = \Scomp_{x_i}$ and $\Ksharp_i =
K_{x_i}\cap K^\sigma$, only depend on the stratum $\sv_i$. Hence $\H_i
= H^\cdot_{K_{x_i}}(\{pt\}) \simeq H^\cdot_{\Scomp_i}(\{pt\}) \otimes
H^\cdot_{\Ksharp_i}(\{pt\})$.

Now we build a morphism from $H^\cdot_{\Scomp_i}(\{pt\})$ to
$\fAintro_i$. Let $D_v$, $v\in V$, be the irreducible $G$-stable
divisors of $\sv$ and $\O_v$ the $G$-orbits such that $D_v =
\ovl{\O_v}$.  We denote by $\Scomp_v$ the stabiliser in $\Scomp$ of
$\O_v\cap Z$. Then $\Scomp_i = \prod_{v\in \Delta_i} \Scomp_v$, where
$\Delta_i = \{ v\in V;\: \sv_i \subset D_v\}$, and
$H^\cdot_{\Scomp_i}(\{pt\}) \simeq \C[\Xi_v;\: v\in \Delta_i]$, with
$\deg \Xi_v =2$.  Let $\delta_v$ be the $G$-equivariant fundamental
class of $D_v$ in $\sv$ and $\xi_v \in \Omega^2_{E\times_K \sv}$ a
representative of $\delta_v$.  For $i\in I$, we define
$f_i:H^\cdot_{\Scomp_i}(\{pt\}) \to \fAintro_i$, $\Xi_v \mapsto
\xi_v|_{E\times_K V_i}$.

For the factor $\smash{\Ksharp_i}$ in the decomposition, using the
Cartan model for the $K$-equivariant cohomology, we prove that we have
quasi-isomorphisms:
$$
\sect(E/\Ksharp_i; \Omega_{E/\Ksharp_i}) \xfrom{g_i}  W_i  \xto{h_i}
H^\cdot_{\Ksharp_i}(\{pt\}),
$$
where the $W_i$ are intermediate dg-algebras, which are subalgebras
of the Weil algebra of the Lie algebra of $K$.  These
quasi-isomorphisms only depend on the choice of a connection on $E$.
They are compatible with the natural maps from
$H^\cdot_{\Ksharp_i}(\{pt\})$ to $H^\cdot_{\Ksharp_j}(\{pt\})$ induced
by inclusions $\Ksharp_j \subset \Ksharp_i$ (and the similar maps
between the de Rham algebras). Finally, the maps $f_i \otimes g_i$ and
$f_i \otimes h_i$ give compatible quasi-isomorphisms between the
$\fAintro_i$ and the $\H_i$.

\smallskip

4) By steps 2 and 3, the categories $\D^b_{G,c}(\sv)$ and $\D_\H
\langle \H_{\phi(\ovl{\O})} \rangle$ are equivalent. Now we may apply
the ``principle'' explained above to the category $\D_\H$.  Indeed, we
can build explicit projective resolutions, using \v{C}ech coverings of
$I$.  Let $P_{\O} \to \H_{\phi(\ovl{\O})}$ be such a resolution,
$P^\cdot =\bigoplus P_{\O}$ and $R^\cdot = \Hom(P^\cdot,P^\cdot)$. To
conclude the proof of the theorem it is sufficient to show that
$R^\cdot$ is quasi-isomorphic to its cohomology. The sheaf $\H$ itself
has differential $0$ and this fact can be used, as in~\cite{L95}, to
endow $R^\cdot$ with a graduation different from the canonical one.
With this graduation, we prove that $H^i(R^\cdot)$ vanishes for
$i\not=0$.  Hence $R^\cdot$ is concentrated in degree $0$, and
quasi-isomorphic to its cohomology, by the natural morphisms $R \from
\tau_{\leq 0} R \to H(R)$.

\medskip

This was the idea for the case where all isotropy groups are
connected.  In general, simply taking the direct image to $I$ as above
would send some local systems to objects quasi-isomorphic to $0$ (for
example, let $H$ be a finite group, $V$ a non-trivial irreducible
representation of $H$, $L$ the $H$-equivariant local system on the
point corresponding to $V$, then we have $H^\cdot_H(\{pt\};L) =0$ --
recall that we work with coefficients in $\C$). So we modify step 2 as
follows.

\smallskip

2'.a) Let $L_E$ be the sheaf induced by $L= \bigoplus_{(\O,\rho)}
L_\O^\rho$ on $E\times_K \sv$. We will replace the dg-algebra
$\fAintro = \psi_*(\Omega_{E\times_K \sv})$, representing
$R\psi_*(\C_{E\times_K \sv})$, by an algebra representing $R\psi_*
\Rhom(L_E,L_E)$.  Let us consider a variety $Y$ endowed with two
sheaves $L$, $L'$ which are local systems on subvarieties $Z$, $Z'$
and $0$ outside (here $Y = E_k \times_K V_i$, for some $k$ and $i$,
and $L=E_k \times_K L^\rho_\O|_{V_i}$, $L'= E_k\times_K
\smash{L^{\rho'}_{\O'}|_{V_i} }$).  Here is how we represent
$\RHom(L,L')$. We assume that $Z\cap Z'$ is a smooth subvariety of
$Z'$.  We consider tubular neighbourhoods: $T_1$ of $Z\cap Z'$ in $Y$,
$T$ of $Z\setminus T_1$ in $Y\setminus T_1$, $T'$ of $Z'\setminus T_1$
in $Y\setminus T_1$.  We choose them so that $T \cap T' =\emptyset$
and $T_1 \cap \ovl{T}$ is a tubular neighbourhood of $T_1 \cap Z$ (and
the same for $T'$).  We extend $L$ to a local system $L_1$ on $T_1
\sqcup T$ and extend $L_1$ by $0$ outside $T_1 \sqcup T$; we define
$L'_1$ from $L'$ similarly.  Then we show that $\RHom(L,L') \simeq
\RHom(L_1,L'_1)$ and that the complex of sheaves $\Rhom(L_1,L'_1)$ is
actually a sheaf (concentrated in degree $0$).  Hence we may represent
$\RHom(L,L')$ by the complex $\sect(Y; \Omega_Y \otimes
\hom(L_1,L'_1))$.

We want to make this procedure work not only for two sheaves $L$, $L'$
as above, but for all pairs $L_\O^\rho$, $\smash{L_{\O'}^{\rho'}}$,
simultaneously.  For this, we prove that we can decompose $\sv$ into
``tubes'' $T_i$, such that $\sv = \bigsqcup_{i\in I} T_i$, with the
following properties.  For a $G$-orbit $\O$ and a representation
$\rho$ of $\tau_\O$, we set $Z_\O^\rho = \{x\in \sv; (L_\O^\rho)_x
\not= 0\}$ and $T_\O^\rho = \bigsqcup_{\{i; Y_i \subset Z_\O^\rho\}}
T_i$. Then $L_\O^\rho$ has an extension, $L_\O'^\rho$, to $T_\O^\rho$
and we have, for any other pair $(\O',\rho')$:
$$
\smash{ \RHom(L_\O^\rho,L_{\O'}^{\rho'}) \simeq
  \RHom(L_\O'^\rho,L_{\O'}'^{\rho'}) 
\quad \text{and}\quad
  H^i(\Rhom(L_\O'^\rho,L_{\O'}'^{\rho'}) ) =0, \quad \text{for
    $i\not=0$}.  }
$$
We set $L'_E= E\times_K \bigoplus_{(\O,\rho)} L_\O'^\rho$ and
define $\psi': E\times_K \sv \to I$ by $\psi'(E\times_K T_i) = \{i\}$.
Then $\fAintro' = \psi'_*(\Omega_{E\times_K \sv} \otimes
\hom(L'_E,L'_E))$ is a sheaf of dg-algebras on $I$, such that
$H^\cdot(\fAintro'_i) \simeq
\Ext^\cdot_{\D_K(V_i)}(L|_{V_i},L|_{V_i})$.  We also have a direct
image functor $\gamma:\D^+_K(\sv) \to \D_{\fAintro'}$, $F \mapsto
\psi'_*(\Omega_{E\times_K \sv} \otimes \hom(L'_E,J))$, for an
injective resolution $F\to J$ of $F$.  Setting $M_\O^\rho =
\gamma(L_\O^\rho)$, $\gamma$ gives an equivalence of categories
between $\D^b_{G,c}(\sv)$ and $\D_{\fAintro'} \langle M_\O^\rho
\rangle$.

\smallskip

2'.b) This procedure almost replaces step 2 above: we have built a
sheaf of dg-algebras on $I$, whose derived category is equivalent to
$\D^b_{G,c}(\sv)$. But our sheaf $\fAintro'$ is not so easy to handle,
because the ``tubes'' $T_i$ are not intrinsically related to the data.
We first built a second sheaf $\fBintro$, quasi-isomorphic to
$\fAintro'$ as follows.  Let us consider again the variety $Y$ endowed
with local systems $L$, $L'$ on subvarieties $Z$, $Z'$ extended to
$L_1$, $L_1'$ on $T_1 \sqcup T$, $T_1 \sqcup T'$. Let us assume
moreover that $Z\cap Z'$ has a neighbourhood $T_2$ containing
$\ovl{T_1}$ and such that the inclusion $Z\cap Z' \subset T_2$ is a
homotopy equivalence.  Hence we may extend $L_1$ and $L'_1$ to local
systems, $L_2$ and $L'_2$ on $T_2$. For $c= \codim^{\C}_{Z'}Z\cap Z'$,
we have a ``Gysin isomorphism'' $\Ext^\cdot(L|_{Z\cap Z'},L'|_{Z\cap
  Z'}) \isoto \Ext^{\cdot +2c}(L, L')$ given by the multiplication by
the fundamental class, $\delta$, of $Z\cap Z'$ in $Z'$. Let $\xi \in
\sect(Y;(\Omega_Y^{2c})_{T_1})$ be a representative of $\delta$. Then
the multiplication by $\xi$ gives a quasi-isomorphism $\sect(T_2;
\Omega_{T_2} \otimes \hom(L_2,L'_2)) \to \sect(Y; \Omega_Y \otimes
\hom(L_1,L'_1))[2c]$. Our sheaf $\fBintro$ is defined as follows.  For
$i\in I$, recall that the inclusion $\sv_i \subset V_i$ is a homotopy
equivalence, so that the $L_\O^\rho|_{\sv_i}$ extend as local systems,
$L_{\O,i}^\rho$, to $V_i$. We set $L''_i = E\times_K
\bigoplus_{(\O,\rho)} L_{\O,i}^\rho$ and $\fBintro_i = \sect(E\times_K
V_i;\Omega_{E\times_K V_i} \otimes \hom(L'',L''))$. The above Gysin
isomorphism implies that we have a quasi-isomorphism $\fBintro \to
\fAintro'$.  Now, as a sheaf, $\fBintro$ is intrinsically defined from
the data of $\sv$, the local systems $L_\O^\rho$ and the
stratification.  However we also have to understand what the
multiplicative structure becomes through the Gysin isomorphism: for a
third local system $L''$ on a subvariety $Z''$, we have a composition
$\Ext^\cdot(L,L') \times \Ext^\cdot(L',L'') \to \Ext^\cdot(L,L'')$.
Its counterpart on the extensions groups $\Ext^\cdot(L|_{Z\cap
  Z'},L'|_{Z\cap Z'})$,\dots is also given by the composition, but
twisted by fundamental classes of some subvarieties obtained from
intersections of $Z$, $Z'$, $Z''$.  More precisely, let us consider
$\Delta,\Delta',\Delta'' \subset V$, such that $Z = \bigcap_{v\in
  \Delta} E\times_K D_v$, $Z'$, $Z''$ being obtained in the same way
from $\Delta', \Delta''$. For $W\subset V$, we set $\delta(W) =
\prod_{v\in W} \delta_v$.  The Gysin isomorphism, from
$\Ext^\cdot(L|_{Z\cap Z'},L'|_{Z\cap Z'})$ to $\Ext^{\cdot +2c}(L, L')$
is then given by the multiplication by $\delta(\Delta \setminus
\Delta')$, and the product
$$
\Ext^\cdot(L|_{Z\cap Z'},L'|_{Z\cap Z'}) \times
\Ext^\cdot(L'|_{Z'\cap Z''},L''|_{Z'\cap Z''}) \to
\Ext^\cdot(L|_{Z\cap Z''},L''|_{Z\cap Z''})
$$
is defined by the composition of the cup-product on $Z\cap Z' \cap
Z''$ and the multiplication by $\delta((\Delta' \setminus (\Delta\cup
\Delta'') ) \cup ( (\Delta\cap \Delta'') \setminus \Delta' ) )$.  The
product in $\fBintro$ is defined by a similar formula, where the
classes $\delta_v$ are replaced by their representatives $\xi_v$
already introduced above. We prove that we can choose a chain of
quasi-isomorphisms as in step 3 compatible with the product. The final
step 4 is the same as in the case of connected isotropy groups.

\medskip

Here is the plan of the paper. In section~\ref{prelim} we recall some
facts about equivariant derived categories, Weil algebras and
constructible sheaves.  In section~\ref{loc_sys_croosdiv} we construct
the dg-algebras $\fAintro'$ and $\fBintro$ of steps 2'.a and 2'.b
above.  The main result of this section is
proposition~\ref{prop:eq_cat1}. In section~\ref{symmetricvarieties} we
recall some results of~\cite{BDP90} on symmetric varieties and use
them to prove that the hypothesis of proposition~\ref{prop:eq_cat1}
are satisfied.  Sections~\ref{formality_of_de_Rham_algebra}
and~\ref{findemo} are devoted to the proofs of steps 3 and 4.

\subsubsection*{Notations}
Notations for functors on sheaves are taken from~\cite{KS94}.  For a
topological space $X$, we denote by $\D(X)$ (resp. $\D^b(X)$) the
(resp. bounded) derived category of sheaves of $\C$-vector spaces on
$X$.  If $X$ is a real analytic manifold, we denote by
$\D^b_{\R-c}(X)$ the subcategory of $\D^b(X)$ formed by complexes with
real constructible cohomology.  The constant sheaf of group $M$ on $X$
is denoted $M_X$.  The direct and inverse images by a map $i:X\to Y$
are denoted $i_*$ and $i^{-1}$.  If $X$ and $Y$ are separated and
locally compact, $i_!$ denotes the direct image with proper supports.
If $i$ is the embedding of a locally closed subset $X$ of $Y$, and
$F\in \D(Y)$, we set $F_X = i_!i^{-1}F$ and for a group $M$, $M_X =
(M_Y)_X$.  We will also use $\sect_X(F)$, which is the subsheaf of $F$
given by the sections with support in $X$, when $X$ is closed, and the
sheaf $U\mapsto F(U\cap X)$ when $X$ is open; in general we have
$\sect_{X\cap X'} F = \sect_X \sect_{X'}F$.  The homomorphisms sheaf
is denoted $\hom(\cdot,\cdot)$.  We recall that $\rsect_X(F) \simeq
\Rhom(\C_X,F)$.  For $F\in\D^b(X)$, we set $F^* = \Rhom(F,\C_X)$.  We
will sometimes use the notation, for a subset $Z\subset Y$ and $F\in
\D(Y)$, $\rsect(Z;F) = \rsect(Z;F|_Z)$.

For a triangulated category $\D$ and objects $M_\alpha \in\D$, we
denote by $\D\langle M_\alpha \rangle$ the triangulated subcategory
generated by the $M_\alpha$, i.e. the smallest triangulated
subcategory of $\D$ containing the $M_\alpha$. For a complex
$M^\cdot$, $n\in\Z$, $M[n]$ denotes the complex $M^i[n] = M^{i+n}$
with differential $(-1)^n d_M$.

For a manifold $X$, we denote by $\Omega_X$ the de Rham complex of
$X$. If not specified, the cohomology of a space is taken with
coefficients in $\C$.

\subsubsection*{Acknowledgements}
The author would like to thank Michel Brion for many explanations and
useful comments.

\section{Preliminaries}
\label{prelim}
\subsection{Equivariant derived categories}
\label{eq_der_cat}
In this section we recall some results of~\cite{BL94} and~\cite{L95}
about equivariant derived categories and (sheaves of) dg-algebras.  We
consider, for the convenience of exposition, a linear Lie group $G$
with finitely many connected components (in
section~\ref{symmetricvarieties}, $G$ will be complex semi-simple). We
let $K$ be a maximal compact subgroup of $G$. Hence the variety $G/K$
is isomorphic to an affine space $\R^k$.

We consider a sequence of embeddings of $G$-manifolds $E_i\subset
E_{i+1}$ with free actions such that $H^k(E_i)=0$ for $0<k<i$
(since $G$ is a linear group, one may choose Stiefel varieties for the
$E_i$).  We set $E=\bigcup_{i\in \N} E_i$, endowed with the limit
topology.  The bounded equivariant derived category of $X$,
$\D^b_G(X)$, is the category formed by the triples
$F=(F_X,\ovl{F},\beta)$ where $F_X\in \D^b(X)$, $\ovl{F}\in
\D^b(E\times_G X)$ and $\beta$ is an isomorphism between the inverse
images of $F_X$ and $\ovl{F}$ on $E\times X$.  The morphisms from $F$
to $F' =(F'_X,\ovl{F'},\beta')$ are the pairs of morphisms
$(u_X,\bar{u})$, $u_X:F_X\to F'_X$, $\bar{u}:\ovl{F} \to \ovl{F'}$
commuting with $\beta$ and $\beta'$.  It is shown in~\cite{BL94} that
$\D^b_G(X)$ is independent of the choice of $E$.  (The reason to
assume that $E$ is a limit of manifolds is to be able to define
functors such as the proper direct image or the extraordinary inverse
image.) If $X$ is a real analytic manifold, we denote by
$\D^b_{G,c}(X)$ the subcategory of $\D^b_G(X)$ formed by the triples
$F$ above such that $F_X$ has real constructible cohomology sheaves.
By~\cite{BL94}, lemma 2.9.2, the forgetful functor
$F=(F_X,\ovl{F},\beta)\mapsto \ovl F$ identifies $\D^b_{G,c}(X)$ as a
full subcategory of $\D^b(E\times_G X)$.

Let $\bar q:E\times_K X \to E\times_G X$ be the quotient map.  The
restriction functor $R_{K,G}: \D^b_{G,c}(X) \to \D^b_{K,c}(X)$,
$(M_X,\ovl{M},\beta) \mapsto (M_X,{\bar{q}}^{-1} (\ovl M), \beta)$, is
fully faithful. Indeed for $F$, $F'\in \D^b_{G,c}(X)$ we have an
isomorphism $\Hom_{\D^b_{G,c}(X)}(F,F') \simeq H^0(E\times_G X;
\Rhom(\ovl{F},\ovl{F'}))$, and the corresponding isomorphism in
$\D^b_{K,c}(X)$. Since the fibres of $\bar q$ are acyclic, we have,
$\forall M\in \D^b_{G,c}(X)$, $H^0(E\times_K X; {\bar q}^{-1}(\ovl M))
\simeq H^0(E\times_G X;\ovl M)$.  This implies the claim because of
the isomorphism $R_{K,G}(\Rhom(F,F')) \simeq
\Rhom(R_{K,G}(F),R_{K,G}(F'))$.

For a $G$-orbit, $\O\simeq G/H$, of $X$, we denote by $\tau_\O =H/H^0$
the group of connected components of the isotropy group $H$.  The
$G$-equivariant sheaves with support $\O$ are in correspondence with
the representations of $\tau_\O$. (Let us recall that the objects of
$\D^b_G(X)$ concentrated in degree $0$ correspond to the
$G$-equivariant sheaves on $X$.) For a representation $\rho$ of
$\tau_\O$, let $L_\O^\rho$ denote the corresponding local system on
$\O$. Let us assume that $G$ has finitely many orbits in $X$ and that
the $\tau_\O$ are finite. Let $i_\O:\O\to X$ be the inclusion of the
orbit $\O$.  For $F\in\D^b_{G,c}(X)$, $i_\O^{-1}F$ decomposes as $
i_\O^{-1}F \simeq \oplus_i L_\O^{\rho_i}[d_i]$. We consider the
morphisms $u:F\to R(i_\O)_*i_\O^{-1}F$ and $v:\oplus_i (i_\O)_*
L_\O^{\rho_i}[d_i] \to R(i_\O)_*i_\O^{-1}F$ and denote by $F_u$ (resp.
$F_v$) the third object of a distinguished triangle built on $u$
(resp. $v$).  If $\O$ is open in $\supp F$, then $\supp F_u$ and
$\supp F_v$ contain less orbits than $\supp F$. We deduce that the
category $\D^b_{G,c}(X)$ is generated by the ${i_\O}_*L_\O^\rho$,
where $\O$ runs over the $G$-orbits and $\rho$ over the irreducible
representations of $\tau_\O$.  Hence the restriction functor
identifies $\D^b_{G,c}(X)$ with the subcategory of $\D^b_{K,c}(X)$
generated by the ${i_\O}_*L_\O^\rho$ (viewed as objects of
$\D^b_{K,c}(X)$).

As in the introduction, we let $\Omega_{E_i\times_K \sv}^d$ denote the
sheaf of $\Cinf^\infty$-forms of degree $d$ on $E_i\times_K \sv$, and
we set $\Omega_{E\times_G X}^d = \varprojlim_i \Omega_{E_i\times_G
  X}^d$.  The complex $\Omega_{E\times_K \sv}^\cdot$ is a soft
resolution of $\C_{E\times_G X}$ by a sheaf of differential graded
anti-commutative algebras.
  
\medskip

We call dg-algebra a non-negatively graded algebra over $\C$,
$\fAintro = \oplus_{i\in \N} \fAintro^i$, endowed with a differential
$d$ of degree $1$ such that, for any homogeneous elements $a,b$,
$d(ab) = (d a)b +(-1)^{\deg a}a\,d b$.  In~\cite{L95} Lunts considers
a sheaf of dg-algebras, $\fAintro$, over a topological set $I$ with
finitely many points and defines the derived category $\D_\fAintro$ as
follows.  We denote by $\M_\fAintro$ the category of sheaves of
dg-modules over $\fAintro$. A morphism $f:M\to M'$ is called a
quasi-isomorphism if $\forall i\in I$, $f_i:M_i\to M'_i$ is a
quasi-isomorphism. We consider $\K_\fAintro$, the category with the
same objects as $\M_\fAintro$ and with sets of morphisms quotiented by
the null-homotopic morphisms. Then $\D_\fAintro$ is the localisation
of $\K_\fAintro$ by the quasi-isomorphisms.  There is a substitute for
the notion of projective object in this framework: an object
$P\in\M_\fAintro$ is said $K$-projective (see~\cite{S88}
and~\cite{BL94}, p. 74) if $\Hom_{\K_\fAintro}(P,\cdot) =
\Hom_{\D_\fAintro}(P,\cdot)$.

Let $A$ be the graded algebra underlying $\fAintro$ and $\M_A$ the
category of (non differential) graded $A$-modules. For $M,N\in
\M_\fAintro$, we set $\Hom^n(M,N) = \Hom_{\M_A}(M,N[n])$ and, for
$f\in \Hom^n(M,N)$, $df = d_N\circ f -(-1)^nf\circ d_M$. This turns
$\Hom^n(M,N)$ into a complex, and we obtain a bifunctor from
$\M_\fAintro^{op}\times \M_\fAintro$ to the category of complexes of
abelian groups. Denoting by $\RHom$ its derived functor, we have
$\Hom_{\D_\fAintro}(M,N) \simeq H^0 \RHom(M,N)$.

Here is how to obtain $K$-projectives in $\M_\fAintro$. For an open
subset $U$ of $I$ and an $\fAintro$-module $F$, let $F_U$ be the
extension by $0$ of $F|_U$. For $M\in \M_\fAintro$, we have
$\Hom^\cdot(\fAintro_U,M) = \sect(U;M)$.  For a point $i\in I$, we
denote by $U_i$ the smallest open subset of $I$ containing $i$.  These
fundamental open sets generate the topology of $I$. For a sheaf $F$ on
$I$ we have $F_i = F(U_i)$. Hence the functor of sections over $U_i$
is exact and the $\fAintro$-module $\fAintro_{U_i}$ is $K$-projective.
One may deduce that the category $\K_\fAintro$ has enough
$K$-projectives and hence also enough $K$-flat objects
(see~\cite{L95}, proposition~1.7.4). Let $\phi:\fAintro\to \fBintro$
be a morphism of sheaves of dg-algebras on $I$ such that $\forall i\in
I$, $H(\phi_i): H(\fAintro_i) \to H(\fBintro_i)$ is an isomorphism.
By~\cite{L95}, proposition~1.11.2, the functors of restriction and
extension of scalars induce an equivalence of categories
$\D_{\fAintro} \simeq \D_{\fBintro}$.

\subsection{Formality of classifying spaces}
An important point in the proof of theorem~\ref{thm} is the fact that
some de Rham algebras are formal, i.e. quasi-isomorphic to their
cohomology algebras. We will use in particular the following
consequence of the results of~\cite{C51}: for a compact Lie group $K$
with universal bundle $E$ (given as above by a sequence of
$K$-manifolds) the de Rham algebras $\sect(E/H; \Omega^\cdot_{E/H})$
are formal in a compatible way for all subgroups $H\subset K$ (see
lemma~\ref{lem:formalisersousgroupes} below).

Let us first recall the definition of the Weil algebra $W(\liek)$ of a
Lie algebra $\liek$, as explained in~\cite{C51} (see
also~\cite{GS99}).  As a graded $\C$-algebra, $W(\liek) =
\Lambda^\cdot(\liek^*) \otimes S^\cdot(\liek^*)$, where
$\Lambda^\cdot(\liek^*)$ denotes the exterior algebra of $\liek^*$ and
elements of $\Lambda^1(\liek^*) \simeq \liek^*$ have degree $1$, and
$S^\cdot(\liek^*)$ denotes the symmetric algebra of $\liek^*$ and
elements of $S^1(\liek^*) \simeq \liek^*$ have degree $2$.  The
algebra $W(\liek)$ is endowed with a differential, $\delta$, of degree
$1$, and derivations, for any $x\in \liek$, $i(x)$ of degree $-1$,
$\theta(x)$ of degree $0$.  They satisfy the relations, for $x, y \in
\liek$:
\begin{gather}
  \label{eq:relCartan1}
  \theta([x,y]) = \theta(x) \, \theta(y) - \theta(y) \,\theta(x), \\
 i([x,y]) = \theta(x) \, i(y) - i(y) \, \theta(x), \\
  \label{eq:relCartan3}
\theta(x) = i(x)  \,\delta  + \delta \, i(x).
\end{gather}
They are defined as follows. First we note that, for a connected Lie
group $K$, with Lie algebra $\liek$, $\Lambda^\cdot(\liek^*)$ is
identified with the left invariant subalgebra of $\Omega^\cdot (K)$
and inherits the differential $d_\Lambda$, the contraction
$i_\Lambda(x)$ by the vector field associated to $x\in \liek$, and the
Lie derivative $\theta_\Lambda(x)$. Explicitly, for $(x,x') \in
\liek\times \liek^*$, we have $i_\Lambda(x)(x') = \langle x,x'
\rangle$, $\theta_\Lambda(x) (x') = -ad_x^t(x')$, and, for dual basis
$(x_i)$ of $\liek$, $(x'_i)$ of $\liek^*$, we have $d_\Lambda =
\frac{1}{2} \sum_i x'_i \theta_\Lambda(x_i)$.

Now we define $i(x)$, $\theta(x)$ and $\delta$ on $W(\liek)$.  Since
they are derivations, they are uniquely determined by their values on
the generators of $W(\liek)$. In the following formulas, $x\in \liek$,
$x' \in \Lambda^1(\liek^*) \simeq \liek^*$, $\tilde x' \in
S^1(\liek^*) \simeq \liek^*$ (recall that $\deg(x') = 1$, $\deg(\tilde
x') =2$). We let $h:\Lambda^1(\liek^*) \isoto S^1(\liek^*)$ be the
natural isomorphism and we consider dual basis $(x_i)$ of $\liek$,
$(x'_i)$ of $\liek^*$. With these notations we have:
\begin{alignat*}{2}
  i(x) (x'\otimes 1) &= i_\Lambda(x)(x') = \langle x,x' \rangle , &\qquad
  i(x) (1 \otimes \tilde x') &= 0, \\
  \theta(x) (x'\otimes 1)& = -ad_x^t(x') \otimes 1, &\qquad
  \theta(x) (1 \otimes \tilde x') &= -1 \otimes ad_x^t(\tilde x'),  \\
  \delta(x'\otimes 1) &= d_\Lambda(x') \otimes 1 + 1\otimes h(x'),
  &\qquad \delta (1 \otimes \tilde x') &= \textstyle \sum_i x'_i \otimes
  \theta(x_i)(\tilde x').
\end{alignat*}
By~\cite{C51}, theorem~1, we have:
\begin{equation}
  \label{eq:Wacyclique}
  H^0(W(\liek),\delta) = \C, \qquad \forall i>0 \quad
H^i(W(\liek),\delta) =0.
\end{equation}
\begin{definition}
  Let $(A,d_A)$ be a dg-algebra. One says that $\liek$ acts on $A$, if
  $A$ is endowed with two linear maps, $i,\theta: \liek \to
  \operatorname{Der}(A)$, from $\liek$ to the space of derivations of
  $A$, such that, $\forall x\in \liek$, $i(x)$ is of degree $-1$,
  $\theta(x)$ of degree $0$, $i(x)^2=0$ and $i$, $\theta$, $d_A$
  satisfy the relations~\eqref{eq:relCartan1}
  to~\eqref{eq:relCartan3}, with $d_A$ instead of $\delta$. In this
  case, the subspace of ``$\liek$-basic'' elements,
$$
A_{\liek-b} = \{ a\in A;\:\forall x\in \liek, \: i(x)(a) =
\theta(x)(a) =0\},
$$
is a sub-dg-algebra.
\end{definition}
We note that, if $\theta$ is given by differentiation of a $K$-action
in $A$, and $K$ is connected, then the subalgebra of $K$-invariants is
$A^K = \{ a\in A;\: \forall x\in \liek,\, \theta(x)(a) =0\}$ (in
general $A^K$ is not stable by $d_A$).  The main example is given by
the de Rham algebra $\Omega^\cdot (T)$ of the total space of a
$K$-principal fibre bundle, $\tau:T \to B$: for $x\in \liek$, $i(x)$
and $\theta(x)$, are the usual contraction and Lie derivative
associated to the vector field on $T$ induced by $x$.  We have
$\Omega^\cdot (T)_{\liek-b} \simeq \Omega^\cdot (B)$.  For $W(\liek)$,
the elements annihilated by all contractions $i(x)$, $x\in \liek$, are
the elements of $S(\liek^*)$; hence, if $K$ is connected,
$W(\liek)_{\liek-b} \simeq (S(\liek^*))^K$.

For a $K$-principal fibre bundle $T$ as above, recall that a connection on
$T$ is the data of projections, $\forall P\in T$, $\phi_P : T_PT \to
T_P( \tau^{-1}\tau(P))$, such that $\forall k\in K$, $\phi_{kP}$ is
conjugate to $\phi_P$ by the action of $k$ (and the $\phi_P$ vary
differentiably). Since the derivative of the $K$-action naturally
identifies $\liek$ with $T_P( \tau^{-1}\tau(P))$, a connection
corresponds to a morphism $f:\liek^* \to \Omega^1(T)$.  More
generally, for a dg-algebra $(A,d_A)$, with a $\liek$-action, a
``connection'' on $A$ is a linear morphism $f:\liek^* \to A^1$
satisfying:
\begin{equation}
  \label{eq:defconnection}
\forall x\in \liek,\: \forall x'\in \liek^*, \qquad
i(x)(f(x')) = \langle x,x' \rangle ,\quad 
\theta(x) (f(x')) = f(-ad_x^t(x')).  
\end{equation}
We extend naturally $f$ to an algebras morphism, still denoted by $f$,
from $\Lambda^\cdot(\liek^*)$ to $A$. But in general, $f$ does not
commute with the differential.  The algebra $W(\liek)$ has the
following universal property: we may extend $f$ to an algebras
 morphism,
$\bar{f} : W(\liek) \to \Omega^\cdot(T)$, with the following values on
the generators:
$$
\bar{f}(x'\otimes 1) = f(x'),\qquad
\bar{f}(1\otimes h(x')) = d_A(f(x')) - f(d_\Lambda(x')),
$$
commuting with the differentials, the ``contractions'', $i(x)$, and
the ``Lie derivatives'', $\theta(x)$.  In particular, for the
$K$-principal fibre bundle $T$ above, we obtain a morphism $\bar
f:W(\liek) \to \Omega(T)$, and it induces a morphism on the basic
sub-algebras $\tilde f: (S(\liek^*))^K \to \Omega^\cdot (B)$.

The following result can be found in~\cite{C51}, though not explicitly
stated. This is also a particular case of theorem~4.3.1 of~\cite{GS99}.
\begin{theorem}[\cite{C51}, \cite{GS99}]
  \label{thm:cohomelementbasiques}
  Let $H$ be a connected compact Lie group with Lie algebra $\lieh$,
  $A$ a dg-algebra with $\lieh$-action and a connection.  We assume
  that $H^i(A) = 0$ for $i>0$ and $H^0(A) = \C$.  Then
  $H^\cdot(A_{\lieh -b}) \simeq S(\lieh^*)^H$.
\end{theorem}
We return to the situation of a compact connected Lie group $K$,
acting on a universal bundle $E$ which is an increasing union of
$K$-manifolds, $E = \bigcup_i E_i$. We choose compatible connections
on the $E_i$ (i.e. the connection on $E_i$ is the restriction of the
one on $E_{i+1}$). This gives a connection, in the algebraic sense
of~\eqref{eq:defconnection}, $f:\liek^* \to \sect(E;\Omega^1_E)$.  It
induces a dg-algebras morphism $\bar{f}:W(\liek) \to
\sect(E;\Omega^\cdot_E)$, compatible with the contraction $i$ and the
Lie derivative $\theta$.  For a connected subgroup $H \subset K$, with
Lie algebra $\lieh$, the action of $\liek$ on $W(\liek)$ obviously
restricts to an action of $\lieh$.
\begin{lemma}
  \label{lem:formalisersousgroupes}
  With the notations $H$, $K$, $E$, $f$, introduced above, the space
  of $\lieh$-basic elements of $W(\liek)$ is $W(\liek)_{\lieh-b}
  \simeq ( \Lambda^\cdot(\lieh^\bot) \otimes S^\cdot(\liek^*))^H$,
  where $\lieh^\bot \subset \liek^*$ denotes the orthogonal of
  $\lieh$.  The projection $W(\liek)_{\lieh-b} \to S^\cdot(\lieh^*)^H$
  is a quasi-isomorphism. The morphism induced by the connection,
  $W(\liek)_{\lieh-b} \to \sect(E;\Omega^\cdot_E)_{\lieh-b} \simeq
  \sect(E/H; \Omega^\cdot_{E/H})$, also is a quasi-isomorphism.
  
  The normaliser of $H$, $N_K(H)$, acts on $W(\liek)_{\lieh-b}$,
  $S^\cdot(\liek^*)^H$ and $\sect(E/H; \Omega^\cdot_{E/H})$, and the
  above morphisms are $N_K(H)$-equivariant. For another connected subgroup
  $H_1 \subset H \subset K$, with Lie algebra $\lieh_1$, we have a
  commutative diagram
$$
\xymatrix@R=5mm{
 \sect(E/H; \Omega^\cdot_{E/H})  \ar[d] 
& W(\liek)_{\lieh-b}  \ar[d] \ar[r] \ar[l]
& S^\cdot(\lieh^*)^H\ar[d] \\
\sect(E/H_1; \Omega^\cdot_{E/H_1})   & W(\liek)_{\lieh_1-b} \ar[r] \ar[l]
&  S^\cdot(\lieh_1^*)^{H_1}  \ponctuation{,}
}
$$
where the horizontal arrows are quasi-isomorphisms.
\end{lemma}
\begin{proof}
  The $\lieh$-basic elements of $W(k)$ are the elements annihilated by
  all $i(x)$ and $\theta(x)$ for $x\in \lieh$. Since $i(x)$ is a
  derivation and acts trivially on $S^\cdot(\liek^*)$, the set of
  elements of $W(k)$ annihilated by all $i(x)$, $x\in \lieh$, is
  $\Lambda^\cdot(\lieh^\bot) \otimes S^\cdot(\liek^*)$.  Since $H$ is
  connected, the elements annihilated by the $\theta(x)$ are the
  $H$-invariants. Hence we have the description of
  $W(\liek)_{\lieh-b}$ given in the lemma. By this description,
  $W(\liek)_{\lieh-b}$ admits a projection to $S^\cdot(\liek^*)^H$ and
  hence to $S^\cdot(\lieh^*)^H$.
  
  Let us choose an $H$-stable decomposition $\liek = \lieh \oplus V$.
  It induces an $H$-equivariant splitting $g:\lieh^* \to \liek^*
  \simeq W^1(\liek)$.  This is a connection on $W(\liek)$, for the
  $\lieh$-action, in the sense of~\eqref{eq:defconnection}.  Hence it
  gives a morphism of dg-algebras $\bar{g}:W(\lieh) \to W(\liek)$.
  By~\eqref{eq:Wacyclique} and theorem~\ref{thm:cohomelementbasiques},
  the induced morphism $\tilde{g}:S^\cdot(\lieh^*)^H \simeq
  W(\lieh)_{\lieh-b} \to W(\liek)_{\lieh-b}$ is a quasi-isomorphism.
  We note that, by definition, $\tilde{g}$ also is a splitting of the
  projection $q:W(\liek)_{\lieh-b} \to S^\cdot(\lieh^*)^H$, so that
  $q$ is a quasi-isomorphism too.
  
  The composition $f_1 = f \circ g :\lieh^* \to \sect(E;\Omega^1_E)$
  also is a connection on $\sect(E;\Omega^\cdot_E)$, for the
  $\lieh$-action.  Hence it gives a morphism $\bar{f}_1: W(\lieh) \to
  \sect(E;\Omega^\cdot_E)$, and we have $\bar{f}_1 = \bar{f} \circ
  \bar{g}$. By theorem~\ref{thm:cohomelementbasiques} again, the
  induced morphism on the $\lieh$-basic elements, $(\bar{f}_1)_{\lieh
    -b}$ is a quasi-isomorphism.  Since $(\bar{g})_{\lieh -b}$ also
  is, $(\bar{f})_{\lieh -b} : W(\liek)_{\lieh-b} \to
  \sect(E;\Omega^\cdot_E)_{\lieh-b}$ is a quasi-isomorphism, as
  claimed.
  
  The compatibility of the above morphisms with the $N_K(H)$-action
  and the commutativity of the diagram follows from the functoriality
  of the construction.
\end{proof}

\subsection{Constructible sheaves}
\label{subsection:constrsheaves}
Here we recall some results of~\cite{KS94} on constructible (complex
of) sheaves on real analytic manifolds.  Let $Y$ be a real analytic
manifold. We say that a locally finite partition of $Y$ by locally
closed real analytic manifolds, $Y= \bigsqcup_{i\in I} Y_i$, is a
stratification if $\forall i,j \in I$, $Y_i \cap \ovl{Y_j} \not=
\emptyset$ implies $Y_i \subset \ovl{Y_j}$.  For two closed subsets
$A$, $B$ of $T^*Y$, which are conic, i.e. stable by the action of
$\R_{>0}$ in the fibres, we let $A \widehat + B$ be the subset of
$T^*Y$ defined as follows (see~\cite{KS94} definition 6.2.3 and remark
6.2.8): in a local chart $U \simeq \R^d$ of $Y$, $(x,\xi) \in T^*U
\simeq \R^d \times \R^d$ belongs to $A \widehat + B$ iff there exists
sequences $(x_n,\xi_n) \in A$, $(y_n,\eta_n) \in B$ such that $x_n \to
x$, $y_n \to x$, $\xi_n + \eta_n \to \xi$ and $|x_n - y_n| |\xi_n| \to
0$. We let $\pi_Y:T^*Y \to Y$ be the projection and set $T^*_{Y_i}Y =
\{(y,\xi) \in T^*Y;\: y\in Y_i, \; \langle \xi, T_yY_i \rangle =0\}$.
We say that the stratification is a $\mu$-stratification if $\forall
i\not= j \in I$ such that $Y_i \subset \ovl{Y_j}$ we have $(T^*_{Y_j}Y
\widehat + T^*_{Y_i}Y) \cap \pi^{-1}_Y(Y_i) \subset T^*_{Y_i}Y$.  We
note that if $Y= \bigsqcup_{i\in I} Y_i$ is a $\mu$-stratification
then so is the product $Y\times \R^d = \bigsqcup_{i\in I} Y_i \times
\R^d$, and if $S=\bigsqcup_{i\in I} S_i$ is a $\mu$-stratification of
the $d$-sphere then so is the cone over it: $\R^{d+1} = \{0\} \sqcup
(\bigsqcup_{i\in I}\R_{>0} \cdot S_i)$ (the condition is trivial at
the vertex and at other strata the stratification is diffeomorphic to
a product).  Finally, if a submanifold $Z$ of $Y = \bigsqcup_{i\in I}
Y_i$ intersects all strata transversally and the partition is a
$\mu$-stratification, then so is the partition $Z = \bigsqcup_{i\in I}
Z\cap Y_i$.

For a complex of sheaves $F \in \D^b(Y)$, we have the notion of
micro-support, $SS(F)$, which is a closed conic subset of $T^*Y$.  We
refer to definition 5.1.1 of~\cite{KS94} and just recall that, if $Y =
\bigsqcup_{i\in I} Y_i$ is a $\mu$-stratification, and $F$ is
constructible with respect to this stratification, then $SS(F) \subset
\bigsqcup_i T^*_{Y_i}Y$ (see proposition 8.4.1 of~\cite{KS94}). We
denote by $\D^b_{\R-c}(Y)$ the subcategory of $\D^b(Y)$ formed by
complexes with real constructible cohomology (with respect to any
stratification). We will use several times the following results
of~\cite{KS94}.
\begin{lemma}[\cite{KS94}, lemma 5.4.14]
  \label{lemmicrosupport}
  Let $Y$ be a real analytic manifold, $F \in \D^b_{\R-c}(Y)$, $G\in
  \D^b(Y)$ and assume that $SS(F) \cap SS(G) \subset T^*_YY$. Then the
  natural morphism $R\hom(F,\C_Y) \otimes G \to R\hom(F,G)$ is an
  isomorphism.
\end{lemma}
\begin{lemma}[\cite{KS94}, lemma 8.4.7]
  \label{lemKS}
  Let $Y$ be a real analytic manifold, $F\in\D^b_{\R-c}(Y)$, $f:Y\to
  \R$ a real analytic function such that $f|_{\supp F}$ is proper. For
  $\varepsilon >0$ we set $Z= f^{-1}(0)$, $Z_\varepsilon =
  f^{-1}([0,\varepsilon])$, $U_\varepsilon = f^{-1}([0,\varepsilon[)$.
  Then there exists $\varepsilon_0>0$ such that, $\forall
  \varepsilon$, $0 < \varepsilon <\varepsilon_0$, we have the
  isomorphisms
  $$
  H^\cdot_{Z}(Y;F) \isoto H^\cdot_{Z_\varepsilon}(Y;F), \qquad
  H^\cdot(Z_\varepsilon;F) \isoto H^\cdot(U_\varepsilon;F) \isoto
  H^\cdot(Z;F).
  $$
\end{lemma}

\subsection{Local systems outside normal crossings divisors}
\label{normal crossings divisors}
We make here some easy remarks on local systems defined outside normal
crossings divisors.  Let $Y$ be a smooth complex manifold and
$(D_v)_{v\in V}$ a finite family of smooth normal crossings divisors.
We set $U = Y \setminus \bigcup_{v\in V} D_v$.  Local systems (over
$\C$) on $U$ are in bijective correspondence with complex
representations of $\poinc(U)$.  For such a representation, $\rho$, we
denote by $L^\rho$ the associated local system.

We fix $v\in V$ and set $Y_v = Y \setminus \bigcup_{w\not= v} D_w$,
$D'_v = D_v \cap Y_v$. We let $T_v$ be a tubular neighbourhood of
$D'_v$ in $Y_v$, homeomorphic to the normal bundle of $D'_v$ in $Y_v$,
and with a projection $\pi_v :T_v \to D'_v$.  For $x\in T_v \cap U =
T_v \setminus D'_v$, the fibre $\pi_v^{-1} \pi_v(x) \simeq \C$ is
oriented, and we let $\gamma_x$ be a loop in $\pi_v^{-1} \pi_v(x)$
with base-point $x$ and turning $+1$ time around $0$.  Now, let $b\in
U$ be a base-point, $\tau$ a path from $b$ to $x$. The conjugacy class
of the image of $\tau^{-1} \, \gamma_x \, \tau$ in $\poinc(U)$ is
well-defined. We denote it by $C_x$.  If $x'$ is another point of $T_v
\cap U$, and $\gamma:[0,1] \to T_v \cap U$ a path from $x$ to $x'$,
the loops $\gamma_x$ and $\gamma^{-1}\, \gamma_{x'} \, \gamma$ are
homotopic. It follows that $C_x$ is independent of $x$. Hence the
image of $C_x$ by $\rho$ also is well-defined up to conjugacy. We call
it the monodromy of $L^\rho$ around $D_v$.  We quote the following
facts for later use.
\begin{lemma}
  \label{lem:monodromy}
  In the above situation, let $L$ be a local system of finite rank on
  $U$.
  \begin{itemize}
  \item [(i)] If the monodromy of $L$ around $D_v$ is $Id$, then $L$
    extends as a local system, $L'$, to $Y_v$.  For $w \not= v$, the
    monodromy of $L'$ around $D_w$ is the monodromy of $L$ around
    $D_w$.
  \item[(ii)] Let $j:U\to Y$ be the inclusion.  We assume that $\rho$
    factors through a finite quotient of $\poinc(U)$.  If, for each
    $v\in V$, the monodromy of $L$ around $D_v$ has no eigenvalue
    equal to $1$, then $Rj_* L \simeq j_* L \simeq j_! L$.
  \end{itemize}
\end{lemma}
\begin{proof}
  (i) Let $j_v:U \to Y_v$ be the inclusion. It follows from the
  definition that $L' = (j_v)_* L$ has the required properties.
  
  (ii) The assertion is equivalent to $(Rj_* L)_x=0$ for any $x\in
  Y\setminus U$. Since this is a local problem around $x$, we may
  assume that $Y = \C^n$, and we have coordinates $(x_1,\ldots,x_n)$
  such that $x=(0,\ldots,0)$, $D_v=\{x_v=0\}$, $v=1,\ldots,m$,
  $U=X\setminus \bigcup_{v=1,\ldots,m} D_v$. Then $(R^ij_* L)_x=
  \varinjlim_V H^i(V\cap U;L)$, where $V$ runs over neighbourhoods of
  $0$.  We may assume $V$ of the type $V = \{ (\udl{x}); \:\forall i,
  \;|x_i| < \varepsilon \}$. Then $V\cap U$ decomposes as a product
  $V\cap U \simeq \R^{2n - m} \times (S^1)^m$ and $\poinc(V\cap U)$
  acts on $L|_{V\cap U}$ through a finite abelian group.  Hence we may
  decompose $L|_{V\cap U}$ into a sum of irreducible components,
  $L_k$, which are local systems of rank $1$.  Then $L_k \simeq
  \C_{\R^{2n-m}} \boxtimes L_k^1 \boxtimes \cdots \boxtimes L_k^m$,
  for rank $1$ local systems $L_k^j$ on $S^1$.  The monodromy of $L_k$
  around $D_j$ is the monodromy of $L_k^j$ around $S^1$. By
  hypothesis, it is not $1$, so that $L_k^j$ is non trivial and we
  have $H^0(S^1;L_k^j) = H^1(S^1;L_k^j) =0$. The K\"unneth formula
  yields $\forall i$, $H^i(V\cap U;L) =0$, as desired.
\end{proof}

\section{Categories of sheaves and dg-algebras}
\label{loc_sys_croosdiv}
We consider a manifold $Y$ endowed with a finite stratification
$Y=\bigsqcup_{i\in I} Y_i$ by locally closed submanifolds. We denote
by $\pxf:Y\to I$ the natural map and endow $I$ with the quotient
topology. We consider sheaves, $(L_\aaa)_{\aaa\in \epm}$ on $Y$,
constructible with respect to this stratification and which are local
systems of finite rank on $Z_\aaa = \{x\in Y;\: (L_\aaa)_x \not=0 \}$.
We will realize $\D(Y)\langle L_\aaa \rangle$ as a derived category of
dg-modules over a sheaf of dg-algebras, $\fA$, on the finite set $I$
(see proposition~\ref{prop:eq_cat1} below).  This sheaf $\fA$ will be
quasi-isomorphic to $R\pxf_* \Rhom(\oplus L_\aaa, \oplus L_\aaa)$. We
make the following hypothesis on the stratification and the $L_\aaa$.
\begin{assumptions}
\label{assumptionstratification}
Let $Y$ be a complex manifold endowed with a finite
$\mu$-strati\-fi\-ca\-tion, $Y=\bigsqcup_{i \in I} Y_i$, by real
analytic submanifolds. We assume that $Y$ is an analytic open subset
of an analytic manifold, $X$, such that $\ovl{Y}$ is compact and has a
stratification, $\ovl{Y}=\bigsqcup_{i \in I} Y'_i $, satisfying:
$\forall i\in I$, $Y_i =Y \cap Y'_i$ (note that $\ovl{Y}$ has no
additional stratum).  For $i\in I$, we define, as in
section~\ref{eq_der_cat}, $U_i$ to be the smallest open subset of $I$
containing $i$:
\begin{equation}
  \label{eq:defUi}
U_i = \{j\in I;\: Y_i \subset \ovl{Y_j}\}.  
\end{equation}
We consider a finite family of (complex) smooth, connected, normal
crossings divisors, $(D_v)_{v\in V}$, on $Y$.  We assume that the
divisors are union of strata: $D_v = \bigsqcup_{i\in I_v} Y_i$, for
some $I_v \subset I$. We define:
\begin{equation}
  \label{eq:defZdelta}
 \textstyle  
\text{for $\Delta\subset V$},\quad
Z_\Delta = \bigcap_{v\in \Delta} D_v, \qquad \quad
\S = \{ \Delta\subset V;\: Z_\Delta \not= \emptyset \}.
\end{equation}
We also consider a finite family of constructible sheaves,
$(L_\aaa)_{\aaa\in \epm}$ on $Y$, and set $Z_\aaa = \{x\in Y;\:
(L_\aaa)_x \not=0 \}$.  We make the following hypothesis on these
data:
\begin{itemize}
\item [(i)] $\forall i,i'\in I$, $\exists j\in I$, $ U_i\cap U_{i'}
  =U_j$.
\item[(ii)] $\forall i\in I$ there exists a homotopy $h: [0,1] \times
  \pxf^{-1}(U_i) \to \pxf^{-1}(U_i)$ contracting $\pxf^{-1}(U_i)$ to a
  submanifold $Y'_i$ of $Y_i$ and preserving the closures of strata:
  $\forall j\in U_i$, $h([0,1] \times Y_j) \subset \ovl{Y_j}$.
\item[(iii)] $\forall i \in I$, $\forall v\in V$, $\exists j\in I$,
  $U_i \setminus \pxf(D_v) = U_j$.
\item[(iv)] $\forall \aaa \in \epm$, $\exists \Delta_\aaa \in \S$,
  $\exists \Delta'_\aaa \subset (V \setminus \Delta_\aaa)$ such that
  $Z_\aaa = Z_{\Delta_\aaa} \setminus \bigcup_{v\in \Delta'_\aaa}
  D_v$.
\item[(v)] $\forall \aaa \in \epm$, $L_\aaa|_{Z_\aaa}$ is a local
  system on $Z_\aaa$ with monodromy $-Id$ around each $Z_{\Delta_\aaa}
  \cap D_v$, for $v\in \Delta'_\aaa$ (see section~\ref{normal crossings
    divisors}).
\end{itemize}
\end{assumptions}
\begin{example}
  \label{exempletorique}
  We will verify in section~\ref{symmetricvarieties} that the
  decomposition of a symmetric variety given in~\cite{BDP90} satisfies
  these assumptions. A more simple example is given by smooth toric
  varieties: let $T = (\C^*)^l$ be a torus, $D\subset T$ a subgroup
  consisting of order $2$ elements and $T' = T /D$.  Let $Y$ be a
  smooth $T'$-toric variety, with the action of $T$ through $T'$. We
  let $(Y_i)_{i\in I}$ be the stratification given by the $T'$-orbits,
  $(D_v)_{v\in V}$ be the set of $T'$-stable irreducible divisors.
  For a $T'$-orbit $\O$ and $x\in \O$, we set $\tau_\O = T_x/ T_x^0$;
  we have $\tau_\O \simeq (\Z/2\Z)^{c_{\O}}$, for some $c_{\O} \in
  \N$. The irreducible $T$-equivariant local systems on $\O$
  correspond to irreducible representations of $\tau_\O$.  Let $\epm$
  be the set of pairs $\aaa=(\O,\rho)$, where $\O$ is an orbit and
  $\rho$ an irreducible representation of $\tau_\O$. We let
  $\Delta_\aaa\in \S$ be such that $\ovl{\O} = \bigcap_{v\in
    \Delta_\aaa} D_v $ and let $L'_\aaa$ be the local system on $\O$
  given by $\rho$.  Since $\tau_\O$ is a $2$-group, the irreducible
  representations are one dimensional and the elements of $\tau_\O$
  act by $1$ or $-1$.  In particular $L'_\aaa$ has monodromy $Id$ or
  $-Id$ around any divisor $D_v \cap \ovl{\O}$ (for $v$ such that $\O
  \not\subset D_v$ and $\ovl{\O} \cap D_v \not= \emptyset$). We let
  $\Delta'_\aaa$ be the set of $v\in V$ for which this monodromy is
  $-Id$.  Then $L'_\aaa$ extend as a local system to $Z_{\Delta_\aaa}
  \setminus \bigcup_{v\in \Delta'_\aaa} D_v$ and we let $L_\aaa$ be
  the extension by $0$ of this local system. Then the assumptions
  above are satisfied in this situation.
\end{example}
\begin{remarks}
  \label{rem:apresassumptions}
  1) In fact we do not use the complex structure; only the geometry of
  the intersections of the $D_v$ matters. In particular, the strata
  $Y_i$ are not assumed to be complex.
  
  2) In view of lemma~\ref{lem:monodromy}, hypothesis (iv) and (v)
  have the following consequences: let $j_\aaa :Z_\aaa \to Y$ be the
  inclusion.  Then $L_\aaa \simeq R(j_\aaa)_* (L_\aaa|_{Z_\aaa})
  \simeq (j_\aaa)_* (L_\aaa|_{Z_\aaa}) \simeq (j_\aaa)_!
  (L_\aaa|_{Z_\aaa})$ (or, with different notations, $L_\aaa \simeq
  \rsect_{Z_\aaa}(L_\aaa) \simeq (L_\aaa)_{Z_\aaa}$). For $\aaa, \bbb
  \in \epm$, we will give representatives for the complex
  $\Rhom(L_\aaa,L_\bbb)$. We already note that
  $$
  \Rhom(L_\aaa,L_\bbb) \simeq \Rhom((L_\aaa)_{Z_\aaa},
  \rsect_{Z_\bbb}(L_\bbb)) \simeq \Rhom(L_\aaa,
  \rsect_{Z_\aaa \cap Z_\bbb}(L_\bbb)).
  $$
  In particular, if $Z_\aaa \cap Z_\bbb = \emptyset$, then
  $\Rhom(L_\aaa,L_\bbb) =0$.
  
  3) We note that $Y_i$ is closed in $\pxf^{-1}(U_i)$ (because the
  strata contained in $\ovl{Y_i} \setminus Y_i$ cannot be in
  $\pxf^{-1}(U_i)$) so that $U_i\setminus \{i\}$ is open in $I$. In
  particular, if $U_i = U_j$ then $i=j$. Hence the $j$ in hypothesis
  (i) and (iii) are unique.
  
  4) For any $i\in I$ and any closed subset $J$ of $U_i$, the homotopy
  $h$ of (ii), also contracts $\pxf^{-1}(J)$ to $Y'_i$.  Hence the
  inclusions $Y'_i \subset \pxf^{-1}(J) \subset \pxf^{-1}(U_i)$ are
  homotopy equivalences (i.e. induce isomorphisms on all homotopy
  groups). In particular the inclusion $Y_i \subset \pxf^{-1}(U_i)$ is
  a homotopy equivalence. Hence, for any $\aaa\in \epm$ and $i\in I$
  with $Y_i \subset Z_\aaa$, the local system $L_\aaa|_{Y_i}$ has a
  unique extension to a local system defined on $\pxf^{-1}(U_i)$.  We
  denote this extension by $L_{\aaa,i}$.  We have
  $L_{\aaa,i}|_{\pxf^{-1}(U_i)\cap Z_\aaa} \simeq
  L_\aaa|_{\pxf^{-1}(U_i)\cap Z_\aaa}$ and, for $j$ such that $Y_j
  \subset Z_\aaa$ and $Y_i \subset \ovl{Y_j}$, we have an isomorphism
  $u_{ij}: L_{\alpha,i}|_{\pxf^{-1}(U_j)} \simeq L_{\alpha,j}$.  Since
  $u_{ij}$ is determined by its restriction to $Y_j$, the $u_{ij}$,
  $(i,j) \in I^2$, satisfy the same relations as their restrictions to
  $Z_\aaa$.  In particular they satisfy the cocycle condition which
  says that the $L_{\aaa,i}$ glue together into a local system, say
  $L^1_\aaa$, on $\bigcup_{\{i; Y_i \subset Z_\aaa\}} U_i$.
    
  5) Since the $U_i$ form a basis of the topology of $I$, a sheaf $F$
  on $I$ is determined by its stalks $F_i = F(U_i)$, for all $i\in I$,
  and the restriction maps, $F_i \to F_j$, for all $i,j \in I$ with $i
  \in \ovl{\{j\}}$.  Conversely, the data of groups $F_i$, for all
  $i\in I$, and restriction maps, $f_{ji}:F_i \to F_j$, for all $i,j
  \in I$ with $i \in \ovl{\{j\}}$, satisfying $f_{ij} \circ f_{jk} =
  f_{ik}$ (whenever it makes sense) define a sheaf on $I$.

\end{remarks}
\begin{notations}
  \label{not:IetIprime}
  We introduce the following notations, for $\aaa, \bbb\in \epm$:
\begin{gather*}
  Z_{\aaa\bbb} = Z_\aaa \cap Z_\bbb , \qquad \qquad d_{\aaa\bbb} =
  \codim^\C_{Z_\bbb} Z_{\aaa\bbb}, \qquad \qquad
  I_{\alpha\beta} = \pxf(Z_{\aaa\bbb}), \\
  \Delta'_{\aaa\bbb} = (\Delta'_\aaa \setminus \Delta'_\bbb) \cup
  (\Delta'_\bbb \setminus \Delta'_\aaa), \qquad
  \textstyle I'_{\alpha\beta} = \pxf \bigl( Z_{\Delta_\aaa \cup
    \Delta_\bbb} \setminus (\bigcup_{v\in \Delta'_{\aaa\bbb}} D_v)
  \bigr)   \setminus I_{\aaa\bbb},
\end{gather*}
and, for $i \in I_{\aaa\bbb}$, $L_{\aaa,i}$, $L_{\bbb,i}$ as in
remark~\ref{rem:apresassumptions}~(4), we introduce the following
sheaf on $\pxf^{-1}(U_i)$: $\Omega_{\aaa\bbb, i} =
\Omega_{\pxf^{-1}(U_i)} \otimes \hom(L_{\aaa,i},L_{\bbb,i})$.
\end{notations}
Let us prove the following facts:
\begin{itemize}
\item[(a)] if $I_{\aaa\bbb} = \emptyset$ then $I'_{\aaa\bbb} =
  \emptyset$; in any case $I'_{\aaa\bbb} \subset
  \ovl{I_{\aaa\bbb}}$.
\item[(b)] $I_{\aaa\bbb} \sqcup I'_{\aaa\bbb}$ is open in
  $\ovl{I_{\aaa\bbb}}$.
\item[(c)] $\forall i \in I'_{\aaa\bbb}$, $\exists ! j \in
  I_{\aaa\bbb}$, such that $U_i \setminus \pxf (\bigcup_{v\in
    \Delta'_\aaa \cap \Delta'_\bbb} D_v) = U_j$.
\end{itemize}
For (a) we note that $I_{\aaa\bbb} = \emptyset$ means $ (\bigcap_{v\in
  \Delta_\aaa} D_v \setminus \bigcup _{w\in \Delta'_\aaa} D_w) \cap
(\bigcap_{v\in \Delta_\bbb} D_v \setminus \bigcup _{w\in \Delta'_\bbb}
D_w) = \emptyset$, which is equivalent to $\bigcap_{v\in \Delta_\aaa
  \cup\Delta_\bbb} D_v \subset \bigcup _{w\in \Delta'_\aaa
  \cup\Delta'_\bbb} D_w$, or also to $(\Delta_\aaa \cup\Delta_\bbb)
\cap (\Delta'_\aaa \cup\Delta'_\bbb) \not= \emptyset$.  Since
$\Delta_\aaa \cap\Delta'_\aaa = \emptyset$ and $\Delta_\bbb
\cap\Delta'_\bbb = \emptyset$, this implies that $(\Delta_\aaa
\cup\Delta_\bbb) \cap \Delta'_{\aaa\bbb} \not= \emptyset$, and then
$Z_{\Delta_\aaa \cup \Delta_\bbb} \setminus (\bigcup_{v\in
  \Delta'_{\aaa\bbb}} D_v) = \emptyset$. In particular, $I'_{\aaa\bbb}
= \emptyset$ as claimed, and $I'_{\aaa\bbb} \subset
\ovl{I_{\aaa\bbb}}$.  Now we note that $Z_{\aaa\bbb}$ is open in
$Z_{\Delta_\aaa \cup \Delta_\bbb}$. Arguing locally around each
connected component of $Z_{\Delta_\aaa \cup \Delta_\bbb}$, we deduce
that, if $I_{\aaa\bbb} \not= \emptyset$, we also have
$\ovl{I_{\aaa\bbb} } = \pxf( Z_{\Delta_\aaa \cup \Delta_\bbb} )
\supset I'_{\aaa\bbb}$.

For (b), we have $\ovl{I_{\aaa\bbb}} \setminus (I_{\aaa\bbb} \sqcup
I'_{\aaa\bbb}) = \pxf ( Z_{\Delta_\aaa \cup \Delta_\bbb} \cap
(\bigcup_{v\in \Delta'_{\aaa\bbb}} D_v) )$ and this is closed.
  
For (c), applying hypothesis (iii) of
assumptions~\ref{assumptionstratification} several times, we know that
there exists a unique $j\in I$ such that $U_i \setminus \pxf
(\bigcup_{v\in \Delta'_\aaa \cap \Delta'_\bbb} D_v) = U_j$.  We note
that $U_i \cap \ovl{I_{\aaa\bbb} } \not= \emptyset$ and
$\ovl{I_{\aaa\bbb} } \not\subset \pxf (\bigcup_{v\in \Delta'_\aaa \cap
  \Delta'_\bbb} D_v)$ (remark that $I_{\aaa\bbb} \not= \emptyset$),
hence $U_j \cap \ovl{I_{\aaa\bbb} } \not= \emptyset$. This implies $j
\in \ovl{I_{\aaa\bbb} }$.  But $i\not\in \ovl{I_{\aaa\bbb}} \setminus
(I_{\aaa\bbb} \sqcup I'_{\aaa\bbb})$ which is closed, hence $j\not\in
\ovl{I_{\aaa\bbb}} \setminus (I_{\aaa\bbb} \sqcup I'_{\aaa\bbb})$ as
well.  We thus obtain $j \in I_{\aaa\bbb}\sqcup I'_{\aaa\bbb}$.  Since
$j \not\in \pxf (\bigcup_{v\in \Delta'_\aaa \cap \Delta'_\bbb} D_v)$,
this implies $j \in I_{\aaa\bbb}$.

\smallskip

\begin{definition}
\label{env:defin_fB}
For $i \in I'_{\aaa\bbb}$, we denote by $i(\aaa,\bbb)$ the element $j
\in I_{\aaa\bbb}$ given by assertion (c) above.  We define a sheaf
$\fA^{\aaa\bbb}$ on $I$ by its stalks at $i\in I$:
\begin{equation}
  \label{eq:def_fB}
\fA^{\aaa\bbb}_i =
\begin{cases}
  \sect(\pxf^{-1}(U_i);\Omega_{\aaa\bbb, i})) \,[-2d_{\aaa\bbb}] &
  \text{if
    $i\in I_{\aaa\bbb}$},  \\
  \fA^{\aaa\bbb}_j & \text{if $i \in I'_{\aaa\bbb}$
    and $j =i(\aaa,\bbb)$}, \\
  0 & \text{if $i\not\in I_{\aaa\bbb} \sqcup I'_{\aaa\bbb}$},
\end{cases}
\end{equation}
and the natural restriction maps.  Let us check that this is indeed a
sheaf. The first case ($i\in I_{\aaa\bbb}$) defines a sheaf, say $\fA'
= \pxf_*(\Omega_Y \otimes \hom(L^1_{\aaa},L^1_{\bbb})) $ (with
$L^1_{\aaa}$, $L^1_{\bbb}$ as in remark~\ref{rem:apresassumptions}
(4)), on $I_{\aaa\bbb}$. Then the second case defines a sheaf, say
$\fA''$, on $I_{\aaa\bbb} \sqcup I'_{\aaa\bbb}$, as $u_*\fA''$, where
$u$ is the inclusion $u: I_{\aaa\bbb} \to I_{\aaa\bbb} \sqcup
I'_{\aaa\bbb}$.  Finally the third case defines $\fA^{\aaa\bbb}$ as
the extension by $0$ of $\fA''$.

We also define $\fA = \oplus_{\aaa,\bbb\in \epm} \fA^{\aaa\bbb}$.
\end{definition}
\begin{remarks}
  1) The stalks $\fA^{\aaa\bbb}_i$ are defined to be $0$ when the
  stratum $Y_i$ is included in a divisor $D_w$ such that the local
  system $\hom(L_\aaa,L_\bbb)$ (on $Z_{\aaa\bbb}$) has monodromy $-Id$
  around $D_w$. This definition is justified by
  remark~\ref{rem:apresassumptions}, (2) above.
  
  2) For $\fA^{\aaa\bbb}_i \not = 0$ (i.e. $i\in I_{\aaa\bbb} \sqcup
  I'_{\aaa\bbb}$) and $Y_i \subset D_v$, we have: $v\in \Delta'_\alpha
  \Longleftrightarrow v\in \Delta'_\beta$. Hence:
  \begin{equation}
    \label{eq:intersection_ouv_div}
  \textstyle
  \text{for $\fA^{\aaa\bbb}_i \not = 0$,} \qquad
  U_i \cap \bigcup_{v\in \Delta'_\alpha} D_v = U_i \cap \bigcup_{v\in
    \Delta'_\beta} D_v  = U_i \cap \bigcup_{v\in \Delta'_\alpha \cap
    \Delta'_\beta} D_v .
  \end{equation}
\end{remarks}
We will introduce an algebra structure on $\fA$. For $v\in V$, $D_v$
has a fundamental class, $\delta_v\in H^2_{D_v}(Y;\C_Y)$. We choose
representatives, $\xi_v\in \sect(Y;\Omega^2_Y)$, of the $\delta_v$.
For $\Delta,\Delta',\Delta'' \in \S$, we set
\begin{equation}
  \label{eq:nabla}
  \nabla(\Delta,\Delta',\Delta'') = 
 \bigl(\Delta' \setminus (\Delta\cup  \Delta'') \bigr) 
\cup  \bigl( (\Delta\cap \Delta'') \setminus \Delta' \bigr)
\end{equation}
and for $\aaa, \bbb, \ccc \in \epm$, $\eta_{\aaa\bbb\ccc} =
\prod_{v\in \nabla} \xi_v$, where $\nabla = \nabla(\Delta_\aaa,
\Delta_\bbb, \Delta_\ccc)$.  For $\aaa,\bbb,\ccc \in \epm$, we define
a morphism $m^{\aaa\bbb\ccc} : \fA^{\bbb\ccc} \otimes \fA^{\aaa\bbb}
\to \fA^{\aaa\ccc}$ as follows. For $i\in \pxf(Z_\aaa \cap Z_\bbb \cap
Z_\ccc)$, we define a sheaf morphism
\begin{align*}
  n_{\aaa\bbb\ccc}^i : \Omega_{\bbb\ccc, i}
  \otimes \Omega_{\aaa\bbb, i}   &\to \Omega_{\aaa\ccc, i}  \\
  (\tau \otimes v)\otimes  (\sigma\otimes u)  &\mapsto
  (\eta_{\aaa\bbb\ccc} \, \tau\, \sigma) \otimes (v\circ u),
\end{align*}
where $\sigma, \tau$ are sections of $\Omega_{\pxf^{-1}(U_i)}$ and
$u,v$ sections of $\hom$ sheaves.  We set $m^{\aaa\bbb\ccc}_i =
\sect(\pxf^{-1}(U_i); n_{\aaa\bbb\ccc}^i )$.  This definition extends
to other $i\in I$, either by restriction to the case $i\in \pxf(Z_\aaa
\cap Z_\bbb \cap Z_\ccc )$ or, trivially, when one of the terms is
$0$. Indeed, if $i\in I\setminus \pxf(Z_\aaa \cap Z_\bbb \cap Z_\ccc
)$ satisfies $\fA^{\bbb\ccc}_i \not= 0$, $\fA^{\aaa\bbb}_i \not=0$
and $\fA^{\aaa\ccc}_i \not= 0$, then we have,
by~\eqref{eq:intersection_ouv_div},
$$
\textstyle
U_i \setminus \pxf (\bigcup_{v\in \Delta'_\bbb \cap \Delta'_\ccc}
D_v) = U_i \setminus \pxf (\bigcup_{v\in \Delta'_\aaa \cap
  \Delta'_\bbb} D_v) = U_i \setminus \pxf (\bigcup_{v\in \Delta'_\aaa
  \cap \Delta'_\ccc} D_v).
$$
It follows that, for the same $j\in \pxf(Z_\aaa \cap Z_\bbb \cap
Z_\ccc )$, we have $\fA^{\bbb\ccc}_i = \fA^{\bbb\ccc}_j$,
$\fA^{\aaa\bbb}_i = \fA^{\aaa\bbb}_j$, $\fA^{\aaa\ccc}_i
=\fA^{\aaa\ccc}_j$. This allows one to define $m^{\aaa\bbb\ccc}_i$.
By definition, these morphisms $m^{\aaa\bbb\ccc}_\cdot$ commute with
the restriction maps and we obtain a sheaves morphism
$m^{\aaa\bbb\ccc} : \fA^{\bbb\ccc} \otimes \fA^{\aaa\bbb} \to
\fA^{\aaa\ccc}$, as claimed. (The justification for the definition of
this product is given in section~\ref{algebrastructure}.)

Now we define a product $m$ on $\fA = \oplus_{\aaa,\bbb\in \epm}
\fA^{\aaa\bbb}$ by $m = \oplus m^{\aaa\bbb\ccc}$.  One checks that $m$
is an associative product using the straightforward identity:
\begin{equation}
  \label{eq:identiteensembliste}
\eta_{\aaa\bbb\ccc} \, \eta_{\aaa\ccc\ddd} 
= \eta_{\bbb\ccc\ddd} \, \eta_{\aaa\bbb\ddd}.
\end{equation}
Hence $\fA$ is a sheaf of dg-algebras on $I$.  For $\aaa \in \epm$,
$N_\aaa = \oplus_{\aaa' \in \epm}\fA^{\aaa'\aaa}$ has a natural
structure of $\fA$-module defined in the same way as the product of
$\fA$. The result of this section is the following equivalence of
categories.
\begin{proposition}
  \label{prop:eq_cat1}
  Let $Y=\bigsqcup_{i\in I} Y_i$ be a stratified complex analytic
  manifold, endowed with normal crossings divisors $D_v$, $v\in V$,
  and sheaves $L_\aaa$, $\aaa\in \epm$, satisfying
  assumptions~\ref{assumptionstratification}. For a choice of forms
  $\xi_v \in \sect(Y;\Omega^2_Y)$, we define a sheaf of dg-algebras
  $\fA$ on $I$, and $\fA$-modules $N_\aaa$ as above.
  
  Then, there exists a choice of $\xi_v$ such that we have an
  equivalence of categories between $\D(Y)\langle L_\aaa \rangle$ and
  $\D_\fA \langle N_\aaa \rangle$ sending $L_\aaa$ to $N_\aaa$.
\end{proposition}
The proof is given at the end of this section.
\begin{remarks}
\label{rem:espclass}
1) In fact one could prove that two choices of representatives
$\xi_v$, $\xi'_v$ of the $\delta_v$ give quasi-isomorphic sheaves of
dg-algebras, $\fA$, $\fA'$: with $\zeta_v$ such that $\xi_v - \xi'_v =
d\zeta_v$, and replacing $Y$ by $Y\times \C$, endowed with $\xi^+_v
=\xi_v +d(t \zeta_v)$ ($t$ is the coordinate on $\C$) and the data $Y_i
\times \C$, $\D_v \times \C$, $L_\alpha \boxtimes \C_{\C}$, we could
build a third sheaf $\fA^+$ on $I$, quasi-isomorphic to $\fA$ and
$\fA'$. Hence the conclusion of the proposition is valid for any
choice of $\xi_v$, but we will not use this result.

2) The results of this section will be applied to $Y=E\times_K X$,
where $X$ is a symmetric variety under the action of a semi-simple
complex algebraic group $G$, $K$ a suitable maximal compact subgroup
of $G$ and $E$ a universal bundle for $K$. Of course, $E$ is not a
manifold, but we may assume that it is an increasing union of
$K$-manifolds, $E=\bigcup_k E_k$, and consider the de Rham algebra of
$Y$, $\Omega_Y = \varprojlim_k \Omega_{E_k\times_K X}$. The
stratification $Y=\bigsqcup_i Y_i$ and the divisors $D_v$ will be
given by a $K$-invariant stratification of $X$ and $K$-invariant
divisors.  All constructions in this section can be made $K$-invariant
(if we choose $K$-invariant functions $f_i$ in
notations~\ref{notationtubes} below, by averaging under the action of
$K$) and transpose to $Y=E\times_K X$.
\end{remarks}

\subsection{System of tubes}
Our first task in the proof of proposition~\ref{prop:eq_cat1} is to
``compute'' the global sections $\Gamma_{\aaa\bbb i}
=\rsect(\pxf^{-1}(U_i); \Rhom(L_\aaa,L_\bbb))$, for $i\in I$,
$\aaa,\bbb \in \epm$ (to compute just means to find suitable
representatives).  For this we replace the strata $Y_i$ by a system of
``tubes'', $T_i$ (with $T_i$ closed enough to $Y_i$ so that the local
systems $L_\aaa|_{Y_i}$ extend to $T_i$) with the properties: (i)
replacing $L_\aaa$ by its extension, say $L'_\aaa$, to the union of
tubes $\bigcup_{\{i; Y_i \subset Z_\aaa\}} T_i$ doesn't change the
global sections $\Gamma_{\aaa\bbb i}$, (ii) the complexes
$\Rhom(L'_\aaa,L'_\bbb)$ are in fact sheaves.  The precise statement
is given in proposition~\ref{prop:plomberie} below.  The first
property implies that the category $\D(Y)\langle L_\aaa \rangle$ is
equivalent to the $\D(Y)\langle L'_\aaa \rangle$ (see
lemma~\ref{equiv_LetL'} below).  The second property will be used to
define a sheaf of dg-algebras, $\fB$, on $I$, such that $\D(Y)\langle
L'_\aaa \rangle$ is equivalent to a subcategory of $\D_{\fB}$
(definition~\ref{def:faisceauA} and proposition~\ref{prop:eqL'M}).
The proof of proposition~\ref{prop:eq_cat1} will then be achieved by
showing that $\fB$ and $\fA$ are quasi-isomorphic.

\begin{notations}
\label{notationtubes}
First we assume that the finite set indexing the stratification is
$I=\{1,\ldots,n\}$, ordered such that $\dim Y_i \leq \dim Y_{i+1}$,
for $i=1,\ldots,n-1$. Recall that $Y$ is open in an analytic manifold
$X$ and $\ovl{Y}$ is compact, with a stratification
$\ovl{Y}=\bigsqcup_{i \in I} Y'_i $, satisfying: $\forall i\in I$,
$Y_i =Y \cap Y'_i$.  For $i=1,\ldots,n-1$, we choose a neighbourhood
of $Y_i$, $\tilde{Y}_i$, whose closure is a neighbourhood of $Y'_i$ in
$\ovl{Y}$. We also choose real analytic functions $f_i:\tilde{Y}_i \to
\R$, such that $f_i(\tilde{Y}_i) \subset \R_{\geq0}$ and $Y_i =
f_i^{-1}(0) \cap \tilde{Y}_i$. For $k<n$ and
$\varepsilon_1,\ldots,\varepsilon_k >0$, we define
$$
T_1(\varepsilon_1) = \{y\in \tilde{Y}_1; \: f_1(y)\leq
\varepsilon_1 \}, \ldots, T_k(\varepsilon_1,\ldots,\varepsilon_k) =
\{y\in \tilde{Y}_k; \: f_k(y)\leq \varepsilon_k \}\setminus
\bigsqcup_{i<k}T_i.
$$
By abuse of notations we will write $T_i(\udl\varepsilon) =
T_i(\varepsilon_1,\ldots,\varepsilon_i)$ for any $\udl\varepsilon$ of
length greater than $i$. We also set $T_n(\udl\varepsilon) = Y
\setminus \bigsqcup_{i<n}T_i(\udl\varepsilon)$.  For a union of
strata, $Z$, we set:
\begin{equation}
  \label{eq:defTZ}
 \textstyle  T_Z(\udl{\varepsilon}) = 
\bigsqcup_{i\in J} T_i(\udl\varepsilon),
\quad \text{where $J$ satisfies} \quad  Z= \bigsqcup_{i\in J}Y_i.
\end{equation}
\end{notations}
\begin{definition}
  We call ``set of bounds'' a subset $B\in ]0,+\infty[^k$ such that
  $\exists \varepsilon_1^0>0$, $\forall \varepsilon_1 <
  \varepsilon_1^0$, $\exists \varepsilon_2^0>0$, $\forall
  \varepsilon_2 < \varepsilon_2^0$, \dots, $\exists
  \varepsilon_k^0>0$, $\forall \varepsilon_k < \varepsilon_k^0$,
  $(\varepsilon_1,\ldots,\varepsilon_k)\in B$.
\end{definition}
We note that a set of bounds is non-empty and that the intersection of
two sets of bounds is a set of bounds. The aim of this paragraph is to
prove the following result.
\begin{proposition}
  \label{prop:plomberie}
  Let $Y=\bigsqcup_{i=1,\ldots,n} Y_i$ be a stratified analytic
  manifold as above.  Let $Z_1, Z_2 \subset Y$ be locally closed
  subsets of $Y$ which are unions of strata.  Let $\Ll^1$, $\Ll^2$ be
  local systems respectively defined on neighbourhoods of $Z_1$ and
  $Z_2$.  Then there exists a set of bounds $B\subset
  ]0,+\infty[^{n-1}$ such that $\forall \udl\varepsilon\in B$, setting
  $T_Z = T_Z(\udl\varepsilon)$, $\Ll^i$ is defined on $T_{Z_i}$ and we
  have:
  
  (i) There exist natural morphisms $\Ll^i_{T_{Z_i}} \to \Ll^i_{Z_i}$
  and they induce isomorphisms
  \begin{equation}
    \label{eq:iso_LT_LZ}
    \RHom(\Ll^1_{T_{Z_1}},\Ll^2_{T_{Z_2}}) \isoto
    \RHom(\Ll^1_{T_{Z_1}},\Ll^2_{Z_2})
    \isofrom \RHom(\Ll^1_{Z_1},\Ll^2_{Z_2}).
  \end{equation}
  
  (ii) For any locally closed union of strata $Z\subset Y$ such that
  $Z_1$, $Z_2\subset Y\setminus \ovl{Z}$, we have an isomorphism
  \begin{equation}
  \label{eq:iso_LT_LZ_bord}
  \rsect(Y;(\Rhom(\Ll^1_{T_{Z_1}},\Ll^2_{T_{Z_2}}))_{T_Z})
  \simeq \rsect(Y;(\Rhom(\Ll^1_{Z_1},\Ll^2_{Z_2}))_Z).
  \end{equation}
  
  (iii) If $Z_1\subset Z_2$, then $\Rhom(\Ll^1_{T_{Z_1}},
  \Ll^2_{T_{Z_2}})$ is concentrated in degree $0$.
  
  (iv) Let us assume that $Z_2$ and $Z_1\cap Z_2$ are smooth and let
  $\mu:Z_2 \hookrightarrow \ovl{Z_2}$, $\nu:Z_1\cap Z_2
  \hookrightarrow \ovl{Z_1\cap Z_2}$ be the inclusions. We assume that
  $R\mu_*(\Ll^2|_{Z_2}) = \mu_!(\Ll^2|_{Z_2})$ and
  $R\nu_*(\Ll^2|_{Z_1\cap Z_2}) = \nu_!(\Ll^2|_{Z_1\cap Z_2})$. Then,
  for any open union of strata $V\subset Y$,
  \begin{equation}
    \label{eq:iso_LT_LZ_degre0}
    \RHom(\Ll^1_{Z_1}|_V,\Ll^2_{Z_2}|_V) \simeq 
    \rsect(T_V; \hom(\Ll^1_{T_{Z_1}},\Ll^2_{T_{Z_2}})).
  \end{equation}
\end{proposition}
The proof will be given at the end of the paragraph.  We first deduce
the following corollary, which gives a category equivalent to
$\D(Y)\langle L_\aaa \rangle$.

We consider $Y$ and $L_\aaa$, $\aaa \in \epm$, as in
assumptions~\ref{assumptionstratification}. We choose neighbourhoods
of the $Z_\aaa$ on which the local systems $L_\aaa$ may be extended to
local systems $L^+_\aaa$.  For $\aaa, \bbb \in \epm$, we set $Z_1 =
Z_\aaa$, $\Ll^1 = L^+_\aaa$, $Z_2 = Z_\bbb$, $\Ll^2 = L^+_\bbb$.
Then, $Z_2$ and $Z_1 \cap Z_2$ are open subsets of intersections of
some $D_v$, hence smooth. Moreover, by
assumptions~\ref{assumptionstratification} (v), $\Ll^2$ has monodromy
$-Id$ around each irreducible divisor of $\ovl{Z_2} \setminus Z_2$ and
the similar property holds a fortiori for $\Ll^2|_{Z_1\cap Z_2}$.
Hence, by lemma~\ref{lem:monodromy}, the hypothesis of
proposition~\ref{prop:plomberie}, (iv), are verified.
\begin{notations} 
\label{notationsystemdeetubes}
We choose a set of bounds, $B$, such that the conclusions of
proposition~\ref{prop:plomberie} hold for $Z_1 = Z_\aaa$, $\Ll^1 =
L^+_\aaa$, $Z_2 = Z_\bbb$, $\Ll^2 = L^+_\bbb$, for any pair
$(\aaa,\bbb) \in \epm^2$. We fix $\udl\varepsilon \in B$ and set:
$$
T_i = T_i(\udl\varepsilon), \qquad T_\aaa =
T_{Z_\aaa}(\udl\varepsilon), \qquad L'_\aaa = (L^+_\aaa)_{T_\aaa}.
$$
\end{notations}
\begin{corollary}
\label{equiv_LetL'}
  The categories $\D(Y)\langle L_\aaa \rangle$ and $\D(Y)\langle
  L'_\aaa \rangle$ are equivalent.
\end{corollary}
\begin{proof}
  This is a consequence of the natural
  isomorphisms~\eqref{eq:iso_LT_LZ}.  By definition the category
  $\D(Y)\langle L_\aaa \rangle$ is the union of the full
  subcategories $\D_n$, $n\in \N$, where $\D_0$ consists of the
  $L_\aaa[k]$, $\aaa\in \epm$, $k\in \Z$, and $\D_{n+1}$ is obtained
  from $\D_n$ by adding objects $H$ appearing in distinguished
  triangles $F\to G \to H \fintd$, with $F, G\in \D_n$. We write in
  the same way $\D(Y)\langle L'_\aaa \rangle = \bigcup \D'_n$.  We
  assume by induction that we have an equivalence, $\delta_n$, between
  $\D_n$ and $\D'_n$, together with functorial morphisms
  $r_n(F):\delta_n(F) \to F$ such that $\delta_n$ is given on the
  morphisms by composing isomorphisms
  $$
  \Hom_{\D_n}(F,F') \isoto \Hom_{\D(Y)}(\delta_n(F),F') \isofrom
  \Hom_{\D'_n}(\delta_n(F),\delta_n(F'))
  $$
  induced by $r_n(F)$, $r_n(F')$.  (The first step is given
  by~\eqref{eq:iso_LT_LZ}.)  Let $F\xto{u} G
  \to H \fintd$ be a distinguished triangle as above,
  $F'=\delta_n(F)$, $G'=\delta_n(G)$ $u'=\delta_n(u)$ and consider a
  distinguished triangle $F'\xto{u'} G' \to H' \fintd$.  We extend
  the square built on $u$, $u'$, $r_n(F)$ and $r_n(G)$ to a morphism
  of triangles:
  $$
  \xymatrix@R=4mm{
    F' \ar[d]\ar[r] & G' \ar[d]\ar[r] & H' \ar[d]^r\ar[r]^{+1} &  \\
    F \ar[r] & G \ar[r] & H \ar[r]^{+1} & . }
  $$
  We set $H'=\delta_{n+1}(H)$, $r_{n+1}(H) = r$ and we have to define
  the images of the morphisms. First, for $X\in \D_n$, we have long
  exact sequences of homomorphisms groups:
  $$
  \xymatrix@C=5mm@R=4mm{ \Hom(X,F) \ar[d]\ar[r] & \Hom(X,G) \ar[d]\ar[r]
    & \Hom(X,H) \ar[d]\ar[r] &  \\
    \Hom(\delta_n(X),F) \ar[r] & \Hom(\delta_n(X),G)\ar[r]
    & \Hom(\delta_n(X),H) \ar[r] &  \\
    \Hom(\delta_n(X),\delta_n(F)) \ar[r]\ar[u] &
    \Hom(\delta_n(X),\delta_n(G)) \ar[r]\ar[u]
    & \Hom(\delta_n(X),\delta_n(H)) \ar[r]\ar[u] & . \\
  }
  $$
  By the five lemma it gives an isomorphism $\Hom(X,H)\simeq
  \Hom(\delta_n(X),\delta_n(H))$, which we use to define $\delta_n$ on
  $\Hom(X,H)$. In the same way, we may define $\delta_n$ on
  $\Hom(H,X)$, still for $X\in\D_n$. Then we may assume $X\in
  \D_{n+1}$ in the above diagram, and this defines $\delta_n$ on
  $\Hom(X,H)$ for $X,H\in \D_{n+1}$, satisfying the compatibility with
  $r_{n+1}$.
\end{proof}
Now we give some preliminary results before we prove
proposition~\ref{prop:plomberie}.
\begin{lemma}
  \label{tubes_fermes}
  Let $Y=\bigsqcup_{i=1,\ldots,n} Y_i$, $f_i$, be as in
  assumptions~\ref{assumptionstratification} and
  notations~\ref{notationtubes}.  There exists a set of bounds
  $B\subset ]0,+\infty[^{n-1}$ such that $\forall \udl\varepsilon \in
  B$ and any union of strata $Z\subset Y$:
  \begin{itemize}
  \item [(i)] if $Z$ is closed then so is $T_Z(\udl\varepsilon)$, and
    $Z\subset T_Z(\udl\varepsilon)$,
  \item[(ii)] if $Z$ is open then so is $T_Z(\udl\varepsilon)$,
  \item[(iii)] $\forall i_1<\cdots < i_p < n$, and $y\in Y$ such that
    $f_{i_l}(y) = \varepsilon_{i_l}$, we have $df_{i_1} \wedge \ldots
    \wedge df_{i_p} (y) \not=0$. In particular, locally around any
    point $y\in Y$, the partition $Y=\bigsqcup T_i(\udl\varepsilon)$
    is homeomorphic to $\R^d = \{x_1 \leq 0 \} \sqcup \{x_1 >0, \, x_2
    \leq 0 \} \sqcup \{x_1>0,\ldots, x_{q-1} >0, \, x_q \leq 0 \}
    \sqcup \{x_1>0,\ldots, x_q >0\}$, for some $q$.
  \end{itemize}
\end{lemma}
\begin{proof}
  Of course (ii) follows from (i) because $T_{Y\setminus
    Z}(\udl\varepsilon) = Y \setminus T_Z(\udl\varepsilon)$.  We prove
  (i) by induction on $n$, the case $n=1$ or $2$ being obvious. Recall
  that $Y$ is open in a manifold $X$. For $i$ such that $Y_1 \cap
  \ovl{Y_i} =\emptyset$, we also have $\ovl{Y_1}\smash{{}^X} \cap
  \ovl{Y_i}\smash{{}^X}=\emptyset$, because otherwise the
  stratification of $\ovl{Y}\smash{{}^X}$ would have additional
  strata.  Since $\ovl{Y}\smash{{}^X}$ is compact, we deduce
  $d(Y_1,Y_i) >0$.  Hence we may choose $r >0$ smaller than $\min\{
  d(Y_1,Y_i);\: Y_1 \cap \ovl{Y_i} = \emptyset\}$ and $\sup\{
  d(Y_1,y);\: y\in Y_j\}$, for all $j$. Then, for $\varepsilon^0_1$
  such that $T_1(\varepsilon^0_1) \subset \{y\in Y; d(Y_1,y) <r \}$
  and for $0 < \varepsilon_1 < \varepsilon^0_1$, we have
  $T_1(\varepsilon_1) \cap \ovl{Y_i} \not= \emptyset$ if and only if
  $Y_1 \subset \ovl{Y_i}$, and moreover $\forall j\not= 1$,
  $Y_j\not\subset T_1(\varepsilon_1)$.
  
  The induction hypothesis applied to $Y' = Y\setminus
  T_1(\varepsilon_1)$ stratified by the $Y'_i = Y'\cap Y_i$ gives a
  set of bounds $B'(\varepsilon_1)\subset ]0,+\infty[^{n-2}$ for which
  (i) holds in $Y'$. For $i$ such that $Y_1\cap \ovl{Y_i}\smash{{}^Y}
  =\emptyset$ we may choose $\varepsilon_i$ small enough so that
  $T_1(\varepsilon_1) \cap \{y\in \tilde{Y}_i; f_i(y) \leq
  \varepsilon_i\} =\emptyset$. In particular, restricting to a smaller
  set of bounds $B''(\varepsilon_1)$, we may assume that
  $T_1(\varepsilon_1) \cap \ovl{T_i(\udl\varepsilon)} \smash{{}^Y}
  =\emptyset$.  Let $Z\subset Y$ be closed.
  
  If $Y_1\subset Z$ then $T_Z = T_1(\varepsilon_1) \sqcup T_{Z\cap
    Y'}$. Since $T_{Z\cap Y'}$ is closed in $Y'$, $T_Z$ is closed in
  $Y = T_1(\varepsilon_1) \sqcup Y'$.  By induction we also have
  $Z\cap Y' \subset T_Z\cap Y'$ and this implies $Z\subset T_Z$.
  
  If $Y_1\not\subset Z$ then $Z$ only contains strata $Y_i$ such that
  $Y_1 \cap \ovl{Y_i} =\emptyset$, so that $T_1(\varepsilon_1) \cap
  \ovl{Y_i} = \emptyset$. It follows that $Z\subset Y'$ (and $Z$ is
  closed in $Y'$). This also gives $T_1(\varepsilon_1) \cap
  \ovl{(T_Z)}\smash{{}^Y} =\emptyset$ and, since $T_Z$ is closed in
  $Y'$, it is closed in $Y$ too. Finally $Z\subset T_Z$ since this is
  already true in $Y'$.  In conclusion the set of bounds $B =
  \{(\varepsilon_1,\ldots,\varepsilon_{n-1});\: 0< \varepsilon_1 <
  \varepsilon^0_1,\; (\varepsilon_2,\ldots,\varepsilon_{n-1}) \in
  B''(\varepsilon_1) \}$ has the required property.
  
\medskip

  Now we prove by induction on $p$ that there exists a set of bounds
  $B_p$ such that the conclusion of (iii) holds for any $i_1<\cdots <
  i_p < n$. For $p=1$ this is a consequence of the curve selection
  lemma: by contradiction, if the closure of $\{y\in Y_{i_1};\:
  df_{i_1}(y) = 0\}$ intersects $\ovl{Y_1}$, then there exists a real
  analytic curve $\gamma:]-1,1[ \to \ovl{Y}$ such that $\gamma(t) \in
  Y\setminus Y_1$ for $t\not= 0$ and $\gamma(0) \in \ovl{Y_1}$. But
  this implies $df_{i_1}( \gamma(t)) = 0$ so that $f_{i_1}(
  \gamma(t))$ is constant, which is impossible.  Hence there exists
  $\varepsilon^0_{i_1} >0$ such that $0 < f_{i_1}(y) <
  \varepsilon^0_{i_1}$ implies $df_{i_1}(y) \not= 0$.  We take $B_1 =
  \prod_i ]0, \varepsilon^0_i [$.
  
  Assuming (iii) holds for $p$, we consider, for $\udl\varepsilon \in
  B_{p}$ and $i_1 < \cdots < i_p$, the smooth subvariety of $Y$, $Y' =
  \{y\in Y;\: f_{i_l}(y) = \varepsilon_{i_l}, \: l=1,\ldots, p\}$.
  For $i_{p+1} > i_p$, the function $f_{i_{p+1}}$ is not constant on
  $Y'$ and the proof of the first step gives
  $\varepsilon^0_{i_1,\ldots,i_p} (\varepsilon_{i_1},\ldots,
  \varepsilon_{i_p}) > 0$ such that the conclusion holds for
  $(\varepsilon_{i_1},\ldots, \varepsilon_{i_{p+1}})$ with
  $\varepsilon_{i_{p+1}} < \varepsilon^0_{i_1,\ldots,i_p}
  (\varepsilon_{i_1},\ldots, \varepsilon_{i_p})$.  We set $B^p_p =
  B_p$ and, for $k=p+1,\ldots, n$, $B^k_p= \{ \udl\varepsilon \in
  B^{k-1}_p;\: \forall i_1 < \cdots < i_p < k,\, \varepsilon_k <
  \varepsilon^0_{i_1,\ldots,i_p} (\varepsilon_{i_1},\ldots,
  \varepsilon_{i_p}) \}$. Then $B_{p+1} = B^n_p$ is a set of bounds
  with the required property for step $p+1$ and we take $B= B_{n-1}$.
  
  Now, for $\varepsilon \in B$ and $y\in Y$, we let $i_1 < \cdots <
  i_q$ be the indices such that $f_{i_l}(y) = \varepsilon_{i_l}$.
  Since $df_{i_1} \wedge \ldots \wedge df_{i_q} (y) \not=0$ the
  functions $x_l = f_{i_l} - \varepsilon_{i_l}$, $l=1,\ldots,q$, may
  be extended to a coordinates system around $y$ and in any such
  system the description of the partition is the one given in the
  lemma.
\end{proof}
\begin{lemma}
  \label{bons_tubes}
  Let $Y=\bigsqcup_{i=1,\ldots,n} Y_i$, $f_i$, be as in
  assumptions~\ref{assumptionstratification} and
  notations~\ref{notationtubes}.  Let $F_1,\ldots,F_m \in
  \D^b_{\R-c}(Y)$ (i.e. $F_j$ is constructible for some stratification
  of $Y$, not a priori $(Y_i)_{i\in I}$). Then there exists a set of
  bounds $B\subset ]0,+\infty[^{n-1}$ such that $\forall
  \udl\varepsilon \in B$, setting for short $T_i =
  T_i(\udl\varepsilon)$, we have isomorphisms $\forall i=1,\ldots,n$,
  $\forall j= 1,\ldots,m$:
  $$
  H_{Y_i}(Y;F_j) \isoto H_{T_i}(Y;F_j), \qquad H(Y;(F_j)_{T_i})
  \isoto H(Y;(F_j)_{Y_i}).
  $$
\end{lemma}
\begin{proof}
  We prove by induction on $k$ that there exists a set of bounds
  $B\subset ]0,+\infty[^{k}$ such that the first isomorphism holds for
  any $j$ and for $i=1,\ldots,k$. For $k=1$, this is a consequence of
  lemma~\ref{lemKS}.
  
  Let us assume the result is true for $k$, and apply it to $F_j$ and
  $F'_j = \rsect_{Y_{k+1}}F_j$, $j= 1,\ldots,m$. Let $B\subset
  ]0,+\infty[^{k}$ be the set of bounds obtained and
  $(\varepsilon_1,\ldots, \varepsilon_{k})\in B$, $T_i$,
  $i=1,\ldots,k$ as in the lemma. Let us set $T=T_1\sqcup\ldots \sqcup
  T_k$ and $U=Y\setminus T$. By lemma~\ref{lemKS} again, applied to
  $f=f_{k+1}$ and the complexes $\rsect_U(F_j)$, there exists
  $\varepsilon_{k+1}^0 >0$ such that $\forall 0< \varepsilon_{k+1} <
  \varepsilon_{k+1}^0$ we have, setting $T_{k+1} =
  f_{k+1}^{-1}([0,\varepsilon_{k+1}]) \cap \tilde{Y}_{k+1} \cap U$,
  $H_{Y_{k+1}}(U;F_j) \isoto H_{T_{k+1}}(U;F_j)$. Since $T_{k+1}
  \subset U$, we have $H_{T_{k+1}}(U;F_j) \simeq H_{T_{k+1}}(Y;F_j)$
  and we have to prove that $H_{Y_{k+1}}(U;F_j) \simeq
  H_{Y_{k+1}}(Y;F_j)$. Using an excision exact sequence, we are
  reduced to proving the vanishing of $A_j = H_{Y_{k+1}\cap
    T}(Y;F_j)$. Now $A_j \simeq H_T(Y;F'_j)$ and, for $i\leq k$, we
  have
  $$
  H_{Y_i}(Y;F'_j) \isoto H_{T_i}(Y;F'_j) \quad \text{and} \quad
  H_{Y_i}(Y;F'_j) \simeq H_{Y_i\cap Y_{k+1}}(Y;F_j) =0.
  $$
  Let us set $T'_i = T_1\sqcup\ldots\sqcup T_i$; we have
  distinguished triangles $\rsect_{T'_{i-1}}(\cdot) \to
  \rsect_{T'_i}(\cdot) \to \rsect_{T_i}(\cdot) \fintd$. We deduce by
  induction on $i\leq k$ that $H_{T'_i}(Y;F'_j) =0$. For $i=k$ we
  obtain $A_j =0$ as desired and this concludes the proof of the first
  isomorphism.
  
  \smallskip The proof of the second isomorphism is the same, again
  using lemma~\ref{lemKS} (we just note that we apply the induction
  hypothesis to $F_j$ and $F''_j= (F_j)_{Y_{k+1}}$, and
  lemma~\ref{lemKS} to $(F_j)_U$.)
\end{proof}
We still consider $Y$ satisfying
assumptions~\ref{assumptionstratification} and we keep
notations~\ref{notationtubes}. We consider moreover local systems
$L^i$ defined on neighbourhoods of $Y_i$.  The sheaves $L^i_{Y_i}$ are
well-defined and, for $\udl\varepsilon \in ]0,+\infty[^{n-1}$ small
enough, the $L^i_{T_i(\udl\varepsilon)}$ are also well-defined.
\begin{lemma}
  \label{hypersurf_transv}
  Let $Y=\bigsqcup_{i=1,\ldots,n}Y_i$, $L^i$, be as above.  Let
  $U\subset Y$ be an analytic open subset, $U'$ a neighbourhood of
  $\ovl{U}$ and $f:U'\to \R$ an analytic function with $1$ as regular
  value. We assume that the smooth hypersurface $S=\{y\in
  U;\:f(y)=1\}$ meets the strata $U\cap Y_i$ transversally. Let $L$ be
  a sheaf on $Y$ which is a local system in a neighbourhood of $S$. We
  set $U_+= \{y;\:f(y) > 1\}$.  Then there exists a set of bounds
  $B\subset ]0,+\infty[^{n-1}$ such that $\forall \udl\varepsilon \in
  B$, setting $T'_i = U\cap T_i(\udl\varepsilon)$, we have an
  isomorphism $\Rhom(L_S,L^i_{T'_i}|_U) \simeq (L^*\otimes
  L^i|_U)_{S\cap T'_i}[-1]$ and the morphisms
  $$
  \RHom(L_S,L^i_{T'_i\cap U_+}) \to\RHom(L_S,L^i_{T'_i}) \to
  \RHom(L_S,L^i_{Y_i})
  $$
  are isomorphisms.
\end{lemma}
\begin{proof}
  We have $\Rhom(L_S,\C_U) \simeq L^*_S[-1]$ because $S$ is a smooth,
  relatively oriented, hypersurface. The micro-support of $L_S$ is
  $SS(L_S) = T^*_SU$ and we also have the bound $SS(L^i_{Y_i}) \subset
  \bigcup_i T^*_{Y_i}U$. By the transversality hypothesis, we have
  $SS(\C_S) \cap SS(L^i_{Y_i}) \subset T^*_UU$.  Hence, by
  lemma~\ref{lemmicrosupport}, we have isomorphisms
  $$
  \Rhom(L_S,L^i_{Y_i}|_U) \simeq \Rhom(L_S,\C_U)\otimes
  L^i_{Y_i}|_U \simeq (L^*\otimes L^i|_U)_{S\cap Y_i}[-1].
  $$
  For $\varepsilon_i$ small enough, $\partial T'_i$ also is
  transversal to $S$ and, since $SS(L^i_{T'_i})$ is the outer conormal
  of $\partial T'_i$ in $U$, we obtain similarly
  $\Rhom(L_S,L^i_{T'_i}|_U) \simeq (L^*\otimes L^i|_U)_{S\cap
    T'_i}[-1]$.
  
  Hence, to show the last isomorphism, it is sufficient to find a set
  of bounds such that the morphisms $H(S;(L^*\otimes L^i|_U)_{S\cap
    T'_i}) \to H(S;(L^*\otimes L^i|_U)_{S\cap Y_i})$ are isomorphisms.
  But this follows from lemma~\ref{bons_tubes} applied to the
  stratification $S=\bigsqcup_i S\cap Y_i$.
  
  Let us prove that the remaining morphism also is an isomorphism.
  Let us set $U_- = f^{-1}(]-\infty,1[)$. We have
  $\Rhom(\C_{U_-},\C_U) \simeq \C_{\ovl{U}_-}$ and $SS(\C_{U_-})$ is
  the inner conormal to $U_-$.  We deduce as above
  $\Rhom(\C_{U_-},L^i_{T'_i}) \simeq L^i_{T'_i \cap \ovl{U}_-}$ and:
  $$
  \Rhom(L_S,L^i_{T'_i \cap \ovl{U}_-}) \simeq \Rhom(L_S,
  \Rhom(\C_{U_-},L^i_{T'_i})) \simeq \Rhom(L_{S \cap U_-},L^i_{T'_i})
  =0.
  $$
  Using the distinguished triangle $L^i_{T'_i\cap U_+} \to L^i_{T'_i}
  \to L^i_{T'_i\cap \ovl{U}_-} \fintd$, we conclude that
  $\RHom(L_S,L^i_{T'_i\cap U_+}) \simeq \RHom(L_S,L^i_{T'_i})$ as
  desired.
\end{proof}

\begin{proposition}
  \label{prop_tubes}
  Let $Y=\bigsqcup_{i=1,\ldots,n}Y_i$, $L^i$, be as in
  lemma~\ref{hypersurf_transv}. For $\udl\varepsilon \in
  ]0,+\infty[^{n-1}$ and a union of strata $Z\subset Y$, we set:
  $$
  \textstyle
  T_i =T_i(\udl\varepsilon), \quad T_Z= T_Z(\udl\varepsilon),
  \qquad \text{for $i<j$:} \quad  T_{ij} = (\ovl{T_i} \cap \ovl{T_j})
  \setminus (\bigcup_{i<l<j} \ovl{T_l}).
  $$
  There exists a set of bounds $B\subset ]0,+\infty[^{n-1}$ such
  that $\forall \udl\varepsilon\in B$:
  
  (i) for any $i<j$, the morphism $\alpha_{ij}: (L^{i*}\otimes
  L^j)_{T_{ij}} [-1] \to \Rhom(L^i_{T_i}, L^j_{T_j})$ is an
  isomorphism,
  
  (ii)  for any $i, j$, the natural morphisms
  $$
  \RHom(L^i_{T_i},L^j_{T_j}) \xto{a_{ij}} \RHom(L^i_{T_i},L^j_{Y_j})
  \xfrom{b_{ij}} \RHom(L^i_{Y_i},L^j_{Y_j})
  $$
  are isomorphisms,
  
  (iii) for any $i,j,p$ with $i,j> p$, we have an isomorphism
  $$
  \rsect(Y;(\Rhom(L^i_{T_i}, L^j_{T_j}))_{T_{Y_p}}) \simeq
  \rsect(Y;(\Rhom(L^i_{Y_i}, L^j_{Y_j}))_{Y_p}).
  $$
\end{proposition}
\begin{proof}
  (i) and (ii). By lemma~\ref{bons_tubes}, we know that the $b_{ij}$ are
  isomorphisms for $\udl\varepsilon$ in a suitable set of bounds. We
  prove by induction on $k$ that there exists a set of bounds
  $B\subset ]0,+\infty[^{n-1}$ such that, for $i= 1,\ldots, k$ and $j=
  1,\ldots, n$, the morphisms $a_{ij}$ and $\alpha_{ij}$ also are
  isomorphisms. For $k=1$, we set for short $K^j =L^{1*}\otimes L^j$
  and we note that
  \begin{gather*}
    \RHom(L^1_{T_1},L^j_{Y_j}) \simeq \rsect_{T_1}(Y;K^j_{Y_j}),
    \quad\qquad
    \RHom(L^1_{T_1},L^1_{T_1}) \simeq \rsect(Y;K^1_{T_1}), \\
    \RHom(L^1_{T_1},L^1_{Y_1}) \simeq \rsect(Y;K^1_{Y_1}), \quad
    \RHom(L^1_{\Int(T_1)},L^j_{Y_j}) \simeq
    \rsect(\Int(T_1);K^j_{Y_j}).
  \end{gather*}
  Hence by lemma~\ref{lemKS}, we may choose $\varepsilon_1^0$ such
  that $\forall 0<\varepsilon_1 < \varepsilon_1^0$, $a_{11}$ is an
  isomorphism and
  \begin{equation}
    \label{section_int}
    \RHom(L^1_{\Int(T_1)},L^j_{Y_j}) \simeq \rsect(Y_1;K^j_{Y_j}).
  \end{equation}
  Let us prove that we may also find a set of bounds $B'\subset
  ]0,+\infty[^{n-2}$ such that for $(\varepsilon_1,\ldots,
  \varepsilon_{n-1})\in \{\varepsilon_1\} \times B'$ the $a_{1j}$ are
  isomorphisms. For $j\geq 2$, we have $\rsect(Y_1;K^j_{Y_j})=0$
  because $Y_1\cap Y_j=\emptyset$.  Hence the
  identity~\eqref{section_int} and the distinguished triangle
  $L^1_{\Int(T_1)} \to L^1_{T_1} \to L^1_{\partial T_1} \fintd$ imply
  $\RHom(L^1_{T_1},L^j_{Y_j}) \simeq \RHom(L^1_{\partial
    T_1},L^j_{Y_j})$.  We may also assume since the beginning that
  $\varepsilon_1^0$ is small enough so that $\partial T_1$ is smooth.
  Then lemma~\ref{hypersurf_transv}, applied with $S=\partial T_1$,
  yields a set of bounds $B'\subset ]0,+\infty[^{n-2}$ such that
  $\forall (\varepsilon_2, \ldots,\varepsilon_{n-1}) \in B'$, we have
  isomorphisms:
  $$
  \Rhom(L^1_{\partial T_1},L^j_{T_j}) \simeq K^j_{T_{1j}} [-1],
  \quad \RHom(L^1_{\partial T_1},L^j_{T_j}) \simeq \RHom(L^1_{\partial
    T_1},L^j_{Y_j}).
  $$
  (With the notations of lemma~\ref{hypersurf_transv}, we have $T_j
  = T'_j \cap U_+$ and $T_{1j} = S\cap T'_j$.)  Now we just have to
  note that $\Rhom(L^1_{\Int(T_1)},L^j_{T_j}) =0$, because $T_j \cap
  \Int(T_1) =\emptyset$, and use the same distinguished triangle as
  above to conclude.
  
  \medskip Now let us assume the conclusion is true for $k$.  Let $B$
  be the set of bounds given in the statement and
  $(\varepsilon^0_1,\ldots, \varepsilon^0_{n-1})\in B$. Let us set $T=
  T_1\sqcup \ldots \sqcup T_k$, $U=Y\setminus T$.  Arguing as in the
  case $k=1$ on the open subset $U$, we find a set of bounds
  $B'\subset ]0,+\infty[^{n-k-1}$ such that $\forall
  (\varepsilon_1,\ldots, \varepsilon_{n-1})\in B\cap
  (\{(\varepsilon^0_1,\ldots, \varepsilon^0_{k})\} \times B')$ we
  have, with the ``new'' $T_i$ (but the $T_i$ for $i\leq k$ stay the
  same):
  \begin{gather*}
    \Rhom(L^{k+1}_{T_{k+1}}|_U,L^j_{T_j}|_U) \simeq (L^{k+1 *}\otimes
    L^j)_{T'_{k+1,j}}[-1]  \quad \text{for $j> k+1$},  \\
    \RHom(L^{k+1}_{T_{k+1}}|_U,L^j_{T_j}|_U) \isoto
    \RHom(L^{k+1}_{T_{k+1}}|_U,L^j_{Y_j}|_U) \quad \text{for $j\geq
      k+1$},
  \end{gather*}
  where $T'_{k+1,j} = U\cap \partial T_{k+1} \cap (\ovl{T_j} \setminus
  \bigcup_{i<l<j} \ovl{T_l})$. Since $T_{k+1} \subset U$, we have in
  fact $\RHom(L^{k+1}_{T_{k+1}}|_U,L^j_{T_j}|_U) =
  \RHom(L^{k+1}_{T_{k+1}},L^j_{T_j})$ and the same with $Y_j$ instead
  of $T_j$.  Hence $a_{k+1,j}$ is an isomorphism, for $j\geq k+1$.
  Since we have chosen $(\varepsilon_1,\ldots, \varepsilon_{n-1}) \in
  B$, $a_{i,j}$ is an isomorphism for $i\leq k$ and any $j$. It just
  remains to check the case $i=k+1$ and $j\leq k$. But $T_{k+1}$ is
  contained in the open set $U$ which does not meet $T_j$ and $Y_j$.
  Hence $\RHom(L^{k+1}_{T_{k+1}},L^j_{T_j}) =
  \RHom(L^{k+1}_{T_{k+1}},L^j_{Y_j}) =0$ and the assertion is
  trivially true.
    
  Now let us consider $\alpha_{k+1,j}$. Let $u:U\to Y$ be the
  embedding.  Since $T_{k+1} \subset U$, we have:
  $$
  \Rhom(L^{k+1}_{T_{k+1}},L^j_{T_j}) \simeq Ru_*
  \Rhom(L^{k+1}_{T_{k+1}}|_U,L^j_{T_j}|_U) \simeq Ru_*(L^{k+1
    *}\otimes L^j)_{T'_{k+1,j}} [-1].
  $$
  We may assume moreover, up to shrinking $B$, that the conclusions
  of lemma~\ref{tubes_fermes} hold. Hence we have local coordinates
  $(x_i)$ such that $U = \{x_i>0;\: i=1,\ldots,m\}$ and $C=T'_{k+1,j}$
  is a convex locally closed subset of $U$.  It follows that $Ru_*\C_C
  \simeq \C_D$ for some locally closed subset $D\subset Y$. More
  precisely, if $C$ is closed in $U$, we check that $(Ru_*\C_C)_x=\C$
  if $x\in \ovl{C}\smash{{}^Y}$ and is $0$ else.  Hence $D =
  \ovl{C}\smash{{}^Y}$. If we write $C = \ovl{C}\smash{{}^U} \setminus
  C_1$, with $C_1$ closed in $U$, we obtain, by the distinguished
  triangle $\C_C \to \C_{\ovl{C}\smash{{}^U}} \to \C_{C_1} \fintd$,
  that $D = \ovl{C}\smash{{}^Y} \setminus \ovl{C_1}\smash{{}^Y}$.
  When applied to $T'_{k+1,j}= \partial T_{k+1} \cap (\ovl{T_j}
  \setminus \bigcup_{i<l<j} \ovl{T_l}) = (T_{k+1} \cap \ovl{T_j})
  \setminus (\bigcup_{i<l<j} \ovl{T_l})$, this formula for $D$ gives
  the expression of the proposition.
  
  \smallskip
  
  (iii) We first consider $p<i<j$, fixed. We set $F=
  \smash{\Rhom(L^i_{Y_i}, L^j_{Y_j})}$ and $G= \smash{\Rhom(L^i_{T_i},
    L^j_{T_j})}$. By lemma~\ref{bons_tubes} applied to $F$, there
  exists a set of bounds, $B_0$, such that, $\forall
  (\udl{\varepsilon}) \in B_0$, we have $\rsect(Y;F_{T_p}) \simeq
  \rsect(Y;F_{Y_p})$.  Let us also consider the Verdier dual $DF =
  \Rhom(F,\C_Y[d_Y])$ (where $d_Y$ is the real dimension of $Y$). In
  the proof of lemma~\ref{bons_tubes}, at the $p^{th}$ step, when
  $(\varepsilon_1,\ldots,\varepsilon_{p-1})$ are fixed, we can choose
  $\varepsilon_p$ so as to have moreover $\rsect(\Int(T_p); DF) \simeq
  \rsect(T_p \cap Y_p; DF)$ (for this we apply lemma~\ref{lemKS} to
  $DF$ restricted to $Y\setminus \bigsqcup_{k<p} T_k$). We also set
  $S_p = T_p \setminus \Int(T_p)$.  Let us prove that
\begin{equation}
  \label{eq:iiietape0}
  \rsect(Y;F_{Y_p}) \simeq \rsect(\ovl{S_p};F_{S_p}).
\end{equation}
We note that $\rsect_{Y_p}(F)\simeq \Rhom(L^i_{Y_i \cap
  Y_p},L^j_{Y_j}) =0$. Since $F$ is constructible we get $(DF)_{Y_p}
\simeq D(\rsect_{Y_p}(F)) =0$, so that $\rsect(\Int(T_p); DF) \simeq
\rsect(T_p \cap Y_p; DF) =0$.  By Poincar\'e-Verdier duality, this
gives the vanishing of the cohomology with compact supports,
$\rsect_c(\Int(T_p); F)=0$. Since $\ovl{T_p}$ is compact, this is
equivalent to $\rsect(Y;F_{\Int(T_p)})=0$.  By the distinguished
triangle $F_{\Int(T_p)} \to F_{T_p} \to F_{S_p} \fintd$, we deduce
$\rsect(Y; F_{T_P}) \simeq \rsect(Y;F_{S_p})$.  This
yields~\eqref{eq:iiietape0} because $\rsect(Y; F_{T_P}) \simeq
\rsect(Y; F_{Y_P})$ and $\rsect(Y;F_{S_p}) \simeq
\rsect(\ovl{S_p};F_{S_p})$.
  
On the other hand, $\supp G \subset \ovl{T_i} \cap \ovl{T_j}$ and,
for $k>p$ we have $\ovl{T_k} \cap T_p \subset S_p$,
so that $\rsect(Y;G_{T_p}) \simeq \rsect(\ovl{S_p};G_{S_p})$.
Together with~\eqref{eq:iiietape0}, this shows that (iii) for our
$p,i,j$, is equivalent to
  \begin{equation}
    \label{eq:iiietape1}
    \rsect(\ovl{S_p};F_{S_p}) \simeq \rsect(\ovl{S_p}; G_{S_p}).
  \end{equation}
  We let $\S_p$ be the smallest set of subsets of $\ovl{S_p}$ stable
  by taking intersections and complements and containing the
  $\ovl{S_p} \cap T_k$, for $k<p$. This set $\S_p$ is finite and we
  denote by $\S_{p,min}$ the set of its minimal elements (for the
  inclusion relation); then any element of $\S_p$ is a union of
  elements of $\S_{p,min}$.  By lemma~\ref{tubes_fermes} (iii), up to
  shrinking the set of bounds $B_0$, any $S\in \S_{p,min}$ is a locally
  closed submanifold of $Y$, transversal to every stratum $Y_i$ that
  it meets.  Let us first prove:
 \begin{equation}
    \label{eq:iiietape2}
    \forall S\in \S_{p,min}, \qquad
 \rsect(S;F) \simeq \rsect(S; G).
  \end{equation}
  By part (i) of the proposition, we have $G|_{S} \simeq (L^{i*}\otimes
  L^j)_{T_{ij}\cap S} [-1]$.  Let us stratify $S$ by $S =
  \bigsqcup_{k>p}(S \cap Y_k)$ and define $T^S_k$, $T^S_{kl}$,
  similarly as $T_k$, $T_{kl}$, with the functions $f_m|_{S}$. We have
  $\ovl{T^S_k} = \ovl{T_k} \cap S$. Hence (i) applied to $S$ gives
  $$
  \Rhom(L^i_{T^S_i},L^j_{T^S_j}) \simeq (L^{i*}\otimes L^j)_{T^S_{ij}}
  [-1] \simeq (L^{i*}\otimes L^j)_{T_{ij}}|_{S} [-1] \simeq
  G|_{S}.
  $$
  Now, part (ii) of the proposition applied to $S$ gives, up to
  shrinking $B_0$ again:
  $$
  \rsect(S; G) \simeq \RHom(L^i_{T^S_i},L^j_{T^S_j}) \simeq
  \RHom(L^i_{Y_i \cap S},L^j_{Y_j \cap S}) \simeq \rsect(S;F),
  $$
  where the last isomorphism follows from the transversality of $S$
  and $Y_i$, $Y_j$. This is~\eqref{eq:iiietape2}. We will deduce:
\begin{equation}
  \label{eq:iiietape3}
  \forall V\in \S_p, \: \text{$V$ closed in $\ovl{S_p}$}, \qquad
\rsect(V;F) \simeq \rsect(V;G).
\end{equation}
Since $F$ is constructible with respect to the stratification
$\bigsqcup Y_i$ and any $S\in \S_{p,min}$ is transversal to every
stratum $Y_i$, we have the isomorphisms, $\forall S\in\S_{p,min}$,
$u_S: \rsect(\ovl{S};F) \isoto \rsect(S;F)$.
  
Recall that $G \simeq (L^{i*}\otimes L^j)_{T_{ij}} [-1]$. For any
$S\in \S_{p,min}$, the inclusion $(T_{ij} \cap S) \subset T_{ij}\cap
\ovl{S}$ is an equivalence of homotopy (by lemma~\ref{tubes_fermes}
(iii)), so that we also have an isomorphism $v_S: \rsect(\ovl{S};G)
\isoto \rsect(S;G)$.
  
Now let us prove~\eqref{eq:iiietape3}. Let us write $V = V_1
\sqcup\ldots \sqcup V_r$, with $V_i \in \S_{p,min}$, and argue by
induction on $r$. For $r=1$, our assertion is~\eqref{eq:iiietape2}. We
may assume that $V_r$ is of maximal dimension (among the $V_i$), so
that $V'= V \setminus V_r$ is closed and the induction hypothesis
applies to $V'$.  We have $\ovl{V_r} = V$ or $V'' = \ovl{V_r} \cap V'$
is a closed union of less than $r$ subsets $V_i$ and the induction
hypothesis also applies to $V''$.  By~\eqref{eq:iiietape2} and the
isomorphisms $u_{V_r}$, $v_{V_r}$, we have $\rsect(\ovl{V_r};F) \simeq
\rsect(\ovl{V_r};G)$.  We conclude by the Mayer-Vietoris distinguished
triangle
$$
\rsect(V;F) \to \rsect(\ovl{V_r};F) \oplus \rsect(V';F) \to
\rsect(V'';F) \fintd,
$$
and the similar one for $G$. 

Now we can prove~\eqref{eq:iiietape1} (and thus (iii) for our
$i,j,p$). We have an excision distinguished triangle,
$$
\rsect(\ovl{S_p};F_{S_p}) \to \rsect(\ovl{S_p};F) \to
\rsect(\ovl{S_p} \setminus S_p;F) \fintd,
$$
and a similar one for $G$.
The last two terms of these distinguished triangles are isomorphic
because of~\eqref{eq:iiietape3} applied to $V=\ovl{S_p}$ and
$V= \ovl{S_p} \setminus S_p$. Hence the first terms are also
isomorphic, as desired.

\smallskip
  
For $i=j$, the same proof works, replacing the isomorphism $G \simeq
(L^{i*}\otimes L^j)_{T_{ij}} [-1]$ by $G \simeq (L^{i*}\otimes
L^i)_{\ovl{T_i}}$: we still have $\smash{\Rhom(L^i_{T^S_i},L^i_{T^S_i})
\simeq G|_{S}}$ and $(\ovl{T_{i}} \cap S) \subset \ovl{T_{i}}\cap
\ovl{S}$ is an equivalence of homotopy.  For $j< i$, we have $F =
\Rhom(L^i_{Y_i}, L^j_{Y_j}) = 0$ because $Y_i$ is open in $\ovl{Y_i}$
and $Y_j \subset (\smash{\ovl{Y_i}} \setminus Y_i)$. In the same way
$G=0$ and (iii) is trivial.  

\smallskip

We let $B_{i,j,p}$ be the subset of $B_0$ formed by the
$\udl\varepsilon$ such that (iii) holds for $p,i,j$. This is a set of
bounds and the intersection of all $B_{i,j,p}$, for $i,j >p$, gives us
the desired set of bounds.

\end{proof}

\begin{proof}[Proof of proposition~\ref{prop:plomberie}]
  Let us set, for $i=1,\ldots,n$, $k= 1,2$, $L^{ki} = \Ll^k$ if $Y_i
  \subset Z_k$ and $L^{ki} =0$ else. We set $L^i = L^{1i} \oplus
  L^{2i}$; this is a local system defined in a neighbourhood of $Y_i$.
  We first choose a set of bounds $B$ such that the conclusions of
  lemma~\ref{tubes_fermes} and proposition~\ref{prop_tubes} hold.
  
\medskip

  (i) Let us prove~\eqref{eq:iso_LT_LZ}. We begin with the case where
  $Z_1 =Y_i$ is a single stratum.  Let $Z' \subset {Z_2}$ be a closed
  subset of $Z_2$ (which is a union of strata). Set $W= {Z_2}\setminus
  Z'$. We have $T_{Z_2} = T_{Z'} \sqcup T_W$ and $T_{Z'}$ is closed in
  $T_{Z_2}$. This gives a morphism of distinguished triangles: $
  \def\objectstyle{\scriptstyle} \def\labelstyle{\scriptscriptstyle}
  \vcenter{\xymatrix@R=3mm{ \C_{T_W} \ar[d]\ar[r] & \C_{T_{Z_2}}
      \ar[d]\ar[r] &   \C_{T_{Z'}}\ar[d]\ar[r]^{\smash{+1}} &  \\
      \C_W \ar[r] & \C_{Z_2} \ar[r] & \C_{Z'}\ar[r]^{+1} & }} $.
  Tensoring this diagram by $\Ll^2$ and applying the functor
  $\RHom(L^i_{T_i},\cdot)$, where $Y_i$ is any stratum of $Y$, we
  obtain the morphism of distinguished triangles:
  $$
  \xymatrix@R=5mm{ \RHom(L^i_{T_i},\Ll^2_{T_W}) \ar[d]\ar[r] &
    \RHom(L^i_{T_i},\Ll^2_{T_{Z_2}}) \ar[d]^{f_i} \ar[r]
    &   \RHom(L^i_{T_i},\Ll^2_{T_{Z'}})  \ar[d]\ar[r]^-{+1} &  \\
    \RHom(L^i_{T_i},\Ll^2_W) \ar[r] & \RHom(L^i_{T_i},\Ll^2_{Z_2})
    \ar[r] & \RHom(L^i_{T_i},\Ll^2_{Z'}\ar[r]^-{+1}) &
    \ponctuation{.} }
  $$
  It follows that the morphism $f_i$ in this diagram is an
  isomorphism, for any $i=1,\ldots,n$.  Indeed this is true if ${Z_2}$
  consists of one stratum, say $Y_j$, by proposition~\ref{prop_tubes}:
  we have here $L^j = L^{1j} \oplus L^{2j}$, and $f_i$ is the
  $L^{2j}$-component of isomorphism $a_{ij}$ of this proposition. Then
  the above diagram allows an induction on the number of strata of
  $Z_2$.
  
  The same reasoning, on the first argument of $\RHom(\cdot,\cdot)$,
  gives the similar isomorphism with $L^i_{T_i}$ replaced by
  $\Ll^1_{T_{Z_1}}$.  The same proof, using the isomorphisms $b_{ij}$
  of proposition~\ref{prop_tubes}, yields the second isomorphism
  of~\eqref{eq:iso_LT_LZ}.
    
  \medskip
  
  (ii) We prove~\eqref{eq:iso_LT_LZ_bord} first in the case where
  $Z_1=Y_i$ and $Z_2=Y_j$.  Let us set $F= \Rhom(\Ll^1_{Y_i},
  \Ll^2_{Y_j})$, $G= \Rhom(\Ll^1_{T_i}, \Ll^2_{T_j})$ and $Z'=Z\cap
  \ovl{Y_i}\cap \ovl{Y_j}$.  We see that $F_Z = F_{Z'}$ and $G_{T_Z} =
  G_{T_{Z'}}$, hence we may assume that $Z$ consists of strata $Y_p$
  with $p<i,j$ (recall that $Z_1$, $Z_2\subset Y\setminus \ovl{Z}$, so
  that $p\not= i$, $p\not= j$). 
  
  Let us argue by induction on the number of strata of $Z$.  By
  proposition~\ref{prop_tubes}, \eqref{eq:iso_LT_LZ_bord} is true when
  $Z=Y_p$ with $p<i,j$.  Let us choose $Y_p\subset Z$ such that it is
  open in $Z$.  We have an excision distinguished triangle
  $\rsect(Y;F_{Y_p}) \to \rsect(Y;F_Z) \to \rsect(Y;F_{Z\setminus
    Y_p})\fintd$, and a similar one with $G_{T_Z}$. By induction
  $\rsect(Y;F_{Z\setminus Y_p}) \simeq \rsect(Y;G_{T_{Z\setminus
      Y_p}})$, and by the first step $\rsect(Y;F_{Y_p}) \simeq
  \rsect(Y;G_{T_p})$.  Hence we obtain $\rsect(Y;F_Z) \simeq
  \rsect(Y;G_{T_Z})$ as desired.
  
  Going from a single stratum to arbitrary locally closed sets $Z_1$,
  $Z_2$, is the same as in the proof of~\eqref{eq:iso_LT_LZ}.
  
  \medskip
  
  (iii) Let us first assume that $Z_1$ consists of a single stratum.
  Since the statement is local on $Y$, we may take coordinates as in
  lemma~\ref{tubes_fermes} and assume $\Ll^1 =\Ll^2 = \C_{\R^d}$. We
  set $C_l= \{x_1>0,\ldots, x_{l-1} >0, \, x_l \leq 0 \}$ and assume
  that $Z_1 = C_i$.  We set $U=\{ x_1>0,\ldots, x_{i-1} >0 \}$ and let
  $u:U\to \R^d$ be the inclusion.  Since $Z_2$ is locally closed and
  $Z_1 \subset Z_2$, there exists $l\geq i$ such that $U\cap T_{Z_2} =
  C_i \sqcup \ldots \sqcup C_l$. Since $C_i \subset U$, we have:
  $$
  \Rhom(\C_{C_i},\C_{T_{Z_2}}) \simeq Ru_*\Rhom(\C_{C_i}
  ,\C_{T_{Z_2}}|_U ) \simeq Ru_* (\C_{C_i \setminus C_{il}}),
  $$
  where $C_{il} = \ovl{C_{i+1} \sqcup \ldots \sqcup C_l}$.  Now
  $Ru_*(\C_{C_i \setminus C_{il}}) \simeq \C_{\ovl{C_i} \setminus
    \ovl{C_{il}}}$ is concentrated in degree $0$.
  
  We deduce the result for an arbitrary locally closed subset
  $Z_1\subset Z_2$ using an excision distinguished triangle.

   \medskip
 
  (iv) We first note that we may assume $V=Y$. Indeed, by
  lemma~\ref{bons_tubes} applied to the complex $G =
  \Rhom(\Ll^1_{Z_1},\Ll^2_{Z_2})$, we may choose a set of bounds such
  that $\RHom(\Ll^1_{Z_1}|_V,\Ll^2_{Z_2}|_V) \simeq
  \RHom(\Ll^1_{Z_1}|_{T_V}, \Ll^2_{Z_2}|_{T_V})$. Hence,
  if~\eqref{eq:iso_LT_LZ_degre0} is true for global sections, applying
  it to $Y=T_V$, with the induced stratification and local systems,
  gives the result for any open $V$.
  
  Let us set $F=\hom(\Ll^1_{T_{Z_1}},\Ll^2_{T_{Z_2}})$, $U =
  Y\setminus (\ovl{Z_1} \setminus Z_1)$ and show that we may replace
  $Y$ by $T_U$.  Since $T_{Z_1} \subset T_U$, we have $F\simeq
  \sect_{T_U}(F)$. But $F|_{T_U} \simeq (\Ll^{1*} \otimes
  \Ll^2)_{T_{Z_1\cap Z_2}}$ and locally around any point of
  $\ovl{T_U}$, the inclusions $T_{Z_1\cap Z_2} \subset T_U \subset Y$
  are homeomorphic to inclusions of convex subsets of $\R^d$. Hence
  $\sect_{T_U}(F) \simeq \rsect_{T_U}(F)$ and $\rsect(Y;F) \simeq
  \rsect(T_U;F)$. We also have $G\simeq \rsect_U(G)$, hence
  $\RHom(\Ll^1_{Z_1},\Ll^2_{Z_2}) \simeq \rsect(U;G)$. By
  lemma~\ref{bons_tubes} applied to $G$, we have, up to shrinking the
  set of bounds, $\rsect(U;G) \simeq \rsect(T_U;G) \simeq
  \RHom(\Ll^1_{Z_1}|_{T_U}, \Ll^2_{Z_2}|_{T_U})$.  Hence we may
  replace $Y$ by $T_U$ and assume $Z_1$ is closed.
  
  By assertion (iii) proved above, we have, setting
  $F'=\Rhom(\Ll^1_{T_{Z_1}},\Ll^2_{T_{Z_2}})$, $F'_{T_{Z_2}} \simeq
  F_{T_{Z_2}}$. Since $Z_1$ is closed we also have an exact sequence
  $0 \to \sect_{T_{Z_1}} \Ll^2_{T_{Z_2}} \to \Ll^2_{T_{Z_2}}$, so that
  $F_x =0$ for $x\not\in T_{Z_2}$, and $F_{T_{Z_2}} \simeq F$.  Hence
  $\rsect(Y;F) \simeq \rsect(Y;F'_{T_{Z_2}})$ and, setting
  $W=\ovl{Z_2} \setminus Z_2$ (which is closed), we have the
  distinguished triangles:
  \begin{gather*}
    \rsect(Y;F'_{T_W}) \to \rsect(Y;F') \to  \rsect(Y;F) \fintd, \\
    \Rhom(\Ll^1_{T_{Z_1\cap W}}, \Ll^2_{T_{Z_2}}) \to F' \to
    \Rhom(\Ll^1_{T_{Z_1\cap (Y\setminus W)}}, \Ll^2_{T_{Z_2}}) \fintd.
  \end{gather*}
  By~\eqref{eq:iso_LT_LZ}, $\rsect(Y;F') \simeq
  \RHom(\Ll^1_{Z_1},\Ll^2_{Z_2})$, thus the first triangle implies
  that~\eqref{eq:iso_LT_LZ_degre0} is equivalent to the vanishing of
  $\rsect(Y;F'_{T_W})$. Using the second triangle, this vanishing
  follows from the two sequences of isomorphisms below:
\begin{align}
  \label{van1}
  \rsect(Y;\Rhom(\Ll^1_{T_{Z_1\cap W}}, \Ll^2_{T_{Z_2}})_{T_W})
  &\simeq \RHom(\Ll^1_{T_{Z_1\cap W}}, \Ll^2_{T_{Z_2}})  \\
  \label{van2}
  &\simeq \RHom(\Ll^1_{Z_1\cap W}, \Ll^2_{Z_2})  \\
  \label{van3}
  &\simeq \RHom(\Ll^1_{Z_1\cap W}, \rsect_{Z_2}(\Ll^2_{Z_2})) = 0,
\end{align}
\begin{equation}
  \begin{split}
  \label{van4}
  \rsect(Y;\Rhom(\Ll^1_{T_{Z_1\cap (Y\setminus W)}},&
  \Ll^2_{T_{Z_2}})_{T_W})  \\
  &\simeq \rsect(Y;\Rhom(\Ll^1_{Z_1\cap (Y\setminus W)},
  \Ll^2_{Z_2})_{W}) =0.
  \end{split}
\end{equation}
Let us explain these isomorphisms: \eqref{van1} is true because the
application of the functor $(\cdot)_{T_W}$ does not change anything
(since the support of the complex of sheaves is included in $T_W$),
\eqref{van2} follows from~\eqref{eq:iso_LT_LZ}, \eqref{van3} comes
from the hypothesis $R\mu_*(\Ll^2|_{Z_2}) = \mu_!(\Ll^2|_{Z_2})$ and
the vanishing follows from $(Z_1\cap W) \cap Z_2 =\emptyset$.

The first isomorphism in~\eqref{van4} follows
from~\eqref{eq:iso_LT_LZ_bord}.  The smoothness hypothesis gives
$\Rhom(\Ll^1_{Z_1\cap (Y\setminus W)}, \Ll^2_{Z_2}) \simeq
R\nu_*((\Ll^{1*} \otimes \Ll^2)|_{Z_1\cap Z_2})$.  Since $\Ll^1$ is
defined on $Z_1$, which is closed, we also have $R\nu_*((\Ll^{1*}
\otimes \Ll^2)|_{Z_1\cap Z_2}) \simeq \nu_!  ((\Ll^{1*} \otimes
\Ll^2)|_{Z_1\cap Z_2})$; this implies the desired vanishing and
concludes the proof of~\eqref{eq:iso_LT_LZ_degre0}.
\end{proof}

\subsection{Direct image to $I$}
We consider a stratified analytic manifold $Y=\bigsqcup_{i \in I} Y_i$
endowed with normal crossings divisors $(D_v)_{v\in V}$ and sheaves
$(L_\aaa)_{\aaa\in \epm}$, satisfying
assumptions~\ref{assumptionstratification}.  As in
notations~\ref{notationsystemdeetubes}, we fix $\udl\varepsilon$ in a
set of bounds such that the conclusions of
proposition~\ref{prop:plomberie} hold for any pair $(\aaa,\bbb) \in
\epm^2$.  We keep the notations $T_i$, for $i\in I$, and $T_\aaa$,
$L'_\aaa$, for $\aaa\in \epm$.  We will give an equivalence between
$\D(Y)\langle L'_\aaa \rangle$ and a derived category of dg-modules on
$I$ following the construction of~\cite{L95}.

\begin{definition}
\label{def:faisceauA}
We define a sheaf of dg-algebras on $Y$, $\Omega$, by
$$
\Omega = \oplus_{(\aaa,\bbb) \in \epm^2} \Omega_{\aaa\bbb}, \quad
\text{where, for $\aaa,\, \bbb\in \epm$,} \quad \Omega_{\aaa\bbb}
= \Omega_Y \otimes \hom(L'_\aaa, L'_{\bbb}).
$$
The product $m_{\aaa\bbb\ccc} :\Omega_{\bbb\ccc} \otimes
\Omega_{\aaa\bbb} \to \Omega_{\aaa\ccc}$, for $\aaa, \bbb, \ccc \in
\epm$, is induced by the product of forms and the composition in the
endomorphisms sheaves.

Considering the partition $Y=\bigsqcup_{i\in I} T_i$, we define
$\pxf':Y\to I$ (like $\pxf: Y=\bigsqcup_{i \in I} Y_i \to I$) by
$\pxf'(T_i) = \{i\}$.  We set :
$$
\fB^{\aaa\bbb} = \pxf'_*(\Omega_{\aaa\bbb}), \qquad \fB=
\pxf'_*(\Omega) \simeq \oplus_{(\aaa,\bbb) \in \epm^2}
\fB^{\aaa\bbb} .
$$
These are sheaves of differential graded $\C$-vector spaces on $I$,
and $\fB$ is a sheaf of dg-algebras.  We also define a direct image
functor $\gamma:\D^+(Y) \to \D_\fB$.  For $F\in \D^+(Y)$, we choose an
injective resolution of $F$, say $F\to R_F$ ($R_F$ and the morphism
$F\to R_F$ depending functorially on $F$), and we set:
$$
\gamma(F) = \pxf'_*(\oplus_{\aaa\in \epm}\, \Omega_Y \otimes
\hom(L'_\aaa,R_F)).
$$
\end{definition}
We note that $\gamma(F)$ has a natural structure of $\fB$-dg-module,
defined by multiplication of forms and composition of homomorphisms
sheaves, like the multiplicative structure of $\fB$.

We also have the following description of the cohomology of the
sections of $\fB$. Since $\Omega_Y$ is a soft resolution of $\C_Y$,
$\Omega_{\aaa\bbb}$ is a soft resolution of $\hom(L'_\aaa,
L'_{\bbb})$.  Let $U\subset I$ be an open subset, $V = \pxf^{-1}(U)$,
$V'={\pxf'}^{-1}(U) = T_V$, and let $\aaa, \bbb\in \epm$.  By (iii)
and (iv) of proposition~\ref{prop:plomberie}, we obtain:
\begin{equation}
  \label{eq:section_de_A}
  H^i(\pxf'_*(\Omega_{\aaa\bbb})(U)) \simeq
  \Ext^i((L'_\aaa)|_{V'}, (L'_{\bbb})|_{V'}) \simeq
  \Ext^i((L_\aaa)|_{V}, (L_{\bbb})|_{V}).
\end{equation}
\begin{proposition}
  \label{prop:eqL'M}
  With the above notations, we set $M_\aaa = \gamma(L'_\aaa)$, for
  $\aaa\in \epm$. The functor $\gamma$ induces an equivalence of
  categories between $\D(Y)\langle L'_\aaa \rangle$ and
  $\D_{\fB}\langle M_\aaa \rangle$.
\end{proposition}
\begin{proof}
  Since our categories are respectively generated by the $L'_\aaa$
  and the $M_\aaa$, it is sufficient to prove that $\gamma$ gives
  isomorphisms, $\forall \aaa, \bbb\in \epm$, $\forall p\in\Z$:
  \begin{equation}
  \label{eq:bij_morphismes}
  \Hom_{\D(Y)}(L'_\aaa,L'_{\bbb}[p]) \simeq 
  \Hom_{\D_\fB}(M_\aaa,M_{\bbb}[p]).
  \end{equation}
  Indeed an inductive argument as in the proof of
  corollary~\eqref{equiv_LetL'} (but easier because, here, the functor
  giving the equivalence is a priori defined) implies that $\gamma$
  also gives a bijection between $\Hom_{\D(Y)}(L_1,L_2)$ and
  $\Hom_{\D_\fB}(\gamma(L_1),\gamma(L_2))$ for any objects $L_1$,
  $L_2$ of $\D(Y)\langle L'_\aaa \rangle$. Let us
  prove~\eqref{eq:bij_morphismes}. We set $L'=\oplus_{\aaa\in \epm}
  L'_\aaa$ and $M=\gamma(L') = \oplus_{\aaa\in \epm} M_\aaa$.
  Then~\eqref{eq:bij_morphismes} is equivalent to: 
  \begin{equation}
    \label{eq:iso_Lprime_M}
    \forall p\in\Z \qquad
\Hom_{\D(Y)}(L',L'[p]) \simeq \Hom_{\D_\fB}(M,M[p]) .
  \end{equation}
  For $\aaa\in \epm$, let $L'_\aaa \to R_\aaa$ be the chosen
  injective resolution of $L'_\aaa$. It induces a morphism of
  differential graded sheaves, $f^{\aaa\bbb}$, from
  $\fB^{\aaa\bbb} = \pxf'_*(\Omega_{\aaa\bbb})$ to
  $\fB'^{\aaa\bbb} = \pxf'_*(\Omega_Y \otimes \hom(L'_\aaa,
  R_{\bbb}) )$.  Since $R_\bbb$ is injective we have the
  isomorphisms below in $\D(Y)$, which give the cohomology of sections
  of $\fB'^{\aaa\bbb}$, for an open set $U\subset I$, and
  $V'=\pxf'^{-1}(U)$:
  \begin{gather*}
    \Omega_Y \otimes \hom(L'_\aaa, R_{\bbb}) \simeq
    \hom(L'_\aaa, R_{\bbb}) \simeq
    \Rhom(L'_\aaa, R_{\bbb}),  \\
    H^\cdot(\fB'^{\aaa\bbb}(U)) = H^\cdot(V';\Omega_Y \otimes
    \hom(L'_\aaa, R_{\bbb})) \simeq \Ext^i((L'_\aaa)|_{V'},
    (L'_{\bbb})|_{V'}).
  \end{gather*}
  By~\eqref{eq:section_de_A} this means that $f^{\aaa\bbb}$ is a
  quasi-isomorphism of differential graded sheaves.  Summing over all
  pairs $(\aaa,\bbb)$, we obtain a quasi-isomorphism of $\fB$-modules
  between $\fB$ and $\gamma(\oplus_{\bbb\in \epm} L'_\beta )= M$.
  Hence we obtain:
  \begin{equation}
    \label{eq:iso_fA_M}
 \Hom_{\D_\fB}(M,M[p]) \simeq  \Hom_{\D_\fB}(\fB,\fB[p]).
  \end{equation}
  We have seen that $\forall i\in I$, $\fB_{U_i}$ is $K$-projective.
  By (i) of assumptions~\ref{assumptionstratification}, any
  intersection of open sets of the type $U_i$ still is of this kind,
  hence $\forall i_1,\ldots, i_n \in I$, $\fB_{U_{i_1}\cap \ldots \cap
    U_{i_n}}$ is $K$-projective.  Let us put any total order on $I$;
  we obtain a $K$-projective resolution of $\fB$ by taking the total
  complex of the following \v{C}ech-like complex of complexes:
  $$
  \cdots\to \oplus_{i_1<i_2<i_3\in I} \fB_{U_{i_1}\cap U_{i_2}\cap
    U_{i_3}} \to \oplus_{i_1<i_2\in I} \fB_{U_{i_1}\cap U_{i_2}} \to
  \oplus_{i_1\in I} \fB_{U_{i_1}} \to 0,
  $$
  with the usual differential $(da)_{i_1,\ldots,i_r} = \sum_{i_1 <
    \cdots < i_k < j < i_{k+1} <\cdots < i_r} (-1)^k
  a_{i_1,\ldots,j,\ldots,i_r}$.  For an open set $U\subset I$, we have
  $\RHom(\fB_U,\fB) \simeq \rsect(U;\fB)$ and, for $U=U_i$, the
  functor $\sect(U_i;\cdot) = (\cdot)_i$ is exact. Hence:
  $$
  \Hom_{\D_{\fB}}(\fB_{U_i},\fB[p]) \simeq H^p(\RHom(\fB_{U_i},\fB))
  \simeq H^p(\rsect(U_i;\fB)) \simeq H^p(\sect(U_i;\fB)).
  $$
  By definition $\sect(U_i;\fB) = \sect(\pxf'^{-1}(U_i);\Omega)$
  and we obtain that $\Hom_{\D_\fB}(M,M[p])$ is the $p^{\mathrm{th}}$
  cohomology group of the total complex of the double complex:
  $$
  0\to \oplus_{i_1\in I} \sect(\pxf'^{-1}(U_{i_1});\Omega) \to
  \oplus_{i_1<i_2\in I} \sect(\pxf'^{-1}(U_{i_1})\cap
  \pxf'^{-1}(U_{i_2});\Omega)\to\cdots .
  $$
  This is a \v{C}ech resolution of the complex $\Omega$, which is
  formed by soft sheaves. Hence $\Hom_{\D_\fB}(\fB,\fB[p])\simeq
  H^p(Y;\Omega)$, and, by~\eqref{eq:section_de_A}, this last group is
  isomorphic to $\Ext^p(L',L')$.  In view of~\eqref{eq:iso_fA_M}, this
  gives~\eqref{eq:iso_Lprime_M} and concludes the proof.
\end{proof}

\subsection{Gysin isomorphism and product}
In the previous paragraph, we have obtained a sheaf of dg-algebras,
$\fB$, on $I$ such that $\D(Y)\langle L_\aaa \rangle$ is equivalent to
a subcategory of $\D_{\fB}$. In
section~\ref{formality_of_de_Rham_algebra}, we will construct a
sequence of quasi-isomorphisms between $\fB$ and its cohomology.  For
this we will in particular replace sets like $U_{i\aaa\bbb} = \{ x\in
\pxf'^{-1}(U_i);\: \hom(L'_\aaa,L'_{\bbb})_x \not= 0\}$ by homotopy
equivalent ones. But before that, we note that $\fB_i$ is not
immediately quasi-isomorphic to sections of $\Omega_{U_{i\aaa\bbb}}
\otimes \hom(L'_\aaa,L'_{\bbb})$. Indeed $U_{i\aaa\bbb}$ is not closed
in $\pxf'^{-1}(U_i)$.  The cohomology of $\fB_i$ is
$\Ext^\cdot(L'_\aaa|_{\pxf'^{-1}(U_i)} , L'_\bbb|_{\pxf'^{-1}(U_i)}) =
\Ext^\cdot(L_\aaa|_{\pxf^{-1}(U_i)},L_\bbb|_{\pxf^{-1}(U_i)})$, but
the cohomology of $\Omega_{U_{i\aaa\bbb}} \otimes
\hom(L'_\aaa,L'_{\bbb})$ is $\Ext^\cdot(L'_\aaa|_{U_{i\aaa\bbb}} ,
L'_\bbb|_{U_{i\aaa\bbb}}) = \Ext^\cdot(L_\aaa|_{Y_i} ,L_\bbb|_{Y_i}
)$.  In our situation they are isomorphic under a ``twisted'' Gysin
isomorphism.  We build this isomorphism at the level of the de Rham
algebras and describe the algebra structure. Then we use this
description to obtain a quasi-isomorphism between $\fB$ and the sheaf
$\fA$ defined in~\ref{env:defin_fB}.

\subsubsection{Gysin isomorphism}
Let us first consider the usual Gysin isomorphism. If $M$ is an
oriented manifold and $N$ a closed oriented submanifold of codimension
$c$, we have an isomorphism $\rsect_N(\C_M) \simeq \C_N[-c]$. On the
global sections it induces $H^\cdot(N;\C_N) \simeq
H^{\cdot+c}_N(M;\C_M)$. Let us choose, by lemma~\ref{lemKS}, two open
tubular neighbourhoods $U$, $V$ of $N$, such that $\ovl{U} \subset V$
and
$$
H^\cdot(V;\C_V) \isoto H^\cdot(U;\C_U) \isoto H^\cdot(N;\C_N), \quad
H^\cdot_N(M;\C_M) \isoto H^\cdot_{\ovl{U}}(M;\C_M).
$$
If we assume moreover that the boundary of $U$ is smooth, we have
$\rsect_{\ovl{U}}\C_M \simeq \C_U$, so that $H^\cdot_N(M;\C_M) \simeq
H^\cdot(M;\C_U)$. Let $\delta(N,M) \in H^c_N(M;\C_M)$ be the
fundamental class of $N$ in $M$. We may choose a representative
$\xi(N,M) \in \sect(M;\Omega^c_M)$ such that $\supp \xi(N,M) \subset
U$.  Then, the multiplication by $\xi(N,M)$ induces a well-defined
quasi-isomorphism of $\sect(N;\Omega_N)$-dg-modules between
$\sect(V;\Omega_V)$ and $\sect(M;(\Omega_M)_U)[c]$.  More generally,
given local systems $L_1$ on $M$, and $L_2$ on $V$, we obtain
isomorphisms between extension groups and a corresponding quasi-isomorphism
between de Rham complexes:
\begin{gather}
  \label{gysin_extension_1}
  \cdot \delta(N,M): \Ext^\cdot_{\D(N)}(L_2|_N,L_1|_N) \isoto
  \Ext^{\cdot +c}_{\D(M)}((L_2)_{N},L_1),  \\
  \label{gysin_extension_2}
  \Ext^\cdot_{\D(V)}(L_2,L_1|_V) \isoto \Ext^{\cdot
    +c}_{\D(M)}((L_2)_{\ovl{U}},L_1), \\
  \label{gysin_de_Rham}
  \begin{split}
    \cdot \xi(N,M): \sect(V;\Omega_V \otimes & \hom(L_2,L_1|_V)) \\
    & \smash{\xto{qis}} \sect(M; \Omega_M \otimes
    \hom((L_2)_{\ovl{U}},L_1))[c].
  \end{split} 
\end{gather}
We will use~\eqref{gysin_de_Rham} to build a second sheaf of
dg-algebras on $I$, quasi-isomorphic to $\fB$.  For this we also need
to describe the product structure: the composition $\hom(L_2,L_1)
\otimes \hom(L_3,L_2) \to \hom(L_3,L_1)$ induces a product on the
right hand side of the above isomorphisms and we want to understand
the corresponding product on the left hand side.

\subsubsection{Algebra structure}
\label{algebrastructure}
This paragraph mainly has an heuristic purpose, in order to justify
the definition of the product $m$ of the sheaf of dg-algebras $\fA$
introduced in~\ref{env:defin_fB} above.  We consider a complex
manifold $Y$, endowed with normal crossings divisors $D_v$, $v\in V$,
and sheaves $L_\aaa$, $\aaa\in \epm$, satisfying
assumptions~\ref{assumptionstratification}.  We keep the notations
of~\ref{assumptionstratification}, in particular for $Z_\Delta$,
$\Delta\in \S$, $Z_\aaa = Z_{\Delta_\aaa} \setminus \bigcup_{v\in
  \Delta'_\aaa} D_v$, for $\aaa\in \epm$. We set also:
$$
Z_{\aaa\bbb} = Z_\aaa \cap Z_\bbb, \qquad
Z_{\aaa\bbb\ccc} = Z_\aaa \cap Z_\bbb \cap Z_\ccc.
$$
We fix $\aaa, \bbb, \ccc \in \epm$ and, up to restricting ourselves
to an open subset of $Y$, we assume that $Z_\aaa$, $Z_\bbb$, $Z_\ccc$
are closed.

For $v\in V$, $D_v$ has a fundamental class $\delta_v\in
H^2_{D_v}(Y;\C_Y)$. For $\Delta\in \S$ the fundamental class of
$Z_\Delta$ in $Y$ is $\delta_\Delta = \prod_{v\in \Delta} \delta_v$.
It belongs to $H^{2|\Delta|}_{Z_\Delta}(Y;\C_Y)$. For $\Delta'\subset
\Delta$, we have $Z_\Delta \subset Z_{\Delta'}$ and $Z_\Delta$ is the
transversal intersection of $Z_{\Delta'}$ and $Z_{\Delta \setminus
  \Delta'}$. Hence the fundamental class $\delta(Z_\Delta,
Z_{\Delta'}) \in H^{2d}_{Z_\Delta}(Z_{\Delta'};\C_{Z_{\Delta'}})$,
with $d=|\Delta|-|\Delta'|$, is the image of $\delta_{\Delta\setminus
  \Delta'} \in H^{2d}_{Z_{\Delta\setminus \Delta'}}(Y;\C_Y)$ in
$H^{2d}_{Z_\Delta}(Z_{\Delta'};\C_{Z_{\Delta'}})$.  By abuse of
notations we will write $\delta_{\Delta\setminus \Delta'}$ for
$\delta(Z_\Delta, Z_{\Delta'})$.

We define $\varepsilon_{\aaa\bbb} = \delta(Z_{\aaa\bbb}, Z_{\bbb})$.
We have:
$$
\varepsilon_{\aaa\bbb} = \delta_{\Delta_\aaa \setminus
  \Delta_{\bbb}} \in
H^{2d_{\aaa\bbb}}_{Z_{\aaa\bbb}}(Z_{\bbb};\C_{Z_{\bbb}}), \qquad
\text{with} \quad d_{\aaa\bbb} = |\Delta_\aaa \setminus
\Delta_{\bbb}|.
$$
  We remark that $\Ext^\cdot_{\D(Y)}(L_\aaa,L_{\bbb})
\simeq \Ext^\cdot_{\D(Y)}((L_\aaa)_{Z_{\aaa\bbb}}, L_{\bbb})$, so
that, by~\eqref{gysin_extension_1}, multiplication by
$\varepsilon_{\aaa\bbb}$ gives an isomorphism:
\begin{equation}
  \label{eq:isom_gysin_LaLb}
  \cdot \varepsilon_{\aaa\bbb} : 
  \Ext^i_{\D(Z_{\aaa\bbb})}(L_\aaa|_{Z_{\aaa\bbb}},
  L_{\bbb}|_{Z_{\aaa\bbb}}) \isoto
\Ext^{i+2d_{\aaa\bbb}}_{\D(Y)}(L_\aaa,L_{\bbb}).
\end{equation}
From now on we forget the subscripts $\D(Z_\cdot)$ if there is no
ambiguity. We want to understand the product of two extension classes
in $\Ext^\cdot(L_\aaa,L_{\bbb})$ and $\Ext^\cdot(L_{\bbb},L_{\ccc})$
in terms of the corresponding classes on the left hand side of the
above formula (see diagram~\eqref{eq:diagA} below).  
Let us introduce some notations:
\begin{alignat}{3}
  \Delta_1 &= \Delta_{\bbb} \setminus (\Delta_{\aaa} \cup
  \Delta_{\ccc}) , &\qquad d_1 &= |\Delta_1|, &\qquad
  \varepsilon_{\aaa\bbb\ccc} &= \delta(Z_{\aaa\bbb\ccc}, Z_{\aaa\ccc})
  = \delta_{\Delta_1}, \\
  \Delta_2 &= (\Delta_{\aaa} \cap \Delta_{\ccc}) \setminus
  \Delta_{\bbb}, &\qquad d_2 &=|\Delta_2|, &\qquad
  \varepsilon^+_{\aaa\bbb\ccc} &=  \delta_{\Delta_2}.
\end{alignat}
We first restrict the classes in $\Ext^\cdot(L_\aaa|_{Z_{\aaa\bbb}},
L_{\bbb}|_{Z_{\aaa\bbb}})$ and $\Ext^\cdot(L_{\bbb}|_{Z_{\bbb\ccc}},
L_{\ccc}|_{Z_{\bbb\ccc}})$ to $Z_{\aaa\bbb\ccc}$ and make the product,
obtaining a morphism:
$$
r: \Ext^i(L_\aaa|_{Z_{\aaa\bbb}},
L_{\bbb}|_{Z_{\aaa\bbb}}) \otimes
\Ext^j(L_{\bbb}|_{Z_{\bbb\ccc}}, L_{\ccc}|_{Z_{\bbb\ccc}})
\to \Ext^{i+j}(L_\aaa|_{Z_{\aaa\bbb\ccc}},
L_{\ccc}|_{Z_{\aaa\bbb\ccc}}).
$$
Now $L_\aaa$ and $L_{\ccc}$ both restrict to local systems on
$Z_{\aaa\ccc}$ (and $Z_{\aaa\bbb\ccc}$), so that we may identify
$\Ext^\cdot(L_\aaa|_{Z_{\aaa\ccc}}, L_{\ccc}|_{Z_{\aaa\ccc}})$ with
$H^\cdot(Z_{\aaa\ccc};L_\aaa^*\otimes L_{\ccc})$ (and the same on
$Z_{\aaa\bbb\ccc}$). Hence multiplication by
$\varepsilon_{\aaa\bbb\ccc}$ gives a morphism:
$$
\cdot \varepsilon_{\aaa\bbb\ccc}:
\Ext^\cdot(L_\aaa|_{Z_{\aaa\bbb\ccc}},
L_{\ccc}|_{Z_{\aaa\bbb\ccc}}) \to
\Ext^{\cdot+2d_1}(L_\aaa|_{Z_{\aaa\ccc}},
L_{\ccc}|_{Z_{\aaa\ccc}}).
$$
The ``correcting term'' $\varepsilon^+_{\aaa\bbb\ccc}$ is
introduced to have the identity~\eqref{eq:produit_epsilon} below (in
$H^\cdot_{\Z_{\aaa\bbb\ccc}}(Y;\C_Y)$). Multiplication by
$\varepsilon^+_{\aaa\bbb\ccc}$ gives morphism~\eqref{mult_correctif}.
\begin{gather}
    \label{eq:produit_epsilon}
  \varepsilon_{\aaa\bbb} \cdot \varepsilon_{\bbb\ccc} =
  \varepsilon^+_{\aaa\bbb\ccc} \cdot \varepsilon_{\aaa\bbb\ccc}\cdot
  \varepsilon_{\aaa\ccc},  \\
  \label{mult_correctif}
  \cdot \varepsilon^+_{\aaa\bbb\ccc}:
  \Ext^{\cdot}(L_\aaa|_{Z_{\aaa\ccc}}, L_{\ccc}|_{Z_{\aaa\ccc}}) \to
  \Ext^{\cdot+2d_2}(L_\aaa|_{Z_{\aaa\ccc}}, L_{\ccc}|_{Z_{\aaa\ccc}}).
\end{gather}
Because of~\eqref{eq:produit_epsilon}, the composition
of~\eqref{mult_correctif} with $\varepsilon_{\aaa\bbb\ccc} \circ r$
gives the commutative diagram
\begin{equation}
  \label{eq:diagA}
  \makebox[0pt]{$\vcenter{\xymatrix@C=7mm{
        \Ext^i(L_\aaa|_{Z_{\aaa\bbb}}, L_{\bbb}|_{Z_{\aaa\bbb}})
        \otimes \Ext^j(L_{\bbb}|_{Z_{\bbb\ccc}},
        L_{\ccc}|_{Z_{\bbb\ccc}}) \ar[d]^\wr_{\varepsilon_{\aaa\bbb}
          \otimes \varepsilon_{\bbb\ccc}} \ar[r]_-p & \Ext^{k}
        (L_\aaa|_{Z_{\aaa\ccc}}, L_{\ccc}|_{Z_{\aaa\ccc}})
        \ar[d]_\wr^{\varepsilon_{\aaa\ccc}}
        \\
        \Ext^{i+2d_{\aaa\bbb}}_{\D(Y)}(L_\aaa,L_{\bbb}) \otimes
        \Ext^{j+2d_{\bbb\ccc}}_{\D(Y)}(L_{\bbb},L_{\ccc}) \ar[r] &
        \Ext^l_{\D(Y)}(L_{\aaa},L_{\ccc}) \ponctuation{,}  }}$}
\end{equation}
where $k = i+j+ 2(d_1 +d_2)$, $l= i+j+ 2( d_{\aaa\bbb} +
d_{\bbb\ccc})$, $p = (\cdot \varepsilon^+_{\aaa\bbb\ccc}) \circ (\cdot
\varepsilon_{\aaa\bbb\ccc}) \circ r$ and the bottom arrow is the
composition of extension classes.  This diagram gives us the
composition of extension classes, through isomorphism
\eqref{gysin_extension_1}. The definition of the product of $\fA$
in~\ref{env:defin_fB} is copied from the definition of $p$ above.

\subsection{Proof of proposition~\ref{prop:eq_cat1}}
Recall that quasi-isomorphic sheaves of dg-al\-ge\-bras have equivalent
derived categories (see~\cite{L95}, proposition 1.11.2).  Hence, by
corollary~\ref{equiv_LetL'} and proposition~\ref{prop:eqL'M}, the
proof of proposition~\ref{prop:eq_cat1} will be achieved if we show
that $\fB$ and $\fA$ are quasi-isomorphic.

We still consider $Y=\bigsqcup_{i\in I} Y_i$, $(D_v)_{v\in V}$,
$(L_\aaa)_{\aaa\in \epm}$ satisfying
assumptions~\ref{assumptionstratification}. We keep the notation
$L_{\aaa,i}$ for the local system on $\pxf^{-1}(U_i)$ extending
$L_\aaa|_{Y_i}$, and also notations~\ref{notationsystemdeetubes} for
$T_i$, $T_\aaa$, $L'_\aaa$, as well as the notations of
paragraph~\ref{algebrastructure} for $\delta_v$, $\delta_\Delta$,
$\Delta_1$, $\Delta_2$.

We choose representatives, $\xi_v\in \sect(Y;\Omega^2_Y)$, of the
$\delta_v$, such that $\supp \xi_v \subset \Int(T_{D_v})$ (remember
that $T_{D_v} = \bigsqcup_{\{i,\, Y_i \subset D_v\}} T_i$).  For
$\Delta \subset V$ we define $\xi_\Delta = \prod_{v\in \Delta} \xi_v
\in \sect(Y;\Omega^{2d}_Y)$; it is a representative of $\delta_\Delta$
with $\supp \xi_\Delta \subset \Int(T_{Z_\Delta})$.  For $\aaa, \bbb,
\ccc \in \epm$ we define forms $\eta_{\aaa\bbb}$,
$\eta_{\aaa\bbb\ccc}$, representing the classes
$\varepsilon_{\aaa\bbb}$, $(\varepsilon_{\aaa\bbb\ccc} \cdot
\varepsilon^+_{\aaa\bbb\ccc})$:
$$
\eta_{\aaa\bbb} = \xi_{(\Delta_\aaa \setminus \Delta_{\bbb})},
\qquad \eta_{\aaa\bbb\ccc} = \xi_{\Delta_1} \cdot \xi_{\Delta_2}.
$$
The forms $\eta_{\aaa\bbb\ccc}$ were already introduced when we
defined the product of $\fA$.  By lemma~\ref{lem:qisAB} below, the
multiplications by the forms $\eta_{\aaa\bbb}$ give quasi-isomorphisms
of sheaves $g^{\aaa\bbb}:\fA^{\aaa\bbb} \to \fB^{\aaa\bbb}$.  By the
definitions and the identity (similar to~\eqref{eq:produit_epsilon})
$\eta_{\aaa\bbb} \cdot \eta_{\bbb\ccc} = \eta_{\aaa\bbb\ccc} \cdot
\eta_{\aaa\ccc}$, the morphism $g= \oplus g^{\aaa\bbb}:\fA \to \fB$ is
a morphism of sheaves of dg-algebras.

We have thus obtained a quasi-isomorphism between $\fB$ and $\fA$, and
hence an equivalence of categories between $\D_\fB \langle M_\aaa
\rangle$ and $\D_\fA \langle \fA \otimes^L_{\fB}M_\aaa \rangle$.  We
remark that $M_\aaa$ is $\fB$-flat because $\fB \simeq \oplus_{\aaa\in
  \epm} M_\aaa$. It follows that $\fA \otimes^L_{\fB}M_\aaa \simeq
N_\aaa$, for the $\fA$-module $N_\aaa = \oplus_{\aaa' \in
  \epm}\fA^{\aaa'\aaa}$. This concludes the proof of the proposition.

\begin{lemma}
  \label{lem:qisAB}
  Let us set, for $i\in I$, $V_i = \pxf^{-1}(U_i)$, $V'_i =
  \pxf'^{-1}(U_i) = T_{V_i}$.

  (i) For $\aaa,\bbb \in \epm$, $i\in
  \pxf(Z_{\aaa\bbb})$, we have a well-defined morphism of sheaves on
  $V'_i$:
  $$
  f^i_{\aaa\bbb}: \Omega^\cdot_{V'_i} \otimes
  \hom(L_{\aaa,i},L_{\bbb,i}) \to \Omega^{\cdot +
    2d_{\aaa\bbb}}_{V'_i} \otimes \hom(L'_\aaa, L'_\bbb), \quad
  (\sigma\otimes u) \mapsto (\eta_{\aaa\bbb}\sigma)\otimes u,
  $$
  where $\sigma\in \Omega_{V'_i}$, $u\in \hom(L'_\aaa,L'_{\bbb})$.
  On the global sections, it induces a morphism $g^{\aaa\bbb}_i :
  \fA^{\aaa\bbb}_i \to \fB^{\aaa\bbb}_i$.
  
  (ii) The morphisms $g^{\aaa\bbb}_i$, $i\in \pxf(Z_{\aaa\bbb})$,
  extend to a morphism of differential graded sheaves on $I$,
  $g^{\aaa\bbb}: \fA^{\aaa\bbb} \to \fB^{\aaa\bbb}$, which is a
  quasi-isomorphism.
\end{lemma}
\begin{proof}
  Since $U_i$ is open, we have, by lemma~\ref{tubes_fermes} (i),
  $Y\setminus V_i \subset T_{Y\setminus V_i} = Y\setminus V'_i$.  By
  definition of $T_{Y\setminus V_i}$, this implies in fact $Y\setminus
  V_i \subset \Int( Y\setminus V'_i)$, or, as well, $\ovl{V'_i}
  \subset V_i$.

\medskip

(i) We consider $i\in \pxf(Z_{\aaa\bbb})$ and set $T^i_{\aaa\bbb} = ((
T_\aaa \cap T_\bbb) \setminus \smash{\ovl{(T_\bbb \setminus T_\aaa)}}
) \cap V'_i$. We first prove the isomorphism:
\begin{equation}
  \label{eq:homlaaalbbb}
   \hom(L'_\aaa, L'_\bbb)|_{V'_i} \simeq (\hom(L_{\aaa,i},
  L_{\bbb,i}))_{T^i_{\aaa\bbb}}.
\end{equation}
We have by definition $L'_\aaa |_{V'_i} \simeq (L_{\aaa,i})_{T_\aaa
  \cap V'_i}$.  Since $i\in \pxf(Z_{\aaa\bbb})$, $Z_\aaa \cap V_i$ is
closed in $V_i$ and, by lemma~\ref{tubes_fermes} (i), (ii), $T_\aaa
\cap V'_i$ is closed in $V'_i$. The same holds for $\bbb$.
Since~\eqref{eq:homlaaalbbb} can be checked locally, we may assume
$L_{\aaa,i} = L_{\bbb,i} = \C_{V_i}$. Now, for two closed subsets $M$,
$N$ of a manifold $X$, and $i_N$ the inclusion of $N$ in $X$, we have
$$
\hom(\C_M,\C_N) \simeq (i_N)_* \hom(i_N^{-1}\C_M,\C_N) \simeq
(i_N)_! \hom(\C_{M\cap N}, \C_N).
$$
One checks that $\hom(\C_{M\cap N}, \C_N) \simeq \C_U$, where $U$
is the interior of $M\cap N$ in $N$. This gives the formula for
$T^i_{\aaa\bbb}$.

\medskip

We may write as well $T^i_{\aaa\bbb} = (( T_\aaa \cap T_\bbb)
\setminus \smash{\ovl{(T_\bbb \setminus T_{Z_{\Delta_\aaa}})}} ) \cap
V'_i$.  A form $\omega \in \sect(V'_i;\Omega_Y)$ such that $\supp
\omega \cap \smash{\ovl{T_\bbb \setminus T_{Z_{\Delta_\aaa}}}} =
\emptyset$ belongs in fact to $\sect(V'_i;(\Omega_Y)_{V'_i \setminus
  \smash{\ovl{(T_\bbb \setminus T_{Z_{\Delta_\aaa}})}}})$. Since
$T_\aaa \cap T_\bbb$ is closed in $V'_i$, we have a natural morphism
from $(\Omega_Y)_{V'_i \setminus \smash{\ovl{(T_\bbb \setminus
      T_{Z_{\Delta_\aaa}})}}}$ to $(\Omega_Y)_{T^i_{\aaa\bbb}}$, and
$\omega$ induces an element of
$\sect(V'_i;(\Omega_Y)_{T^i_{\aaa\bbb}})$. This condition on the
support is satisfied by $\eta_{\aaa\bbb}$. Indeed, we have $\supp
\eta_{\aaa\bbb} \subset \Int( \smash{T_{Z_{\Delta_\aaa \setminus
      \Delta_{\bbb}}}} )$, hence it is sufficient to check that
$\smash{T_{Z_{\Delta_\aaa \setminus \Delta_{\bbb}}}} \cap T_{Z_\bbb
  \setminus Z_{\Delta_\aaa}} = \emptyset$.  This is equivalent to
$Z_{\Delta_\aaa \setminus \Delta_{\bbb}} \cap (Z_\bbb \setminus
Z_{\Delta_\aaa}) = \emptyset$, which is obvious.  Hence the
multiplication by $\eta_{\aaa\bbb}$ sends $\Omega_{V_i}$ into
$(\Omega_{V_i})_{T^i_{\aaa\bbb}}[2d_{\aaa\bbb}]$.  In view
of~\eqref{eq:homlaaalbbb}, this gives the morphism $f^i_{\aaa\bbb}$.
Now, remember that:
\begin{align*}
  \fB^{\aaa\bbb}_i &= \sect(V'_i; \Omega_Y \otimes \hom(L'_\aaa,
  L'_\bbb)), \\
  \fA^{\aaa\bbb}_i &=\sect(V_i;\Omega_Y \otimes
  \hom(L_{\aaa,i},L_{\bbb,i})) \,[-2d_{\aaa\bbb}].
\end{align*}
Hence the restriction from $V_i$ to $V'_i$, composed with $\sect(V'_i;
f^i_{\aaa\bbb})$ gives $g^{\aaa\bbb}_i$.

\medskip

(ii) Let us first see that the $g^{\aaa\bbb}_i$ extend to $i\not\in
\pxf(Z_{\aaa\bbb})$.  In view of the definition of $\fA$
(see~\ref{env:defin_fB}), for $i\not\in \pxf(Z_{\aaa\bbb}) \sqcup
\pxf(Z_{\Delta_\aaa \cup \Delta_\bbb} \cap (\bigcup_{v\in \Delta'_\aaa
  \cap \Delta'_\bbb} D_v))$, we have $\fA^{\aaa\bbb}_i=0$ and
$g^{\aaa\bbb}_i$ is trivially defined.  So we assume $i \in \pxf(
Z_{\Delta_\aaa \cup \Delta_\bbb} \cap ( \bigcup_{v\in \Delta'_{\aaa}
  \cap\Delta'_{\bbb}} D_v)$.  We let $j$ be such that $V_i \setminus
(\bigcup_{v\in \Delta'_\aaa \cap \Delta'_\bbb} D_v) = V_j$.

Let us first extend~\eqref{eq:homlaaalbbb} to this case:
\begin{equation}
  \label{eq:homlaaalbbb_bis}
\hom(L'_\aaa, L'_\bbb)|_{V'_i} \simeq (\hom(L_{\aaa,i},
L_{\bbb,i}))_{U^i_{\aaa\bbb}},
\end{equation}
where $U^i_{\aaa\bbb} = (\ovl{ ( T_\aaa \cap T_\bbb)} \setminus
\smash{\ovl{(T_\bbb \setminus T_\aaa)}} ) \cap V'_i$.  By definition,
$L'_\aaa|_{T_{D_v}} =0$ for $v\in \Delta'_\aaa$, hence we have
$L'_\aaa|_{V'_i} = (L'_\aaa)_{V'_j}$. Denoting by $u:V'_j \to V'_i$
the inclusion, we deduce, in view of~\eqref{eq:homlaaalbbb}:
\begin{equation}
  \label{eq:homlaaalbbb_ter}
  \hom(L'_\aaa, L'_\bbb)|_{V'_i} 
\simeq u_*(\hom(L'_\aaa, L'_\bbb)|_{V'_j}) 
\simeq u_*(\hom(L_{\aaa,i}, L_{\bbb,i})_{T^j_{\aaa\bbb}}).
\end{equation}
By lemma~\ref{tubes_fermes} (iii), the inclusions $T^j_{\aaa\bbb}
\subset V'_j \subset V'_i$ are locally homeomorphic to inclusions of
convex subsets of $\R^d$, and~\eqref{eq:homlaaalbbb_bis} follows.

\medskip 

Now we define $g^{\aaa\bbb}_i$. Since $\fA^{\aaa\bbb}_i =
\fA^{\aaa\bbb}_j$, we just have to check that $g^{\aaa\bbb}_j$
factors through the restriction morphism $\fB^{\aaa\bbb}_i \to
\fB^{\aaa\bbb}_j$. As in (i), formula~\eqref{eq:homlaaalbbb_bis}
implies the existence of a morphism
$$
f^i_{\aaa\bbb}: \Omega^\cdot_{V'_i} \otimes
\hom(L_{\aaa,i},L_{\bbb,i}) \to \Omega^{\cdot + 2d_{\aaa\bbb}}_{V'_i}
\otimes \hom(L'_\aaa, L'_\bbb).
$$
We also note, by~\eqref{eq:homlaaalbbb_ter}, that $\supp
\hom(L'_\aaa, L'_\bbb)|_{V'_i} \subset \ovl{V'_j}$.  Since $\ovl{V'_j}
\subset V_j$, we obtain
$$
\fB^{\aaa\bbb}_i = \sect(V'_i \cap V_j; \Omega_Y \otimes
\hom(L'_\aaa, L'_\bbb)).
$$
Hence $\sect(V'_i \cap V_j; f^i_{\aaa\bbb})$ yields a morphism from
$$
B'_{\aaa\bbb j} = \sect(V'_i \cap V_j; \Omega_Y \otimes
\hom(L_{\aaa,i},L_{\bbb,i})) \,[-2d_{\aaa\bbb}]
$$
to $\fB^{\aaa\bbb}_i$.  Composed with the restriction morphism from
$\fA^{\aaa\bbb}_j$ to $B'_{\aaa\bbb j}$, it gives the required
morphism $\fA^{\aaa\bbb}_j \to \fB^{\aaa\bbb}_i$.

Since $g^{\aaa\bbb}_i$ is defined by factorising $g^{\aaa\bbb}_j$, it
is clear that the $g^{\aaa\bbb}_\cdot$ commute with the restriction
maps and define a morphism of sheaves, $g^{\aaa\bbb}: \fA^{\aaa\bbb}
\to \fB^{\aaa\bbb}$.

\medskip

Now, let us check that $g^{\aaa\bbb}$ is a quasi-isomorphism.
We have to prove (with the notations~\ref{not:IetIprime}):
\begin{itemize}
\item [(a)] for $i\in I_{\aaa\bbb}$, $H^\cdot (g^{\aaa\bbb}_i)$
  is an isomorphism,
\item[(b)] for $i \in I'_{\aaa\bbb}$ and $j\in I_{\aaa\bbb}$ such that
  $V_i \setminus (\bigcup_{v\in \Delta'_\aaa \cap \Delta'_\bbb} D_v) =
  V_j$, we have $H^\cdot (\fB^{\aaa\bbb}_i) \simeq H^\cdot
  (\fB^{\aaa\bbb}_j)$,
\item[(c)] for $i\not\in I_{\aaa\bbb} \sqcup I'_{\aaa\bbb})$, $H^\cdot
  (\fB^{\aaa\bbb}_i) =0$.
\end{itemize}
By the definition of $g^{\aaa\bbb}_i$, (a) follows directly from
quasi-isomorphism~\eqref{gysin_de_Rham}.  Let us verify (b). Since
$L_\aaa|_{D_v} =0$ for $v\in \Delta'_\aaa$, we have $L_\aaa |_{V_i} =
(L_\aaa)_{V_j}$. Hence, using~\eqref{eq:section_de_A}, we have the
required isomorphism:
$$
H^\cdot (\fB^{\aaa\bbb}_i) \simeq \Ext^\cdot( L_\aaa|_{V_i},
L_\bbb|_{V_i}) \simeq \Ext^\cdot( L_\aaa|_{V_j}, L_\bbb|_{V_j})
\simeq H^\cdot (\fB^{\aaa\bbb}_j).
$$
Let us prove (c). We first note that
assumption~\ref{assumptionstratification} (ii) implies, for $F\in
\D^b(Y)$, constructible with respect to the stratification
$Y=\bigsqcup_{i\in I} Y_i$, $\forall i \in I$, $H^\cdot(V_i;F) \simeq
H^\cdot(Y_i;F)$.  By~\eqref{eq:section_de_A}, it follows that:
\begin{multline*}
  H^\cdot (\fB^{\aaa\bbb}_i) \simeq \Ext^\cdot( L_\aaa|_{V_i},
  L_\bbb|_{V_i}) \\
  \simeq H^\cdot(V_i;\Rhom(L_\aaa,L_\bbb)) \simeq
  H^\cdot(Y_i;\Rhom(L_\aaa,L_\bbb)|_{Y_i}).
\end{multline*}
Now, we have either $i\not\in \pxf( Z_{\Delta_\aaa \cup \Delta_\bbb})$
or there exists $v_0\in \Delta'_{\aaa\bbb} = (\Delta'_\aaa \cup
\Delta'_{\bbb} ) \setminus (\Delta'_\aaa \cap \Delta'_{\bbb} )$ such
that $Y_i \subset D_{v_0}$.  In the first case $Y_i$ doesn't meet
$\supp(\Rhom(L_\aaa,L_\bbb))$ and the vanishing of $H^\cdot(
\fB^{\aaa\bbb}_i )$ is clear. Let us assume we are in the second case.
For $\aaa \in \epm$, we set $V_\aaa = Y\setminus \bigcup_{v\in
  \Delta'_\aaa} D_v$. By lemma~\ref{lem:monodromy} we have $L_\aaa
\simeq (L_\aaa)_{V_\aaa} \simeq \sect_{V_\aaa}(L_\aaa)$.  Let us set
$V = V_\aaa \cap V_\bbb$ and let $j:V\to Y$ be the inclusion. We
obtain:
\begin{multline*}
  \Rhom(L_\aaa,L_\bbb) \simeq
  \Rhom((L_\aaa)_{V_\aaa},\rsect_{V_\bbb}(L_\bbb))  \\
  \simeq \rsect_V \Rhom(L_\aaa,L_\bbb) \simeq Rj_*
  \Rhom(L_\aaa|_V,L_\bbb|_V).
\end{multline*}
Now $Z_\aaa \cap V$ and $Z_\bbb \cap V$ are closed in $V$ and their
intersection $Z_{\aaa\bbb} \cap V$ is empty or a smooth submanifold of
$Z_\bbb \cap V$ of codimension $2d_{\aaa\bbb}$. Hence
$\Rhom(L_\aaa|_V,L_\bbb|_V)$ is isomorphic to
$K_{\aaa\bbb}[2d_{\aaa\bbb}]$, where $K_{\aaa\bbb}$ is a local system
on $Z_{\aaa\bbb} \cap V$.  By
assumption~\ref{assumptionstratification} (v), the monodromy of
$K_{\aaa\bbb}$ around $D_{v_0}$ is $-Id$.  By
lemma~\ref{lem:monodromy} we obtain $Rj_*K_{\aaa\bbb}|_{D_{v_0}} = 0$
and, a fortiori $\Rhom(L_\aaa,L_\bbb)|_{Y_i} =0$ and $H^\cdot(
\fB^{\aaa\bbb}_i) =0$.
\end{proof}

\section{Symmetric varieties}
\label{symmetricvarieties}
We recall some results of~\cite{BDP90} on regular compactifications of
homogeneous symmetric varieties, in particular the structure of the
decomposition by the $K$-orbit types, for a suitable maximal compact
subgroup $K$.  Then we show that the hypothesis of
proposition~\ref{prop:eq_cat1} are satisfied.

\medskip Let $G$ be a semi-simple algebraic group of adjoint type over
$\C$, $\sigma$ an automorphism of order $2$ of $G$ and $H=G^\sigma$.
Let $T$ be a $\sigma$-stable maximal torus of $G$ containing a maximal
$\sigma$-split torus $\Split$ (i.e. $\forall t\in \Split$, $\sigma(t)
=t^{-1}$). The corresponding root system $\Phi= \Phi(G,T)$ decomposes
as $\Phi=\Phi_0\sqcup \Phi_1$, where $\Phi_0$ denotes the set of roots
fixed by the action of $\sigma$.  One may choose a basis of simple
roots $\Sigma$ such that $\sigma$ exchanges the corresponding positive
roots of $\Phi_1$ with the negative roots of $\Phi_1$.  The non-fixed
roots $\Phi_1$ induce a root system on $\Split$ with basis
$\{\gamma_1,\ldots,\gamma_l\}$ given by the restriction of $\Sigma
\cap \Phi_1$. The corresponding Weyl group is denoted by $W^1$. We set
$D=H\cap \Split$, the subgroup of $\Split$ of elements of order $2$, and
$\Splitbis= \Split/D \simeq T/(T\cap H)$.  The natural map $\Split \to
\Splitbis$ gives an identification between the Lie algebras $Lie(\Split)$
and $Lie(\Splitbis)$.  Let $H^0$ be the identity component of $H$. By
proposition~1 of~\cite{KR71} (see also proposition~7 of~\cite{V74}),
we have $H = D \cdot H^0$. In particular the group of components of
$H$ is a quotient of $D$, hence of the type $H/H^0 \simeq (\Z/2\Z)^a$,
for some $a\in \N$.

\subsection*{Regular compactifications}
Let $X$ be the canonical compactification of $G/H$ described
in~\cite{DP83} and~\cite{DP85}. It can be defined as follows: let
$Gr_n$ be the Grassmann variety of $n$-dimensional subspaces of
$Lie(G)$, for $n=\dim (H)$, with the $G$-action induced by the adjoint
action. Let $x\in Gr_n$ be the point associated to $Lie(H)$.  One can
show that $G\cdot x \simeq G/H$ and that $X$ is isomorphic to the
closure of $G\cdot x$ in $Gr_n$.  It is proved in~\cite{DP83} that $X$
is smooth, $X\setminus (G/H)$ is the union of $l$ smooth, normal
crossings divisors, say $D_i$, $i=1,\ldots,l$, which are closures of
$G$-orbits, and any $G$-orbit closure is the intersection of the $D_i$
containing it. More precisely, the decomposition into $G$-orbits is
identified with the decomposition of the toric variety $\C^l$ into
orbits for the action of $(\C^*)^l$, as follows. The inclusion $\Split
\subset G$ gives an embedding of $\Splitbis$ in $G/H$, and hence in $X$.
The closure of $\Splitbis$ in $X$ is a $\Splitbis$-toric variety, whose fan
can be identified with the subdivision of $Lie(\Splitbis)$ into Weyl
chambers under the action of $W^1$.  We consider the affine space
$\C^l$ associated to the negative Weyl chamber.  Then the $G$-orbits
in $X$ correspond bijectively (by taking the intersection with $\C^l$)
to the $\Splitbis$-orbits in $\C^l$. In particular, there are
$2^l$-orbits and one single closed orbit.  Let $\O \subset X$ be a
$G$-orbit. Then its closure is fibred over a variety of parabolic
subgroups of $G$, with fibre a symmetric variety. More precisely,
there exist a parabolic subgroup $P\subset G$, and a $G$-equivariant
fibration $\ovl{\O} \to G/P$, whose fibre, say $X_{\O}$, has the
following description. There exists a $\sigma$-stable Levi subgroup
$L$ of $P$, such that, denoting by $L'$ the quotient $L/Z(L)$ ($Z(L)$
is the centre of $L$), $X_{\O}$ is the canonical compactification of
$L'/L'^\sigma$.


\medskip

In~\cite{DP85} there is a description of the embeddings of $G/H$ over
$X$.  If $Y$ is such an embedding, with a map $\pi:Y \to X$, the
closure of $\Splitbis$ in $Y$ gives a toric variety, say $Z'$, and $Z=
Z'\cap \pi^{-1}(\C^l)$ is a toric variety over $\C^l$. It is shown
in~\cite{DP85} that this gives a bijection between the embeddings of
$G/H$ over $X$ and the toric varieties over $\C^l$.  For $Z\to \C^l$ a
morphism of toric varieties, let $\pxzx:X_Z\to X$ be the corresponding
embedding of $G/H$. Then $Z=\pxzx^{-1}(\C^l)$ and the $G$-orbits in
$X_Z$ correspond bijectively (by taking the intersection with $Z$) to
the $\Splitbis$-orbits in $Z$.  Moreover $X_Z$ is smooth if, and only if,
$Z$ is smooth and $X_Z$ is complete if, and only if, $Z\to \C^l$ is
proper.  From now on we assume that $X_Z$ is smooth and complete.
These are the symmetric varieties we consider here.
\begin{notations}
  We denote by $V$ the set of irreducible $G$-stable divisors, $D_v$,
  $v\in V$.  Any $G$-orbit closure is the intersection of the $D_v$
  containing it.  For $\Delta\subset V$ such that $\bigcap_{v\in
    \Delta} D_v \not=\emptyset$, we let $\O_\Delta$ be the $G$-orbit
  such that $\ovl{\O_\Delta} = \bigcap_{v\in \Delta} D_v$.  We denote
  by $\S$ the set of $G$-orbits; hence $\S$ is identified with a set
  of subsets of $V$.
\end{notations}

\subsubsection*{Fundamental domain}
In~\cite{BDP90}, we also have a description of the orbits of a
suitable maximal compact subgroup of $G$. Let $K$ be a compact form of
$G$ such that the Cartan involution of $G$ corresponding to $K$
commutes with $\sigma$ and such that $K\cap T$ is a maximal compact
subgroup of $T$. We set $\Scomp = K\cap \Split$.  Let $t_i=
t^{-2\gamma_i}$, $i=1,\ldots, l$, be the characters of $\Splitbis$
associated to the simple roots $\gamma_i$, extended to a coordinates
system on $\C^l$. We set $C=]0,1]^l\subset \Splitbis \subset G/H$. We
consider the closure of $C$ in $X_Z$, $C_{X_Z}= \ovl{C}$.  In
particular, $C_X=[0,1]^l$.  Note that $(\pxzx|_Z)^{-1}([0,1]^l)$ is
closed and contains $C$, so that $C_{X_Z}$ is in fact contained in
$Z$.  Hence $C_{X_Z}$ is mapped to $C_X$ by $\pxzx$.  This is a
fundamental domain for the action of $K$ on $X_Z$ (see~\cite{BDP90}
theorem~27).

\subsubsection*{Stabilisers}
The stabiliser in $K$ of a point of $C_{X_Z}$ is described
in~\cite{BDP90}, p. 27, as follows.  For a point $q\in C_X=[0,1]^l$,
we call its $J$-support the subset of $\{1,\ldots, l\}$ defined by
$J(q)= \{i;\; q_i \not= 1\}$.  Let $\ptutd:\Split\to \Splitbis$ denote the
quotient map. For $q\in ]0,1]^l \subset \Splitbis$, let $\tilde{q}\in
\ptutd^{-1}(q)$ be in the connected component of
$\ptutd^{-1}(]0,1]^l)$ containing $1$. Then the centraliser of
$\tilde{q}$ in $G$ only depends on $J=J(q)$; we set:
$$
\Ksharp_J = K^\sigma \cap Z_G(\tilde{q}).
$$
For $J\subset \{1,\ldots,l\}$, we denote by $C_{X_Z,J}$, or simply
$C_J$ if there is no ambiguity, the subset of $C_{X_Z}$ formed by the
$p$ such that $J(\pxzx(p))= J$. By definition this decomposition of
$C_{X_Z}$ arises from the decomposition of $C_X = [0,1]^l$ according
to the set of coordinates equal to $1$: $C_{X_Z,J} =
(\pxzx|_{C_{X_Z}})^{-1} (C_{X,J})$.

We also consider the partition of $C_{X_Z}$ given by the $G$-orbits:
$C_{X_Z}= \bigsqcup_{\Delta\in \S} \O_\Delta \cap C_{X_Z}$ .  All
points $p\in \O_\Delta \cap C_{X_Z}$ have the same stabiliser in
$\Scomp$.  Indeed $C_{X_Z} \subset Z$, and $\O_\Delta \cap Z$ is a
single $\Splitbis$-orbit of the toric variety $Z$, so that all points
of $\O_\Delta \cap Z$ have the same stabiliser in $\Splitbis$. We set
$\Scomp_\Delta = \Scomp_p$ for any $p\in \O_\Delta \cap C_{X_Z}$.

Now we mix the two partitions above and set, for $\Delta\in\S$ and
$J\subset \{1,\ldots,l\}$:
$$
F_{\Delta,J}= \O_\Delta \cap C_J.
$$
The $F_{\Delta,J}$ are called the ``faces'' of $C_{X_Z}$ (if $X_Z$
is projective, the moment map for the $K$-action identifies $C_{X_Z}$
with a polytope, and the ``faces'' are the usual faces of this
polytope -- see~\cite{BDP90}, p. 25).  We have the following
description of the stabilisers:
\begin{theorem}[theorem~32 and corollary~35 of ~\cite{BDP90}]
  \label{thm:BDP}
  For $p\in C_{X_Z}$ let $F_{\Delta,J}$ be the face containing $p$,
  and let $K_p$ be the stabiliser of $p$ in $K$.  Since $\O_\Delta$ is
  $K$-stable, $K_p$ acts on $N_p = T_pX_Z / T_p \O_\Delta$.  We have:
\begin{gather}
  \label{eq:noyauKpversNp}
\Ksharp_J = \ker( K_p \to GL(N_p) ),  \\
  \label{eq:splitting}
  K_p = \Scomp_\Delta \cdot \Ksharp_J  , \qquad
\Scomp_\Delta  \cap \Ksharp_J= D .
\end{gather}
\end{theorem}
In particular $K_p$ only depends on the face to which $p$ belongs.
For a face $F$, we set $K_F = K_p$, for any $p\in F$. Let us make some
remarks on the decomposition of $C_{X_Z}$ into faces.

For $X_Z=X$, we have $Z=\C^l$, $C_X = [0,1]^l$. We may
identify $V$ with $\{1,\ldots,l\}$ so that, for $\Delta\subset V$,
$\O_\Delta\cap \C^l = \{(t_1,\ldots,t_l);$ $t_i =0$ iff $i\in
\Delta\}$. Hence:
$$
\text{for $X_Z =X$}:\quad F_{\Delta,J} = \{(t_1,\ldots,t_l)\in
[0,1]^l; \:\text{$t_i =0$ iff $i\in \Delta$ and $t_i=1$ iff $i\notin
  J$}\}.
$$
The image of a $G$-orbit of $X_Z$, say $\O_\Delta$, $\Delta\subset
V$, by $\pxzx$ is a $G$-orbit of $X$, say $\O_{\Sigma}$, $\Sigma
\subset \{ 1,\ldots,l\}$.  The restriction
$\pxzx|_{\O_\Delta}:\O_\Delta \to \O_{\Sigma}$ is a fibration, and
gives a bijection between the faces $F_{\Delta,J}$ of $C_{X_Z}\cap
\O_{\Delta}$ and the faces $F_{\Sigma,J}$ of $C_X\cap \O_{\Sigma}$.
In particular $C_{X_Z}\cap \O_{\Delta}$ has a unique closed face,
$F_{\Delta,\Sigma}$, and a unique open face, $F_{\Delta,
  \{1,\ldots,l\}}$.  We also deduce that $F_{\Delta,J} \subset
\ovl{F_{\Delta,J'}}$ if and only if $J\subset J'$.

We have similarly $F_{\Delta,J} \subset \ovl{F_{\Delta',J}}$ if and
only if $\Delta' \subset \Delta$. Indeed if $F_{\Delta,J} \subset
\ovl{F_{\Delta',J}}$, then certainly $\O_\Delta \cap
\ovl{\O_{\Delta'}} \not=\emptyset$ and hence $\Delta' \subset \Delta$.
Conversely, assume $\Delta' \subset \Delta$ and let $p\in
F_{\Delta,J}$.  Since $p\in \ovl{\O_{\Delta'}}$, we may write
$p=\lim_{n} p_n$, with $p_n\in \O_{\Delta'}$. Now each $p_n$ is itself
a limit of points of $\O_\emptyset$, which is identified with
$(\C^*)^l$ by $\pxzx$.  Let us write $p_n = \lim_i q^n_i$, with $q^n_i
=(t^n_{1,i},\ldots,t^n_{l,i})$.  Since $p\in C_J$, for any
$\varepsilon>0$, there exists a neighbourhood $U$ of $p$, such that
$q^n_i \in U$ implies $|t^n_{j,i} -1| < \varepsilon$, $\forall j
\not\in J$. Hence, up to restricting to a subsequence, we may assume
that, for each $n$ great enough, $\forall j\not\in J$, $\lim_i
t^n_{j,i} = s^n_j$ for some $s^n_j$, with $|s^n_j -1| \leq
\varepsilon$. We set $s^n_j =1$ for $j\in J$.  Then, for the element
$g_n=(\udl s^n)\in (\C^*)^l$, the point $p'_n = g^{-1}_n\cdot p_n =
\lim_i g^{-1}_n \cdot q^n_i$ is in $\O_{\Delta'} \cap\C_J$. The $p'_n$
converge to $p$, and we have $p\in \ovl{F_{\Delta',J}}$, as required.

Combining both characterisations for the inclusions of closures of
faces, we have: $F_{\Delta,J} \subset \ovl{F_{\Delta',J'}}$ is
equivalent to $\Delta' \subset \Delta$ and $J\subset J'$.  Now we
summarise the notations and properties introduced so far and add some
others.
\begin{properties}
  For $\Delta\in\S$, $J\subset \{1,\ldots,l\}$, we have $F_{\Delta,J}=
  \O_\Delta \cap C_{X_Z,J}$. We let $\F$ denote the set of faces, $\F=
  \{ F_{\Delta,J};\: F_{\Delta,J} \not=\emptyset \}$. We let $\pcxf:
  C_{X_Z} \to \F$ be the natural map induced by this partition and
  endow $\F$ with the quotient topology.  For $\Delta\in\S$, $\O_\Delta
  \cap C_{X_Z}$ has a unique closed face, which we denote by
  $F_{\Delta,J_\Delta}$, $J_\Delta \subset \{1,\ldots,l\}$. As in
  section~\ref{eq_der_cat}, for a face $F$ of $C_{X_Z}$, we let $U_F$
  be the smallest open subset of $\F$ containing $F$. We set $\uu_F =
  \pcxf^{-1}(U_F)$; this is an open subset of $C_{X_Z}$ which contains
  $F$ as its unique closed face.
\begin{gather}
  \label{descriptionstabilisers}
  K_{F_{\Delta,J}} = \Scomp_\Delta  \cdot \Ksharp_J  , \\
  \label{inclusiostabil}
  \Delta' \subset \Delta \: \Longrightarrow\: \Scomp_{\Delta'} \subset
  \Scomp_\Delta, \qquad
  J\subset J' \:\Longrightarrow\: \Ksharp_{J'} \subset \Ksharp_J, \\
  \label{inclusionadhfaces}
  F_{\Delta,J} \subset \ovl{F_{\Delta',J'}}
  \qquad\Longleftrightarrow\qquad
  \Delta' \subset \Delta \quad \text{and} \quad J\subset J', \\
    \textstyle
  U_F= \{F' \in \F;\: F\subset \ovl{F'}\}, \qquad \uu_F=
  \bigsqcup_{F\subset \ovl{F'}}F', \\
  \textstyle
  \label{descriptionUF}
  \uu_{F_{\Delta,J}} \:  = \: \bigsqcup_{\Delta'\subset \Delta, J'\supset J}
  F_{\Delta',J'} \: = \: (\bigsqcup_{\Delta'\subset \Delta}\O_{\Delta'})
  \cap (\bigsqcup_{J'\supset J} C_{J'}), \\
  \label{intersectionUF}
  U_{F_{\Delta,J}} \cap U_{F_{\Delta',J'}} = U_{F_{\Delta\cap
      \Delta',J\cup J'}}.
\end{gather}
The first equality in~\eqref{descriptionUF} follows
from~\eqref{inclusionadhfaces}. The second follows directly by
applying the definition of the faces, and it
implies~\eqref{intersectionUF}.
\end{properties}

\subsection{Stratification by the faces}
\label{stratificationfaces}
\begin{lemma}
\label{lem:stratfaces}
  We keep the notations introduced above. Let $X'_Z=X_Z /K$ be the
  topological quotient, $p_Z:X_Z \to X'_Z$ the quotient map and
  $q_Z=p_Z|_{C_{X_Z}}$. We consider on $C_{X_Z}$ the topology induced
  by its inclusion in $X_Z$. We have:
  \begin{itemize}
  \item [(i)] the map $q_Z$ is a homeomorphism,
  \item[(ii)] the partition of $C_{X_Z}$ by the faces, $C_{X_Z} =
    \bigsqcup_{F\in \F} F$, satisfies: if $F\cap \ovl{F'}
    \not=\emptyset$ then $F\subset \ovl{F'}$,
  \item[(iii)] the induced partition of $X_Z$, $X_Z = \bigsqcup_{F\in
      \F} K\cdot F$, is a $\mu$-stratification; it satisfies the same
    inclusions relations for the closures of strata as the partition
    of $C_{X_Z}$ by the faces.
  \end{itemize}
\end{lemma}
\begin{proof}
  (i) Since $C_{X_Z}$ is a fundamental domain for the $K$-action,
  $q_Z$ is bijective.  It is continuous by definition. For an open
  subset $U\subset C_{X_Z}$, $C_{X_Z}\setminus U$ is compact because
  $C_{X_Z}$ is. Hence $q_Z(C_{X_Z}\setminus U)$ is compact too, and
  $q_Z(U)$ is open. This proves that $q^{-1}_Z$ is continuous too.
  
  (ii) Let $F_{\Delta,J}$, $F_{\Delta',J'}$ be two faces such that
  $F_{\Delta,J} \cap \ovl{F_{\Delta',J'}} \not=\emptyset$. By
  definition of the faces, this implies $\O_\Delta \cap
  \ovl{\O_{\Delta'}} \not=\emptyset$ and $C_J \cap \ovl{C_{J'}}
  \not=\emptyset$. Hence $\Delta'\subset \Delta $ and $J\subset J'$;
  we conclude by~\eqref{inclusionadhfaces}.
  
  (iii) By (i), we have for any subset $C\subset C_{X_Z}$:
  $$
  \ovl{K\cdot C} = p_Z^{-1} ( \ovl{ p_Z(K\cdot C)} ) = p_Z^{-1} (
  \ovl{q_Z(C)}) = p_Z^{-1} (q_Z( \ovl{C}))= K\cdot \ovl{C}.
  $$
  This implies that our partition of $X_Z$ is a stratification and
  the last assertion. Let us verify the ``$\mu$-condition'' for two
  strata $K\cdot F$, $K\cdot F'$, with $F\subset \ovl{F'}$
  (see~\ref{subsection:constrsheaves}). This condition is local around
  a point of $K\cdot F$ and we may restrict ourselves to $K \cdot
  \uu_F$.  In $K\cdot \uu_F$, our stratification coincides,
  by~\eqref{descriptionstabilisers}, with the partition by the
  $K$-orbit types. By the existence of slices for action of compact
  groups this partition is locally trivial: a point $x\in K\cdot F$
  has a neighbourhood of the type $\R^d \times E_\lambda$, where
  $E_\lambda$ is a representation of $K_x$ with $0$ as unique fixed
  point, and the partition is induced by the partition of $E_\lambda$
  in $K_x$-orbit types. In $E_\lambda$, the $\mu$-condition for the
  strata $\{0\}$ and $(K\cdot F') \cap E_\lambda$ is trivially
  satisfied; hence it holds for $K\cdot F$ and $K\cdot F'$.
\end{proof}
\begin{notations}
  \label{not:stratifXZ}
  We stratify $X_Z$ as in the above lemma. We denote by $\pxf:X_Z \to
  \F$ the continuous map defined by $\pxf(K\cdot F) = F$, for $F\in
  \F$.  We set $\vv_F = K\cdot \uu_F =\pxf^{-1}(U_F)$.  From now on we
  denote by $E$ a universal bundle for $K$ which is an increasing
  union of manifolds, $E=\bigcup_k E_k$. Since $\poinc(E) = 1$, for
  any subgroup $H\subset K$ we have $\poinc(E/H) \simeq H/H^0$. We let
  $\pxef: E\times_K X_Z \to \F$ be the map induced by $\pxf$.
\end{notations}
We are in the setting of assumption~\ref{assumptionstratification},
with $Y=X_Z$ in a $K$-equivariant way, or $Y=E\times_K X_Z$ (see
remark~\ref{rem:espclass}), $I=\F$, and a set of normal crossings
divisors $D_v = \ovl{\O_{\{v\}}}$, $v\in V$ (we will introduce local
systems $L_\aaa$ in the next paragraph).  Let us verify conditions
(i)--(iii). The first condition follows directly
from~\eqref{intersectionUF}. By~\eqref{descriptionUF}, we have, for
$\Delta \in \S$ and $v\in \Delta$,
\begin{equation}
\label{eq:verifassumiii}
\begin{split}
  \textstyle \uu_{F_{\Delta,J}} \setminus D_v & \textstyle =
  ((\bigsqcup_{\Delta'\subset \Delta}\O_{\Delta'}) \setminus D_v) \cap
  (\bigsqcup_{J'\supset J} C_{J'})  \\
  & \textstyle = (\bigsqcup_{v\not \in\Delta'\subset \Delta }
  \O_{\Delta'}) \cap (\bigsqcup_{J'\supset J} C_{J'}) =
  \uu_{F_{\Delta\setminus \{v\},J}},
\end{split}
\end{equation}
and this gives condition (iii), the case $v\not\in \Delta$ being
trivial. Condition (ii) of~\ref{assumptionstratification} follows from
the next lemma.
\begin{lemma}
  \label{lem:contraction_faces}
  Let $F$ be a face.

  (i) There exists a homotopy $h:[0,1]\times \uu_F \to \uu_F$, such that
  $h_1=id$ and $h_0$ is the projection of $\uu_F$ to a point of $F$,
  with the property that the closures of the faces of $\uu_F$ are stable
  under $h_t$, $\forall t\in [0,1]$.
  
  (ii) It induces a $K$-equivariant homotopy $\bar{h}:[0,1]\times
  \vv_F \to \vv_F$ contracting $\vv_F$ to a $K$-orbit $K/K_F \subset
  K\cdot F$ and preserving the closures of strata.
\end{lemma}
\begin{proof}
  (i) By definition, the faces of $C_{X_Z}$ are contained in
  $\Splitbis$-orbits of $Z$.  Let $Y$ be the $\Splitbis$-orbit
  containing $F$ and let $U$ be the union of the $\Splitbis$-orbits
  whose closure contains $Y$.  Then $\uu_F$ is included in $U$. There
  exists a one-dimensional subtorus $i:\C^* \hookrightarrow \Splitbis$
  of $\Splitbis$ contracting $U$ to $Y$, i.e.  $\forall y\in U$
  $\lim_{t\to 0} i(t)\cdot y \in Y$. Note that for $t\in ]0,1]$, $x\in
  U\cap C_{X_Z}$, $i(t)\cdot x$ still is in $U\cap C_{X_Z}$ and
  belongs to the same face as $x$. Hence the homotopy $h_1:[0,1]\times
  \uu_F \to \uu_F$, $(t,x) \mapsto i(t)\cdot x$ if $t\not=0$, $(0,x)
  \mapsto \lim_{t\to 0}i(t)\cdot x$, contracts $\uu_F$ to $\uu_F\cap
  Y$ and preserves the closure of the faces.
  
  Now $Y$ is fibred over $\pxzx(Y)$ with fibre $(\C^*)^k$, a quotient
  of $\Splitbis$.  It follows that $C_{X_Z}\cap Y$ is fibred over
  $C_X\cap \pxzx(Y)$, with fibre $W$, where $W$ is a subset of $\R_{>
    0}^k$ stable by multiplication by $]0,1]^k$. In particular $W$ is
  contractible.  Moreover this fibration is a bijection on the faces.
  Hence there exist sections $s:C_X\cap \pxzx(Y) \to C_{X_Z}\cap Y$ of
  $\pxzx|_{C_{X_Z}\cap Y}$, and for any such $s$, we may find a
  homotopy contracting $C_{X_Z}\cap Y$ to $\im(s)$, which preserves
  the fibres and thus the faces.
  
  Finally, $\im(s) \simeq C_X\cap \pxzx(Y)$ is a union of standard
  faces in $[0,1]^l$, isomorphic to $]0,1]^m\times \{0\}^{l-m}$, up to
  a permutation of coordinates.  Then $\im(s)\cap F \simeq
  ]0,1[^n\times \{1\}^{m-n} \times \{0\}^{l-m}$ and $\im(s)\cap \uu_F
  \simeq ]0,1[^n\times ]0,1]^{m-n} \times \{0\}^{l-m}$.  Hence there
  exists a third homotopy contracting $\im(s)\cap \uu_F$ to a point of
  $F$ also preserving the closure of the faces.
  
\smallskip

  (ii) By~\eqref{inclusiostabil} and~\eqref{inclusionadhfaces}, for
  any faces $F_1$, $F_2$ with $F_1\subset \ovl{F_2}$, we have $K_{F_2}
  \subset K_{F_1}$.  Since $\forall F'\in U_F$, $h([0,1] \times F')
  \subset \ovl{F'}$, we obtain $\forall (t,x) \in [0,1] \times \uu_F$,
  $K_x \subset K_{h(t,x)}$. Hence it makes sense to define $\bar{h} :
  [0,1] \times \vv_F \to \vv_F$, $(t,k\cdot x) \mapsto k\cdot h(t,x)$,
  for $t\in [0,1]$, $k\in K$, $x\in \uu_F$.  It is clear that
  $\bar{h}_0 = id$ and $\bar{h}_1$ is a $K$-equivariant projection to
  $K/K_F$. Let us verify that $\bar{h}$ is continuous.
  
  Let $(t_n,y_n)$ be any sequence in $[0,1] \times \vv_F$ converging
  to $(t,y)$. Let us see that a subsequence of $\bar{h}(t_n,y_n)$
  converges to $\bar{h}(t,y)$.  For each $n$ there exists a unique
  $x_n\in \uu_F$ such that $y_n \in K\cdot x_n$. With the notations of
  lemma~\ref{lem:stratfaces}, we have $x_n = q_Z^{-1}(p_Z(y_n))$.
  Since these maps are continuous, the sequence $x_n$ converges to $x
  =q_Z^{-1}(p_Z(y))$.  Let us write $y_n = k_n \cdot x_n$, $k_n \in
  K$. Since $K$ is compact, we may assume, up to restriction to a
  subsequence, that $k_n$ converges to $k\in K$. Then
  $\bar{h}(t_n,y_n) = k_n \cdot h(t_n,x_n)$ converges to $\bar{h}(t,y)
  = k\cdot h(t,x)$, as desired.
\end{proof}
Let us give some immediate consequences of this lemma.  We consider
the map $p_F:\vv_F \to K/K_F$, $k\cdot x \mapsto \ovl{k}$, where $k\in
K$, $x\in \uu_F$. This makes sense because $\forall x\in \uu_F$, $K_x
\subset K_F$ (recall that $\uu_F = \bigsqcup_{F\subset \ovl{F'}} F'$).
With the notations of lemma~\ref{lem:contraction_faces}, let $x_0 \in
F$, be the point of $F$ such that $h(1,U'_F) = \{x_0\}$.  By
definition, we have in fact $p_F = \bar{h}(1,\cdot)$ modulo the
identification $K\cdot x_0 = K/K_F$. In particular $p_F$ is
continuous. Since it is $K$-equivariant, it is a fibration over
$K/K_F$, with fibre $K_F\cdot \uu_F$. The homotopy $\bar{h}$ also is
$K$-equivariant and hence contracts each fibre $p_F^{-1}(y)$, $y\in
K/K_F$ to the point $y$.  We let
\begin{equation}
  \label{eq:def_qF}
  q_F: E\times_K \vv_F \to E/K_F
\end{equation}
be the map induced by $p_F$.  This also is a fibration with
contractible fibres.

For a $G$-orbit $\O_\Delta$, with closed face $F_{\Delta,J_\Delta}$,
and a face $F' = F_{\Delta,J'}$ such that $F'\subset \O_\Delta$ (i.e.
$J_\Delta \subset J'$), $\O_\Delta \cap \vv_{F'}$ is closed in
$\vv_{F'}$: indeed if a face $F_{\Delta_1,J_1}$ is included in
$\ovl{\O_\Delta}$ it satisfies $\Delta \subset \Delta_1$, and if it is
in $\vv_{F'}$ it satisfies $\Delta_1 \subset \Delta$; hence any face
of $\ovl{\O_\Delta} \cap \vv_{F'}$ is in $\O_\Delta$.  Since the
homotopy $h$ of lemma~\ref{lem:contraction_faces} preserves the
closures of faces, it follows that $h$ contracts $\O_\Delta\cap
\vv_{F'}$ to $K/K_F$. Hence the maps
\begin{equation}
  \label{eq:homot_equi_orbit}
  E\times_K (\O_\Delta \cap \vv_{F'}) \to E\times_K \vv_{F'} \to
E/K_{F'}
\end{equation}
are homotopy equivalences. In particular the fundamental groups are
the same:
\begin{equation}
  \label{eq:groupfond}
 \poinc(  E\times_K (\O_\Delta \cap \vv_{F'}) ) = 
\poinc ( E/K_{F'} ) = K_{F'}/K_{F'}^0.
\end{equation}

\subsection{Equivariant local systems on $X_Z$}
\label{eq_loc_syst}
Remember that $G$-equivariant local systems on a homogeneous variety
$G/G'$ are in bijective correspondence with representations of the
components group $G'/G'^0$ of $G'$.  For a $G$-orbit $\O$, we denote
by $\tau_\O$ the group of components of a stabiliser, $\tau_\O = G_p/
G^0_p$, for $p\in \O$. We introduce similar notations for the groups
of components of the groups defined up to now:
\begin{gather*}
  \Scompz_\Delta = (\Scomp_\Delta)^0, \qquad\quad
  D_\Delta = D \cap \Scompz_\Delta, \\
  \tau_J = \Ksharp_J / \Ksharp_J^0, \qquad \tau_\Delta = \Scomp_\Delta
  / \Scompz_\Delta = D/D_\Delta, \qquad \tau_F = K_F/K_F^0.
\end{gather*}
For a $G$-orbit $\O_\Delta$, remember that $F_{\Delta,J_\Delta}$
denotes the unique closed face of $\O_\Delta\cap C_{X_Z}$. For $p\in F
= F_{\Delta,J_\Delta}$, $K_p$ is a maximal compact subgroup of $G_p$
(see~\cite{BDP90}, p. 31), hence $\tau_\O \simeq \tau_F$.  Let us
remark that $\tau_J \simeq (\Z/2\Z)^a$, for some $a\in\N$. Indeed,
by~\eqref{eq:splitting}, we have $K_F = \Scomp_\Delta \cdot
\Ksharp_J$, but $K_F$ is connected and $\Scomp_\Delta$ and $\Ksharp_J$
are compact, so that we have as well $K_F = \Scompz_\Delta \cdot
\Ksharp_J^0$.  Hence, for $k\in \Ksharp_J$, there exist $s\in
\Scompz_\Delta$, $k' \in \Ksharp_J^0$ such that $k = sk'$.  Applying
$\sigma$ we find $k = s^{-1} k'$, so that $s=s^{-1}$ and $s\in D$.
Thus we have proved that $\Ksharp_J = D \cdot \Ksharp_J^0$, and the
claim follows.

Let us fix a face $F=F_{\Delta,J}$. By~\eqref{eq:noyauKpversNp},
$\Ksharp_J$ is a normal subgroup of $K_F$.  In fact the identity
component $\Scompz_\Delta$ of $\Scomp_\Delta$ centralises $\Ksharp_J$.
Indeed, $\forall s\in \Scomp_\Delta$, $\forall k \in \Ksharp_J$, we have
$k' = s\,k\,s^{-1} \in \Ksharp_J$. Since $\Ksharp_J \subset
K^\sigma$, we deduce $k' = \sigma(k') = s^{-1}\, k\, s$ and then $k =
s\,k'\,s^{-1} = s^2 \, k\, s^{-2}$. Since any element of the torus
$\Scompz_\Delta$ is a square, our claim follows.
Then~\eqref{eq:splitting} gives the exact
sequences~\eqref{eq:suitesexactesisotropies}, (in which $\Scompz_\Delta
\times \Ksharp_J$ is a direct product).  We deduce the exact
sequences~\eqref{eq:suitesexactescompisotropies}.
\begin{gather}
  \label{eq:suitesexactesisotropies}
  1\to D \to \Scomp_\Delta \ltimes \Ksharp_J \to K_F \to 1, \qquad
  1\to D_\Delta \to \Scompz_\Delta \times \Ksharp_J \to K_F \to 1,  \\
  \label{eq:suitesexactescompisotropies}
  D \to \tau_\Delta \ltimes \tau_J \to \tau_F \to 1, \qquad\quad
  D_\Delta \to \tau_J \to \tau_F \to 1.
\end{gather}
In particular, the groups $\tau_F$ (and then $\tau_{\O}$) are
quotients of the $\tau_J$, hence of the type $(\Z/2\Z)^a$, for some
$a\in\N$.

Let us consider more precisely the group $\tau_\Delta$.  Remember that
$\Scomp_\Delta$ is the stabiliser of a point $p\in C_{X_Z}\subset Z$
in $\Scomp = K\cap \Split$, which is the maximal compact subgroup of
$\Split$ (hence connected). Moreover $D= \{t\in \Split;\: t^2=1\}
\simeq (\Z/2\Z)^l$, $\Split$ acts on $Z$ via $\Splitbis =\Split/D$ and
$Z$ is toric for $\Splitbis$, so that $\Splitbis_p$ is connected. Let
$\Scompbis = \Scomp/D$ be the maximal compact subgroup of $\Splitbis$.
Then $\Scompbis_p$ is connected too and the exact sequences
$$
1\to D \to \Scomp_p \to \Scompbis_p \to 1, \qquad
1\to (\Z/2\Z)^{l'} \to \Scompz_p \to \Scompbis_p \to 1,
$$
where $l' = \dim \Scompz_p = |\Delta|$, show that $\tau_\Delta \simeq
(\Z/2\Z)^{l- |\Delta|}$ and $D_\Delta \simeq (\Z/2\Z)^{|\Delta|}$. For
$\Delta_1 \in \S$, such that $\Delta \subset \Delta_1$, we have
$\Scomp_\Delta \subset \Scomp_{\Delta_1}$ and a natural morphism
$\tau_\Delta \to \tau_{\Delta_1}$ which is surjective with kernel
$(\Z/2\Z)^{|\Delta_1| - |\Delta|}$.
\begin{lemma}
  \label{lem:extension}
  Let $\O_\Delta$ be a $G$-orbit of $X_Z$ and $L_\rho$ a
  $G$-equivariant local system on $\O_\Delta$, corresponding to a
  representation $\rho:\tau_{\O_\Delta} \to GL(V_\rho)$ of
  $\tau_{\O_\Delta}$.  We assume that $L_\rho$ is irreducible.  We set
  $F = F_{\Delta,J_\Delta}$ and let $i_F: \tau_\Delta \to \tau_F =
  \tau_{\O_\Delta}$ be the morphism induced
  by~\eqref{eq:suitesexactescompisotropies}.  For $v\in V\setminus
  \Delta$ such that $\Delta_1 = \Delta \sqcup \{v\} \in \S$, we have
  $\ker(\tau_\Delta \to \tau_{\Delta_1}) \simeq \Z/2\Z$. We let $s_v
  \in \tau_\Delta$ be the generator of this kernel.
  
  Then $\rho(i_F(s_v)) = \pm Id_{V_\rho}$ and this is the monodromy of
  $L_\rho$ around $D_v$.  If $\rho(i_F(s_v)) = Id_{V_\rho}$ then
  $\rho$ induces a representation, $\rho_1$, of $\tau_{\O_{\Delta_1}}$
  and $L_\rho$ extends to a local system, $L_1$, on $\O_\Delta \sqcup
  \O_{\Delta_1}$, such that $L_1|_{\O_{\Delta_1}}$ corresponds to
  $\rho_1$.
\end{lemma}
\begin{proof}
  Recall that the variety $X_Z$ comes with a morphism $\pi:X_Z \to X$.
  The $G$-orbits of $X$ are parameterised by subsets of
  $\{1,\ldots,l\}$ and $J_\Delta \subset \{1,\ldots,l\}$ is determined
  by: $\pi(\O_\Delta) =\O_{J_\Delta}$.  The hypothesis gives
  $\pi(\O_{\Delta_1}) \subset \ovl{\pi(\O_\Delta)}$, so that $J_\Delta
  \subset J_{\Delta_1}$.  We set $F_1 = F_{\Delta_1,J_{\Delta_1}}$.
  By~\eqref{intersectionUF} we deduce that $U_{F_1} \cap U_F =
  U_{F_2}$, where $F_2 = F_{\Delta,J_{\Delta_1}}$. We obtain the
  following commutative diagram
  $$
  \xymatrix@R=4mm{ E\times_K \vv_{F_{\Delta_1,J_{\Delta_1}}} &
    E\times_K \vv_{F_{\Delta,J_{\Delta_1}}} \ar@{_{(}->}[l]
    \ar@{^{(}->}[r] &
    E\times_K \vv_{F_{\Delta,J_\Delta}}   \\
    E\times_K \O_{\Delta_1} \ar@{}[u]|\cup \ar[d] & E\times_K
    (\O_\Delta \cap \vv_{F_{\Delta,J_{\Delta_1}}}) \ar@{_{(}->}[l]_-i
    \ar@{^{(}->}[r] \ar@{}[u]|\cup \ar[d] &
    E\times_K \O_\Delta \ar@{}[u]|\cup \ar[d]  \\
    E/K_{F_{\Delta_1,J_{\Delta_1}}} & E/K_{F_{\Delta,J_{\Delta_1}}}
    \ar@{>>}[l] \ar@{>>}[r] & E/K_{F_{\Delta,J_\Delta}}
    \ponctuation{,}  }
  $$
  where the vertical arrows are homotopy equivalences,
  by~\eqref{eq:homot_equi_orbit}.  Let us consider a ``small'' loop
  $\gamma$ in $E\times_K \O_\Delta$ around $E\times_K \O_{\Delta_1}$,
  as in section~\ref{normal crossings divisors}. Since $\gamma$ is
  small, we may assume that it is included in the neighbourhood
  $E\times_K \vv_{F_1}$ of $E\times_K \O_{\Delta_1}$.  Since $\gamma$
  doesn't meet $E\times_K D_v$ and $\vv_{F_1} \setminus D_v
  =\vv_{F_2}$, by~\eqref{eq:verifassumiii}, $\gamma$ is in fact
  contained in $E\times_K \vv_{F_2}$. Hence it represents a generator
  of the kernel of $\poinc(i)$, the map induced on the fundamental
  groups by the inclusion $i$ of the above diagram.
  
  Let $j:\tau_{F_2} \to \tau_{F_1}$ be the morphism induced by
  $\tau_\Delta \to \tau_{\Delta_1}$.  By~\eqref{eq:groupfond}, $j =
  \poinc(i)$ and by~\eqref{eq:suitesexactescompisotropies} we have the
  commutative diagram:
  $$
  \xymatrix@R=5mm{ D \ar[r]\ar@{=}[d] &\tau_{\Delta} \ltimes
    \tau_{J_{\Delta_1}} \ar[r]\ar[d]
    &\tau_{F_2} \ar[r]\ar[d]^j & 1 \ar@{=}[d] \\
    D \ar[r] &\tau_{\Delta_1} \ltimes \tau_{J_{\Delta_1}} \ar[r]
    &\tau_{F_1} \ar[r] & 1 \ponctuation{.} }
  $$
  We have already seen that $\ker(\tau_\Delta \to \tau_{\Delta_1})
  \simeq (\Z/2\Z)^{|\Delta_1| - |\Delta|} \simeq \Z/2\Z$.  The above
  diagram implies that $\ker(j)$ is $0$ or $\Z/2\Z$ and is
  generated by $i_F(s_v)$.  Since $s_v^2= 1$ and $V_\rho$ is
  irreducible it follows that $\rho(i_F(s_v)) = \pm Id_{V_\rho}$.
  
  We consider $L_\rho$ as well as a $G$-equivariant local system on
  $\O_\Delta$ or as a local system on $E\times_K \O_\Delta$.  Then
  $L_\rho|_{E\times_K \vv_{F_2}}$ corresponds to the representation
  $\rho_2$ of $\tau_{F_2}$ given by $j$ and $\rho$. The representation
  $\rho_2$ gives a representation $\rho_1$ of $\tau_{F_1} =
  \tau_{\O_{\Delta_1}}$ if and only if it sends $\ker(j)$ to
  $Id_{V_\rho}$, i.e. $\rho(i_F(s_v)) = Id_{V_\rho}$.  This also is
  equivalent to the fact that $L_\rho$ extends to $\O_\Delta \sqcup
  \O_{\Delta_1}$, with a restriction to $\O_{\Delta_1}$ corresponding
  to $\rho_1$.
\end{proof}
\begin{definition}
\label{def:defA}
We let $A$ be the set of pairs $\alpha=(\O,\rho)$, where $\O$ is a
$G$-orbit and $\rho:\tau_\O \to GL(V_\rho)$ an irreducible
representation of $\tau_\O$.  For $\alpha=(\O,\rho) \in A$, we let
$\Delta_\alpha \in \S$ be such that $\O = \O_{\Delta_\alpha}$ and,
with the notations of the previous lemma, we set 
$$
\Delta'_\alpha =
\{v\in V \setminus \Delta_\alpha;\: \Delta_\alpha \sqcup\{v\} \in \S
\: \text{and} \: \rho(i_F(s_v)) = -Id\},
\quad Z_\alpha = \smash{ \textstyle \ovl{\O} \setminus \bigcup_{v\in
    \Delta'_\alpha} D_v}.
$$
By the lemma, the $G$-equivariant local system on $\O$
corresponding to $\rho$ extends to a $G$-equivariant local system
$L^0_\alpha$ on $Z_\alpha$.  We extend it by $0$ outside $Z_\alpha$
(keeping the notation $L^0_\alpha$ also for the extension). We denote
by $L_\alpha$ the corresponding local system on $E\times_K X_Z$.
\end{definition}

\subsection{dg-algebras on the set of faces}
By section~\ref{eq_der_cat} we have the equivalences of categories:
$$
\D^b_{G,c}(X_Z) \simeq \D_G(X_Z) \langle L^0_\alpha, \,\alpha\in A
\rangle \simeq \D(E\times_K X_Z) \langle L_\alpha,
\,\alpha\in A \rangle.
$$
We are in the situation of
assumptions~\ref{assumptionstratification} with $Y = E\times_K X_Z$,
stratified by the set of faces $\F$, the subspaces $E\times_K D_v$ and
the local systems $L_\alpha$, $\alpha\in A$.  Conditions (i)--(iii)
were verified in section~\ref{stratificationfaces} and (iv), (v)
follow from definition~\ref{def:defA}.  Hence we may apply
proposition~\ref{prop:eq_cat1}: the category $\D_G(X_Z)\langle
L^0_\alpha \rangle$ is equivalent to a category of dg-modules over
$\F$, $\D_{\fC{}{}} \langle N_\alpha \rangle$, where $\fC{}{}$ is a
sheaf of dg-algebras on $\F$, whose description is recalled below, and
the $N_\alpha$ are $\fC{}{}$-modules.

First, for $\alpha = (\O,\rho), \beta = (\O',\rho') \in A$, we define
a sheaf $\fC{\alpha\beta}{}$ on $\F$ by its stalks at any face $F$.
We have:
$$
\pxf(\ovl{\O}) = \{ F_{\Delta,J}\in \F; \: \Delta_\alpha \subset
\Delta \},
\quad \text{and} \quad
\pxf(Z_\alpha) = \{
F_{\Delta,J}\in \F; \: \Delta_\alpha \subset \Delta \subset
(V\setminus \Delta'_\alpha) \}.
$$
We recall the notations~\ref{not:IetIprime} (with $\F_{\alpha\beta}$ instead
of $I_{\alpha\beta}$):
\begin{equation}
  \label{eq:deffabfprimeab}
\begin{split}
  \F_{\alpha\beta} &= \pxf(Z_\alpha \cap Z_\beta) = \{ F_{\Delta,J}\in
  \F; \: (\Delta_\alpha \cup \Delta_\beta) \subset \Delta \subset
  (V\setminus (\Delta'_\alpha \cup \Delta'_\beta)) \}, \\
  d_{\alpha\beta} &= | \Delta_\alpha \setminus \Delta_\beta |, \qquad
  \qquad \Delta'_{\aaa\bbb} = (\Delta'_\aaa \setminus \Delta'_\bbb)
  \cup
  (\Delta'_\bbb \setminus \Delta'_\aaa),  \\
  \F'_{\alpha\beta} &\textstyle = \pxf \bigl( Z_{\Delta_\aaa \cup
    \Delta_\bbb} \setminus (\bigcup_{v\in \Delta'_{\aaa\bbb}} D_v)
  \bigr)   \setminus \F_{\aaa\bbb}  \\
  &= \{ F_{\Delta,J}\in \F; \: (\Delta_\alpha \cup \Delta_\beta)
  \subset \Delta \subset (V\setminus \Delta'_{\alpha\beta}) \quad
  \text{and} \quad \Delta'_\alpha \cap \Delta'_\beta \cap \Delta \not=
  \emptyset\}.
\end{split}
\end{equation}
For $F\in \pxf(Z_\alpha)$, the restriction to $E\times_K K\cdot F$ of
the local system $L_\alpha$ has an extension to $E\times_K \vv_F$
(recall that $\vv_F =\pxf^{-1}(U_F)$ --
notations~\ref{not:stratifXZ}); we denote it by $L_{\alpha,F}$.
According to the defining formula~\eqref{eq:def_fB} we consider three
cases: (i) $F\in \F_{\alpha\beta}$, (ii) $F\in \F'_{\alpha\beta}$,
(iii) $F\not\in \F_{\alpha\beta} \sqcup \F'_{\alpha\beta}$.  For $F\in
\F_{\alpha\beta}$, we have
$$
\fC{\alpha\beta}{F} = \sect(E\times_K \vv_F; \Omega_{E\times_K
  \vv_F} \otimes \hom(L_{\alpha,F}, L_{\beta,F})) \,[-2d_{\alpha\beta}].
$$
Case (ii) is reduced to (i) as in~\ref{env:defin_fB}, and in case
(iii) we have $\fC{\alpha\beta}{F} =0$.  We set $\fC{}{} =
\oplus_{\alpha,\beta\in A} \fC{\alpha\beta}{}$.  

We denote by $\delta_v \in H^2_{K,D_v}(X_Z;\C_{X_Z})$ the
$K$-equivariant fundamental class of $D_v$ in $X_Z$.  We choose forms
$\xi_v \in \sect(E\times_K X_Z; \Omega_{E\times_K X_Z}^2)$
representing the $\delta_v$, and use them to define a product on
$\fC{}{}$, as in~\ref{env:defin_fB}, turning $\fC{}{}$ into a sheaf of
dg-algebras on $\F$.

The $\fC{}{}$-dg-module $N_\alpha$ is $N_\alpha =
\oplus_{\alpha',\alpha} \fC{\alpha'\alpha}{}$, with a
$\fC{}{}$-structure defined like the product of $\fC{}{}$. Let us set
$L=\oplus_{\alpha\in \epm} L_\alpha$.  By~\eqref{eq:section_de_A} and
lemma~\ref{lem:qisAB} (ii), we have, for a face $F\in \F$, the
isomorphism of algebras:
\begin{equation}
  \label{eq:germecohomo}
H^\cdot(\fC{}{F}) \simeq H^\cdot(\sect(U_F;\fC{}{})) \simeq
\Ext^\cdot_{\D_G(X_Z)}(L|_{V_F},L|_{V_F}).
\end{equation}
We also introduce the sheaf $\H$ on $\F$ given by the cohomology of
$\fC{}{}$, i.e. the sheaf associated to the presheaf $U\mapsto
H^\cdot(\sect(U;\fC{}{}))$. This is a sheaf of dg-algebras on $\F$,
with differential $0$. For a face $F\in \F$ we have $\sect(U_F;\H) =
\H_F = H^\cdot(\fC{}{F})$. We define in the same way the $\H$-module
$\H_\alpha$ associated to $U\mapsto H^\cdot(\sect(U;N_\alpha))$.

\section{Formality of the de Rham algebra}
\label{formality_of_de_Rham_algebra}
We keep the notations introduced in the previous section. Our aim is
to prove the following result.
\begin{proposition}
  \label{prop:formalite}
  There exists a sequence of quasi-isomorphisms of sheaves of
  dg-algebras on $\F$, $\fC{}{} \to \fcun{}{} \from \fcde{}{} \to
  \fctr{}{} \from \fcqu{}{} \to \fcci{}{} \simeq \H$, relating
  $\fC{}{}$ and $\H$. It induces an equivalence of categories between
  $\D^b_{G,c}(X_Z)$ and $\D_\H \langle \H_\alpha,\: \alpha\in A
  \rangle$.
\end{proposition}
We know from the previous section that $\D^b_{G,c}(X_Z)$ is equivalent
to $\D_{\fC{}{}} \langle N_\alpha \rangle$. Now a quasi-isomorphism
between sheaves of dg-algebras induces an equivalence between their
derived categories of dg-modules.  Hence the first part of the
proposition implies that $\D^b_{G,c}(X_Z)$ is equivalent to $\D_\H
\langle M_\alpha \rangle$, where $M_\alpha$ is the image of $N_\alpha$
by the chain of equivalences.  The remainder of this section is
devoted to the construction of a sequence of quasi-isomorphisms as in
the proposition.  It will follow from the construction that $M_\alpha$
is indeed isomorphic to $\H_\alpha$.

\subsection{Decomposition of the cohomology}
By theorem~\ref{thm:BDP}, the isotropy group $K_F$ almost decomposes
as a product, up to a finite subgroup. We deduce a decomposition for
$H^\cdot ( \fC{\alpha\beta}{F} )$.

For $F\in \F$, we have defined in~\eqref{eq:def_qF} a fibration $q_F:
E\times_K \vv_F \to E/K_F$, with contractible fibres $K_F\cdot \uu_F$.
In particular $q_F$ gives an isomorphism between the fundamental
groups. Hence, for $\alpha \in A$ such that $F \subset Z_\alpha$, the
local system $L_{\alpha,F}$ on $E\times_K \vv_F$ is the inverse image
of a local system $L'_{\alpha,F}$ on $E/K_F$.  Setting $M =
\hom(L'_{\alpha,F}, L'_{\beta,F})$, for another $\beta\in A$ with $F
\subset Z_\beta$, we obtain:
\begin{equation}
\label{eq:cohomCab}
\begin{split}
  H^\cdot ( \fC{\alpha\beta}{F} ) &\simeq \Ext^\cdot_{\D(E\times_K
    \vv_F)}
  (L_{\alpha,F}, L_{\beta,F})  \\
  &\simeq \Ext^\cdot_{\D(E/K_F)} (L'_{\alpha,F}, L'_{\beta,F}) \simeq
  H^\cdot ( E/K_F; M).
\end{split}
\end{equation}
The following lemma describes more precisely $H^\cdot ( E/K_F; M)$.

Let us introduce some notations. For a face $F=F_{\Delta,J}$, we
recall that $\Scompz_\Delta$ and $\Ksharp_J$ commute and we consider
the action of $\Scompz_\Delta \times \Ksharp_J$ on $E^3$ by
$(s,k)\cdot (e_1,e_2,e_3) = (sk \cdot e_1, s\cdot e_2, k\cdot e_3)$.
Let also $a_F$ be the group morphism $\Scompz_\Delta \times \Ksharp_J
\to K_F$, $(s,k)\mapsto sk$.  The first projection $E^3 \to E$ is
$a_F$-equivariant and induces the morphism $r^1_F$ below.  In view
of~\eqref{eq:suitesexactesisotropies}, $r^1_F$ is a fibration with
fibre $E^2/D_\Delta$, which is acyclic (i.e.  $H^0(E^2/D_\Delta;\C) =
\C$ and, for $i\not=0$, $H^i(E^2/D_\Delta;\C) = 0$).  The projection
to the last two factors $E^3 \to E^2$ induces in the same way the
morphism $r^2_F$ below, which is a fibration with acyclic fibre $E$.
\begin{equation}
  \label{eq:fibrationsEEE}
E/K_F  \xfrom[E^2/D_\Delta]{\quad r^1_F \quad} 
  E^3/(\Scompz_\Delta \times \Ksharp_J)
\xto[E]{\quad r^2_F \quad} (E/\Scompz_\Delta) \times (E/\Ksharp_J).
\end{equation}
\begin{lemma}
  \label{lem:deccohom}
  We consider a face $F = F_{\Delta,J} \in \F$, $\rho: \tau_F \to
  GL(V_\rho)$ a representation of $\tau_F$, and $M$ the local system
  on $E/K_F$ corresponding to $\rho$. We let $\rho_J:\tau_J \to
  GL(V_\rho)$ be the representation obtained from $\rho$ and the
  morphism $\tau_J \to \tau_F$. We let $M_J$ be the local system on
  $E/\Ksharp_J$ corresponding to $\rho_J$. Then
  $$
  H^\cdot(E/K_F; M) \simeq \C[X_v; \: v\in \Delta] \otimes
  H^\cdot(E/\Ksharp_J; M_J),
  $$
  where the $X_v$ are indeterminates of degree $2$.
\end{lemma}
\begin{proof}
  Since $r^1_F$ is a fibration with acyclic fibres, we have $M \simeq
  R(r^1_F)_* (r^1_F)^{-1} M$. Hence $H^\cdot(E/K_F; M) \simeq
  H^\cdot(E^3/(\Scomp_\Delta \times \Ksharp_J);(r^1_F)^{-1} M)$.  We
  have $\poinc(E/K_F) =\tau_F$, $\poinc(E^3/(\Scompz_\Delta \times
  \Ksharp_J)) =\tau_J$ and the morphism induced by $r^1_F$ on the
  fundamental groups is the morphism of the lemma $\tau_J \to \tau_F$.
  Hence $(r^1_F)^{-1} M$ is the local system corresponding to
  $\rho_J$.
  
  Since $r^2_F$ has a contractible fibre, it gives an isomorphism on
  the fundamental groups. Hence $(r^1_F)^{-1} M \simeq (r^2_F)^{-1}
  (\C_{E/\Scompz_\Delta} \boxtimes M_J)$. This also gives an isomorphism
  on the cohomology groups of $(r^1_F)^{-1} M$ and $\C_{E/\Scompz_\Delta}
  \boxtimes M_J$.  We conclude by the K\"unneth formula and the fact
  that $\Scompz_\Delta$ is the torus $(\C^*)^{| \Delta|}$, so that
  $H^\cdot(E/\Scompz_\Delta;\C) = H^\cdot_{\Scompz_\Delta}(\{pt\};\C)$ is a
  polynomial algebra in $|\Delta|$ variables.
\end{proof}
We describe the local systems $L_{\alpha,F}$ on $E\times_K \vv_F$,
$L'_{\alpha,F}$ on $E/K_F$ and $(L'_{\alpha,F})_J$ on
$\smash{E/\Ksharp_J}$, in terms of representations.

For $\alpha= (\O,\rho) \in A$, let $F_\alpha
=F_{\Delta_\alpha,J_{\Delta_\alpha}}$ be the closed face of $\O\cap
C_{X_Z}$ and $V_\rho$ be the representation space of $\rho$.  We have
seen that $\tau_\O = \tau_{F_\alpha}$, so that we have a morphism
$\tau_{J_{\Delta_\alpha}} \to \tau_\O$. Let $\O_\Delta$ be a $G$-orbit
such that $\O_\Delta \subset \ovl{\O}$ and with closed face
$F_{\Delta,J_\Delta}$.  We recall that $\pxzx(\O) =
\O_{J_{\Delta_\alpha}}$ and $\pxzx(\O_\Delta) = \O_{J_\Delta}$ (where
$\pxzx$ is the map from $X_Z$ to $X$). Hence $J_{\Delta_\alpha}
\subset J_\Delta$.  If $F=F_{\Delta,J}$ is another face of
$\O_\Delta$, we have $J_\Delta \subset J$.  Finally, for any face
$F=F_{\Delta,J}$ such that $F \subset \ovl{\O}$, we have
$J_{\Delta_\alpha} \subset J$, so that $\Ksharp_J \subset
\Ksharp_{J_{\Delta_\alpha}}$ and we obtain a group morphism:
\begin{equation}
  \label{eq:deftjo}
 \text{for $F_{\Delta,J} \subset \ovl{\O}$}, \qquad
 t_J^\O: \tau_J \to \tau_\O.
\end{equation}
We let $\rho_J$ be the representation of $\tau_J$ given by $V_\rho$
and $ t_J^\O$, and we let $L^1_{\alpha,F}$ be the corresponding local
system on $E/\Ksharp_J$.

Now we assume moreover that $F\subset Z_\alpha$. This means, by
lemma~\ref{lem:extension}, that $\rho$ induces a representation, say
$\rho'$, of $\tau_{\O_\Delta}$. Then the representation $\rho_J$ of
$\tau_J$ is given by $\rho'$ and the morphism $\tau_J \to
\tau_{J_\Delta} \to \tau_{\O_\Delta}$.  Since the morphism induced by
$r^1_F$ on the fundamental groups is $\tau_J \to \tau_F$, we obtain
the following relations:
\begin{equation}
  \label{eq:LLprimeL1}
L_{\alpha,F} \simeq q_F^{-1} (L'_{\alpha,F}), \qquad
(r^1_F)^{-1}(L'_{\alpha,F}) \simeq (r^2_F)^{-1}( \C_{E/\Scompz_\Delta}
\boxtimes L^1_{\alpha,F}).
\end{equation}
In lemma~\ref{lem:deccohom}, we have used
$H^\cdot_{\Scompz_\Delta}(\{pt\};\C) \simeq \C[X_v;\: v\in \Delta]$.  The
choice of indexing the indeterminates by $\Delta$ is not arbitrary, as
explained in the next lemma.

For $v\in V$, we have denoted by $\delta_v$ the $G$- (or $K$-)
equivariant fundamental class of $D_v$ in $X_Z$, $\delta_v \in
H^2_{K,D_v}(X_Z;\C)$. Let us also denote by $\delta_v$ its
``restriction'' to any $K$ stable open subset of $X_Z$.  The following
lemma describes the image of $\delta_v \in H^2_K(V_F;\C)$ by the
isomorphism $H^\cdot_K(V_F;\C) \simeq H^\cdot(E/K_F; \C)$ composed with
the isomorphism of lemma~\ref{lem:deccohom}.

Let us first recall the construction of the isomorphism
\begin{equation}
  \label{eq:isomHSalgsym}
H^\cdot(E/\Scompz_\Delta;\C) = H^\cdot_{\Scompz_\Delta}(\{pt\};\C) \simeq
Sym(Lie(\Scompz_\Delta)^*),
\end{equation}
where the elements of $Lie(\Scompz_\Delta)^*$ have degree $2$. A
character $\chi : \Scompz_\Delta \to \C^*$ gives an element $d_\chi \in
Lie(\Scompz_\Delta)^*$ by differentiation. It also gives a one dimensional
representation of $\Scompz_\Delta$, $\C_\chi$, and a line bundle $l_\chi =
E\times_{\Scompz_\Delta} \C_\chi$ over $E/\Scompz_\Delta$.  The above
isomorphism sends $d_\chi$ to the Chern class $c_2(l_\chi)$.  We note
that this Chern class is nothing but the $\Scompz_\Delta$-equivariant
fundamental class of $\{0\}$ in $\C_\chi$.
\begin{lemma}
  \label{lem:delatvXv}
  Let $v\in V$ and $F=F_{\Delta,J}$ be a face such that $F\subset
  D_v$.  For a point $p\in F$, $\Scompz_\Delta$ acts on $T_pX_Z / T_p D_v
  \simeq \C$. Let $\chi_v$ be the corresponding character of
  $\Scompz_\Delta$ and $X_v \in H^2_{\Scompz_\Delta}(\{pt\};\C)$ the
  associated equivariant class.  We have
  $$
  H^\cdot_K(V_F;\C) \simeq H^\cdot_{K_F}(\{pt\};\C) \simeq
  H^\cdot_{\Scompz_\Delta}(\{pt\};\C) \otimes
  H^\cdot_{\Ksharp_J}(\{pt\};\C)
  $$
  and this isomorphism sends $\delta_v$ to $X_v \otimes 1$.
  Moreover $H^\cdot_{\Scompz_\Delta}(\{pt\};\C) \simeq \C[X_v; \: v\in
  \Delta]$.
\end{lemma}
\begin{proof}
  We have seen that the first isomorphism follows from the homotopy
  equivalence $q_F$ (see~\eqref{eq:def_qF}). The second one is a
  special case of lemma~\ref{lem:deccohom}, with $M=\C_{E/K_F}$.
  
  We have an action of $K_F$ on $N_{p,v} = T_pX_Z / T_p D_v$ and
  natural isomorphisms:
  $$
  H^\cdot_K(V_F;\C) \isoto H^\cdot_{K_F}(N_{p,v};\C) \isoto
  H^\cdot_{K_F}(\{pt\};\C).
  $$
  The class $\delta_v \in H^\cdot_K(V_F;\C)$ can be identified with
  the $K_F$-equivariant fundamental class of $\{0\}$ in $N_{p,v}$,
  $\delta_{ \{0\}| N_{p,v}} \in H^2_{K_F, \{0\}}(N_{p,v};\C)$.  Hence
  its image by the natural morphism from $H^\cdot_{K_F}(\{pt\};\C)$ to
  $\smash{ H^\cdot_{\Scompz_\Delta}(\{pt\};\C)}$ is the
  $\Scompz_\Delta$-equivariant fundamental class of $\{0\}$ in $N_{p,v}$,
  i.e. $X_v$.  By~\eqref{eq:noyauKpversNp}, $\Ksharp_J$ acts
  trivially on $N_p$, so that the image of $\delta_v$ by
  $H^\cdot_{K_F}(\{pt\};\C) \to H^\cdot_{\Ksharp_J}(\{pt\};\C)$ is
  $0$. Since $H^2_K(V_F;\C)$ only has the components
  $H^2_{\Scompz_\Delta}(\{pt\};\C)$ and $H^2_{\Ksharp_J}(\{pt\};\C)$, we
  deduce that $\delta_v$ is sent to $X_v \otimes 1$, as claimed.

  Let us set $N_p = T_pX_Z / T_p \O_\Delta$. By theorem~\ref{thm:BDP},
  the kernel of $\Scompz_\Delta \to GL(N_p)$ is finite.  Since $N_p \simeq
  \oplus_{v\in\Delta} N_{p,v}$, it follows that the characters
  $\chi_v$, $v\in\Delta$, are independent. Since $\dim \Scompz_\Delta =
  |\Delta|$ we obtain the last assertion.
\end{proof}

\subsection{Decomposition of the dg-algebras}
We would like to decompose the de Rham complex $\fC{\alpha\beta}{F}$
as we have decomposed its cohomology in lemma~\ref{lem:deccohom}.
However in the sequence of fibrations,
$$
E\times_K \vv_F \xto{q_F} E/K_F  \xfrom{ r^1_F } 
  E^3/(\Scompz_\Delta \times \Ksharp_J)
\xto{ r^2_F } (E/\Scompz_\Delta) \times (E/\Ksharp_J),
$$
the morphism $r^1_F$ goes in the wrong direction, i.e. we have no
natural map from $\sect(E/\Ksharp_J; \Omega_{E/\Ksharp_J})$ to
$\sect(E\times_K \vv_F ;\Omega_{E\times_K \vv_F })$. Hence we first
replace $E\times_K \vv_F$ by the fibre product built on $q_F$ and
$r^1_F$.
\subsubsection{Pull-back to a fibre product}
\begin{lemma}
  \label{lem:defdevpF}
  We keep the notations $q_F$, $r^1_F$, $r^2_F$ defined
  in~\eqref{eq:def_qF} and~\eqref{eq:fibrationsEEE}. For a face
  $F=F_{\Delta,J} \in \F$, we set:
  $$
  \vp_F = (E\times_K \vv_F) \times_{E/K_F} (E^3/(\Scompz_\Delta \times
  \Ksharp_J)).
  $$
  \begin{itemize}
  \item [(i)] For any face $F=F_{\Delta,J} \in \F$, we have fibrations
$$
\nu_F:\vp_F \to E\times_K \vv_F, \qquad\quad 
r_F: \vp_F \to (E/\Scompz_\Delta) \times (E/\Ksharp_J),
$$
$\nu_F$ has fibres homeomorphic to $E^2/D_\Delta$ and $r_F$ has
contractible fibres.
\item[(ii)] For $F_i=F_{\Delta_i,J_i}$, $i=1,2$, with $F_1 \subset
  \ovl{F_2}$, we have a natural morphism $v_{F_1F_2}: \vp_{F_2} \to
  \vp_{F_1}$ and a commutative diagram
\begin{equation}
  \label{eq:diag_proje}
  \vcenter{\xymatrix@R=5mm{
E\times_K \vv_{F_2}  \ar[d]  &
\vp_{F_2} \ar[d] \ar[r]_-{r_{F_2}} \ar[l]^-{\nu_{F_2}}  & 
 (E/\Scompz_{\Delta_2}) \times (E/\Ksharp_{J_2})  \ar[d]  \\
E\times_K \vv_{F_1}  &
\vp_{F_1}  \ar[r]^-{r_{F_1}} \ar[l]_-{\nu_{F_1}} & 
(E/\Scompz_{\Delta_1}) \times  (E/\Ksharp_{J_1}) \ponctuation{.}
}}
\end{equation}
\item[(iii)] For a third face $F_3$ with $F_1 \subset \ovl{F_2} \subset
  \ovl{F_3}$, we have $v_{F_1F_3} = v_{F_1F_2} \circ v_{F_2F_3}$.
  \end{itemize}
\end{lemma}
\begin{proof}
  The proof is more or less tautological.  By definition, $\vp_F$
  comes with two fibrations, $\nu_F:\vp_F \to E\times_K \vv_F$, with
  fibre $E^2/D_\Delta$, and $\mu_F: \vp_F \to E^3/(\Scompz_\Delta \times
  \Ksharp_J)$, with contractible fibre $K_F\cdot \uu_F$.  We set $r_F
  = r^2_F \circ \mu_F$; since $\mu_F$ and $r^2_F$ are fibrations with
  contractible fibres, so is $r_F$. This gives (i).
  
  For (ii), we have the inclusions $K_{F_2} \subset K_{F_1}$,
  $\Ksharp_{J_2} \subset \Ksharp_{J_1}$, $\Scomp_{\Delta_2} \subset
  \Scomp_{\Delta_1}$ and they induce commutative squares of fibrations:
\begin{equation}
  \label{eq:carreqetr}
  \vcenter{\xymatrix@R=5mm{
E^3/ (\Scomp_{\Delta_2} \times \Ksharp_{J_2})  \ar[d]\ar[r]_-{r^1_{F_2}}
&  E/K_{F_2}  \ar[d]  \\
E^3/ (\Scomp_{\Delta_1} \times \Ksharp_{J_1}) \ar[r]^-{r^1_{F_1}} 
&  E/K_{F_1} 
}}
\qquad
 \vcenter{\xymatrix@R=6mm{ E\times_K \vv_{F_2} 
\ar@{^{(}->}[d] \ar[r]_{q_{F_2}} &
   E/K_{F_2}\ar@{>>}[d]   \\
  E\times_K \vv_{F_1}
  \ar[r]^{q_{F_1}}
& E/K_{F_1} \ponctuation{,} }}
\end{equation}
and a similar square corresponding to $r^2_{F_\cdot}$. Now (ii)
follows from these diagrams and the definitions. The proof of (iii) is
similar.
\end{proof}
{\bf Definition of $\fcun{}{}$}. Now we pull back the construction of
$\fC{}{}$ to the $\vp_F$.  For $\alpha=(\O,\rho) \in A$, and $F\in F$, we
set $L^+_{\alpha,F} = \nu_F^{-1}L_{\alpha,F}$.  For $v\in V$, we set
$\xi'_{v,F} = \nu_F^*(\xi_v|_{E\times_K \vv_F}) \in
\sect(\vp_F;\Omega^2_{\vp_F})$.  For two faces $F_1 \subset
\ovl{F_2}$, we have :
\begin{equation}
  \label{eq:xipLplus}
v_{F_1F_2}^{-1} (L^+_{\alpha,F_1}) = L^+_{\alpha,F_2}, 
\qquad \quad 
v_{F_1F_2}^* (\xi'_{v,F_1}) = \xi'_{v,F_2}.
\end{equation}
We introduce a sheaf $\fcun{}{}$ on $\F$, copying the definition of $\fA$
in~\ref{env:defin_fB}.  For $\alpha,\beta \in A$, we define the sheaf
$\fcun{\alpha\beta}{}$ by its stalks at a face $F\in
\F_{\alpha\beta}$:
$$
\fcun{\alpha\beta}{F} = \sect(\vp_F; \Omega_{\vp_F} \otimes
\hom(L^+_{\alpha,F}, L^+_{\beta,F})) \,[-2d_{\alpha\beta}],
$$
and we reduce the case $F\not\in \F_{\alpha\beta}$ to this one, as
in definition~\ref{env:defin_fB}.  The only difference is that the
restriction maps, say from $\fcun{\alpha\beta}{F_1}$ to
$\fcun{\alpha\beta}{F_2}$, for $F_1 \subset \ovl{F_2}$, are
induced by $v_{F_1F_2}: \vp_{F_2} \to \vp_{F_1}$, instead of the
inclusion $E\times_K \vv_{F_2} \subset E\times_K \vv_{F_1}$. This
gives a sheaf by (iii) of lemma~\ref{lem:defdevpF}.

We set $\fcun{}{} = \oplus_{(\alpha,\beta) \in A^2} \fcun{\alpha\beta}{}$
and endow it with an algebra structure as $\fC{}{}$: more precisely, we
define the product on the stalks at a given face $F$ by replacing the
$\xi_v$ in the definition of $\fC{}{}$ by the $\xi'_{v,F}$; the
compatibility of the product and the restriction maps follows
from~\eqref{eq:xipLplus}.

For a face $F$, the inverse image by $\nu_F$ induces a natural
morphism $\nu^*_F:\fC{}{F} \to \fcun{}{F}$. Since $\xi'_{v,F} =
\nu^*_F\xi_v$, we see that the $\nu^*_F$ are morphisms of dg-algebras.
The commutative squares in lemma~\ref{lem:defdevpF} imply that the
$\nu^*_F$ induce a morphism of sheaves on $\F$, say $\nu^*: \fC{}{}
\to \fcun{}{}$.  Let us verify that it is a quasi-isomorphism.  We have
$\hom(L^+_{\alpha,F}, L^+_{\beta,F}) \simeq \nu_F^{-1}
\hom(L_{\alpha,F}, L_{\beta,F})$.  Since $\nu_F$ is a fibration with
fibre $E^2/D_\Delta$, which is acyclic over $\C$, we have, for any
sheaf $L$ on $E\times_K \vv_F$, $R(\nu_F)_* \nu_F^{-1}L \simeq L$.
Hence $H^\cdot (\fcun{\alpha\beta}{F}) = H^\cdot
(\fC{\alpha\beta}{F})$, as required.

\subsubsection{Decomposition}
For $\alpha=(\O,\rho) \in A$ and $F=F_{\Delta,J} \in \F$ with
$F\subset Z_\alpha$, let $L^1_{\alpha,F}$ be the local system on
$(E/\Ksharp_J)$, corresponding to the representation of $\tau_J$
given by $\rho$ and $t_J^\O:\tau_J \to \tau_\O$ (see~\eqref{eq:deftjo}
and after).  Then, by~\eqref{eq:LLprimeL1}, $L^+_{\alpha,F}\simeq
r_F^{-1} (\C_{E/\Scompz_\Delta} \boxtimes L^1_{\alpha,F})$.

For $\alpha,\beta \in A$ and $F\in \F_{\alpha\beta}$, we have, by
lemma~\ref{lem:deccohom}:
$$
H^\cdot (\fcun{\alpha\beta}{F}) \simeq 
H^\cdot(E/\Scompz_\Delta;
\C_{E/\Scompz_\Delta} ) \otimes H^\cdot(E/\Ksharp_J;\hom(L^1_{\alpha,F},
L^1_{\beta,F}) ).
$$
This isomorphism corresponds to a quasi-isomorphism at the level of
de Rham complexes. The product of forms, composed with the inverse
image by $r_F$ gives a quasi-isomorphism:
\begin{equation}
  \label{eq:presquedecompose}
\sect(E/\Scompz_\Delta; \Omega_{E/\Scompz_\Delta} ) \otimes \sect(E/\Ksharp_J;
\Omega_{E/\Ksharp_J} \otimes \hom(L^1_{\alpha,F}, L^1_{\beta,F}) )
\xto{qis} \fcun{\alpha\beta}{F} .
\end{equation}
However the forms $\xi'_{v,F} \in \smash{\Omega^2_{\vp_F}}$ do not
have to be pull-backs of forms by $r_F$ and we have no natural algebra
structure on the sum over $(\alpha,\beta)\in A^2$ of the groups
appearing in the left hand side of~\eqref{eq:presquedecompose}.  For
this we will replace the factor $\sect(E/\Scompz_\Delta;
\Omega_{E/\Scompz_\Delta} )$ by a free anti-commutative algebra
quasi-isomorphic to it.  We will define a sheaf $\fcde{}{}$ on $\F$
(see~\eqref{eq:deffC3} below) as the product of two sheaves: $\fcs{}{}$,
quasi-isomorphic to the de Rham algebra of $E/\Scompz_\Delta$
appearing in~\eqref{eq:presquedecompose}, and $\fck{}{}$, given by the
twisted de Rham complex on $E/\Ksharp_J$.

For $\Delta \subset V$, we introduce the dg-algebras $A(\Delta)$,
$B(\Delta)$ below, which are free anti-commu\-ta\-tive algebras, and a
quasi-isomorphism, $b(\Delta)$, between them. (The reason for
introducing $B(\Delta)$ is to be able to define morphism $f$ below,
which would be impossible with $A(\Delta)$ instead of $B(\Delta)$.)
\begin{gather}
  \notag 
  A(\Delta) =\C[X_v;\: v\in\Delta], \quad B(\Delta) = \C[X_v, Y_w;\:
  v\in V,\: w\in V\prive\Delta],  \\
  \label{eq:defalphadelta}
  \deg X_v =2, \quad \deg Y_w =1,\qquad dX_v = 0, \quad dY_w = X_w,  \\
  \notag 
  b(\Delta):B(\Delta) \to A(\Delta), \qquad \forall v\in \Delta,
  \: X_v \mapsto X_v, \quad \forall w\in V\setminus \Delta, \: X_w
  \mapsto 0, \: Y_w \mapsto 0.
\end{gather}
For $(\alpha,\beta) \in A^2$, and a face $F=F_{\Delta,J} \in
\F_{\alpha\beta}$ we set:
$$
\fcs{\alpha\beta}{F} = B(\Delta) \,[-2d_{\alpha\beta}].
$$
For other faces, we reduce to the case $F \in \F_{\alpha\beta}$,
as in definition~\ref{env:defin_fB}.  The restriction maps are given by
the inclusions $B(\Delta_1) \subset B(\Delta_2)$ if $\Delta_2 \subset
\Delta_1$.  We have a product similar to the product in $\fC{}{}$, as
follows.  For $\alpha, \beta, \gamma \in A$, we set
$\varepsilon_{\alpha\beta\gamma} = \prod_{v\in \nabla} X_v$, where
$\nabla$ is defined in~\eqref{eq:nabla}.  We define
$m_S^{\alpha\beta\gamma}: \fcs{\alpha\beta}{} \otimes
\fcs{\beta\gamma}{} \to \fcs{\alpha\gamma}{}$, $P\otimes Q \mapsto
\varepsilon_{\alpha\beta\gamma}PQ$.

In a similar way, for $(\alpha,\beta) \in A^2$, and a face
$F=F_{\Delta,J} \in \F_{\alpha\beta}$ we set:
$$
\fck{\alpha\beta}{F} = \sect(E/\Ksharp_J;
\Omega_{E/\Ksharp_J}\otimes \hom(L^1_{\alpha,F}, L^1_{\beta,F})).
$$
For other faces, we reduce to the case $F \in \F_{\alpha\beta}$, as
in definition~\ref{env:defin_fB}. The restriction maps are the
following ones: for $U_{F_2} \subset U_{F_1}$ we have a map $r_{12}:
E/\Ksharp_{J_2} \to E/\Ksharp_{J_1}$ and $L^1_{\alpha,F_2} =
r_{12}^{-1}L^1_{\alpha,F_1}$, hence an inverse image morphism
$\fck{\alpha\beta}{F_1} \to \fck{\alpha\beta}{F_2}$.  We also have an
obvious product $m_K^{\alpha\beta\gamma}: \fck{\alpha\beta}{} \otimes
\fck{\beta\gamma}{} \to \fck{\alpha\gamma}{}$ given by the product of
forms and the composition of morphisms.

{\bf Definition of $\fcde{}{}$}.
Now we set:
\begin{equation}
  \label{eq:deffC3}
\forall (\alpha,\beta) \in A^2,\quad
\fcde{\alpha\beta}{} =
\fcs{\alpha\beta}{} \otimes \fck{\alpha\beta}{},
\qquad
\fcde{}{} = \oplus_{(\alpha,\beta) \in A^2} \fcde{\alpha\beta}{}
\end{equation}
We let $m^{\alpha\beta\gamma}$ be the tensor product of
$m_S^{\alpha\beta\gamma}$ and $m_K^{\alpha\beta\gamma}$.  Since, by
definition, the $\varepsilon_{\alpha\beta\gamma}$ satisfy the same
identity~\eqref{eq:identiteensembliste} as the
$\eta_{\alpha\beta\gamma}$, we obtain a product on $\fcde{}{}$ defined
by the sum of the $m^{\alpha\beta\gamma}$.

{\bf Definition of $f: \fcde{}{} \to \fcun{}{}$}.  For $F\in \F$, we
note that $\fcun{}{F}$ is an algebra over
$\sect(\vp_F;\Omega_{\vp_F})$.  We define a morphism
$f^{\alpha\beta}_F:\fcde{\alpha\beta}{F} \to \fcun{\alpha\beta}{F}$ as
the product of $f^{S\alpha\beta}_F:\fcs{\alpha\beta}{F} \to
\sect(\vp_F;\Omega_{\vp_F})$ and
$f^{K\alpha\beta}_F:\fck{\alpha\beta}{F} \to \fcun{\alpha\beta}{F}$,
which are obtained as follows.

For a face $F$, let $s_F: \vp_F \to E/\Ksharp_J$ be the composition of
$r_F$ and the projection to $E/\Ksharp_J$.  For a form $\sigma$ and a
sheaves endomorphism $u$, we set $f^{K\alpha\beta}_F (\sigma\otimes u)
=s_F^*(\sigma) \otimes s_F^{-1}(u)$.

Now we define $f^{S\alpha\beta}_F$.  For $v\in V$, the fundamental
class, $\delta_v$, of $D_v$ in $X_Z$ restricts to $0$ on $X_Z
\setminus D_v$.  Hence the restriction of $\xi_v$ on $E\times_K
(X_Z\setminus D_v)$ is a boundary. Let us choose a form $\zeta_v \in
\sect( E\times_K (X_Z\setminus D_v); \Omega^1_{E\times_K X_Z})$ such
that, on $E\times_K (X_Z\setminus D_v)$ we have $\xi_v = d\zeta_v$.
We set $\zeta'_{v,F} = \nu_F^*(\zeta_v|_{E\times_K \vv_F})$ (we note
that, for $F=F_{\Delta,J}$ and $v\in V$ such that $v\not\in \Delta$,
we have $\vv_F \cap D_v =\emptyset$, so that $\zeta_v$ is defined on
$E\times_K \vv_F$). We set:
$$
\forall v\in V, \:\: f^{S\alpha\beta}_F(X_v) = \xi'_{v,F}, \qquad
\forall w \in V\setminus \Delta, \:\: f^{S\alpha\beta}_F(Y_w) =
\zeta'_{w,F}.
$$
Using lemma~\ref{lem:defdevpF}, one checks that
$f^{S\alpha\beta}_F$ and $f^{K\alpha\beta}_F$ give morphisms of
sheaves, say $f^{S\alpha\beta}$ and $f^{K\alpha\beta}$.  We define
$f^{\alpha\beta} =f^{S\alpha\beta} \otimes f^{K\alpha\beta}$ and $f =
\oplus f^{\alpha\beta}$.  In view of the definitions of the product
and the differentials, $f$ is a morphism of dg-algebras.  By
lemmas~\ref{lem:deccohom} and~\ref{lem:delatvXv} it is a
quasi-isomorphism.

\subsection{Formality of the toric part}
We define a sheaf quasi-isomorphic to $\fcs{}{}$ but with differential
$0$ as follows. For $(\alpha,\beta) \in A^2$, and a face
$F=F_{\Delta,J} \in \F_{\alpha\beta}$, we set:
\begin{equation}
  \label{eq:defCS1}
  \fcsun{\alpha\beta}{F} = A(\Delta) \,[-2d_{\alpha\beta}].
\end{equation}
with the following restriction maps. For two faces
$F_i=F_{\Delta_i,J_i}$, $i=1,2$, such that $F_1 \subset \ovl{F_2}$
(i.e. $U_{F_2}\subset U_{F_1}$), we have $\Delta_2 \subset \Delta_1$,
and the restriction $A(\Delta_1) \to A(\Delta_2)$ sends $X_v$ to $X_v$
for $v\in \Delta_2$ and to $0$ for $v\in \Delta_1 \setminus \Delta_2$.
As for $\fcs{}{}$, other faces are reduced to this case.

We define a product, $m_T^{\alpha\beta\gamma}$, similar to the product
of $\fcs{}{}$.  For $\alpha, \beta, \gamma \in A$, we let
$e_{\alpha\beta\gamma}$ be the section of $\fcsun{\alpha\gamma}{}$
defined by $(e_{\alpha\beta\gamma})_{F_{\Delta,J}} = b(\Delta)
(\varepsilon_{\alpha\beta\gamma})$.  We set $m_T^{\alpha\beta\gamma}:
\fcsun{\alpha\beta}{} \otimes \fcsun{\beta\gamma}{} \to
\fcsun{\alpha\gamma}{}$, $P\otimes Q \mapsto e_{\alpha\beta\gamma}PQ$.

{\bf Definition of $\fctr{}{}$}. We set as in~\eqref{eq:deffC3}:
\begin{equation}
  \label{eq:deffC4}
\forall (\alpha,\beta) \in A^2,\quad
\fctr{\alpha\beta}{} =
\fcsun{\alpha\beta}{} \otimes \fck{\alpha\beta}{},
\qquad
\fctr{}{} = \oplus_{(\alpha,\beta) \in A^2} \fctr{\alpha\beta}{}
\end{equation}
We let $m_C^{\alpha\beta\gamma}$ be the tensor product of
$m_T^{\alpha\beta\gamma}$ and $m_K^{\alpha\beta\gamma}$.  The sum
$m_C= \oplus m_C^{\alpha\beta\gamma}$ defines a product on
$\fctr{}{}$.  The $b(\Delta)$ defined in~\eqref{eq:defalphadelta} give
a quasi-isomorphism of sheaves of dg-algebras $\fcs{}{} \to
\fcsun{}{}$.  This induces a quasi-isomorphism of sheaves of
dg-algebras $\fcde{}{} \to \fctr{}{}$.

\subsection{Formality of $\fck{}{}$ and conclusion}
It remains to prove that the factor $\fck{}{}$ of $\fcde{}{}$ also is
quasi-isomorphic to its cohomology. For this we first give a more
handy expression for $\fck{\alpha\beta}{F}$ (see
formula~\eqref{eq:isom_K} below).

For $J \subset \{1,\ldots,l\}$, let $\pi_J:E/\Ksharp_J^0 \to E/
\Ksharp_J$ be the covering map, with group $\tau_J$, and set $A_J=
(\pi_J)_*(\C_{E/\Ksharp_J^0})$.  Then $A_J$ is a local system on
$E/\Ksharp_J$, considered as a right module over $\C[\tau_J]$,
locally free of rank one. It is also a sheaf of algebras. For $\alpha
= (\O,\rho) \in A$, $F =F_{\Delta,J} \in \F$, such that $F\subset
Z_\alpha$, we have, by~\eqref{eq:deftjo}:
$$
L^1_{\alpha,F} \simeq A_{J} \otimes_{\C[\tau_{J}]} V_\rho.
$$
Hence, for another element $\beta = (\O',\rho') \in A$, with
$F\subset Z_\beta$, we have:
$$
\hom(L^1_{\alpha,F}, L^1_{\beta,F}) \simeq \hom( A_{J}
\otimes_{\C[\tau_{J}]} V_\rho, A_{J} \otimes_{\C[\tau_{J}]}
V_{\rho'}).
$$
We want to ``factorise'' the local systems $A_J$ and the
representation spaces $V_\rho$, $V_{\rho'}$ in this last formula. We
will use the following definition (see~\cite{So01}).
\begin{definition}
  For a group $W$ and a $\C$-algebra $A$ with a right $W$-action by
  algebra automorphisms, we set $A^t[W] = A\otimes \C[W]$ with the
  product, for $a\in A$, $w\in W$, $(a\otimes w)\cdot (a'\otimes w') =
  (a (a'\cdot w))\otimes (w'w)$.  We have a natural embedding
  $\C[W]^{op} \hookrightarrow A^t[W]$, $w\mapsto a_w =1\otimes w$.  It
  induces a structure of right $\C[W] \otimes \C[W]^{op}$-module on
  $A^t[W]$: for $x\in A^t[W]$, $w,w'\in W$, $x\cdot (w\otimes w') =
  a_{w} \, x \, a_{w'}$.
\end{definition}
We set $B_J = \C[\tau_{J}] \otimes \C[\tau_{J}]^{op}$.  We consider
$\hom(A_J,A_J)$ as a sheaf of algebras, where the product is the
composition, and right $B_J$-module by $(\phi \cdot (t\otimes t') )
(a) = (\phi (a\cdot t') )\cdot t$, for $\phi$ a section of
$\hom(A_J,A_J)$, $a$ a section of $A_J$ and $t, t' \in \tau_J$.  With
these definitions, one can check that the map
\begin{equation}
  \label{eq:homAA}
A^t_J[\tau_J] \to \hom(A_J,A_J), \qquad a\otimes t \mapsto (\alpha
\mapsto a (\alpha \cdot t)), 
\end{equation}
where $a$, $\alpha$ are sections of $A_J$ and $t\in \tau_J$, is an
isomorphism of (sheaves of) algebras and right $B_J$-modules.
\begin{lemma}
  For $J\subset \{1,\ldots,l\}$, $\alpha, \beta \in A$, $F =
  F_{\Delta,J} \in \F_{\alpha\beta}$, we have, with the above
  notations:
  \begin{gather}
    \label{eq:homL2alphaL2beta}
    \hom(L^1_{\alpha,F}, L^1_{\beta,F}) \simeq \hom(A_J,A_J)
    \otimes_{B_J} \Hom(V_\rho, V_{\rho'}), \\
    \label{eq:isom_K}
    \fck{\alpha\beta}{F} \simeq \sect(E/\Ksharp_J^0;
    \Omega_{E/\Ksharp_J^0})^t[\tau_J] \otimes_{B_J} \Hom(V_\rho,
    V_{\rho'}).
  \end{gather}
\end{lemma}
\begin{proof}
  Let us prove~\eqref{eq:homL2alphaL2beta}. We will use the following
  fact. Set $R=\C[W]$ for a finite group $W$ and consider left and
  right $R$-modules, $M$, $N$, of finite ranks.  We have a right
  $R$-module structure on $M^* = \Hom_{\C}(M,\C)$ and a left one on
  $N^*$. Then the composition of (vector spaces) morphisms
\begin{equation}
  \label{eq:prod_dual=dual_prod}
  (N\otimes_R M)^* \to (N\otimes_\C M)^*  \isofrom 
  M^*\otimes_\C N^* \to M^* \otimes_R N^*
\end{equation}
is a canonical isomorphism. Indeed, it is compatible with direct sum,
and any $R$-module is semi-simple.  For $N$, $M$ irreducible with
$N\not= M^*$, both $N\otimes_R M$ and $M^* \otimes_R N^*$ are $0$. For
$N = M^*$, both $N\otimes_R M$ and $M^* \otimes_R N^*$ are canonically
identified with $\C$ by the duality contractions $N\otimes M\to \C$
and $M^* \otimes N^*\to \C$. Since the three morphisms
in~\eqref{eq:prod_dual=dual_prod} commute with the duality contraction
their composition corresponds to $id_\C$. Using this we deduce
canonical isomorphisms for left and right $R$-modules, $M_i$, $N_i$,
$i=1$, $2$:
\begin{align}
   \notag  \Hom_\C(N_1\otimes_R M_1, N_2\otimes_R M_2)
  &\simeq (M_1^*\otimes_R N_1^*) \otimes_\C (N_2\otimes_R M_2)  \\
   \label{hom_et_tens}  &\simeq (N_1^*\otimes_\C N_2) \otimes_{R\otimes R^{op}}
   (M_1^*\otimes_\C M_2)  \\
   \notag &\simeq \Hom_\C(N_1,N_2) \otimes_{R\otimes R^{op}}
   \Hom_\C(M_1,M_2).
\end{align}
Since isomorphism~\eqref{hom_et_tens} is canonical, it works as well
for sheaves and we obtain~\eqref{eq:homL2alphaL2beta}.

Now we deduce the second isomorphism.  We remark that $\Omega_{E/
  \Ksharp_J} \otimes A_J$ is isomorphic to $(\pi_J)_*\Omega_{E/
  \Ksharp_J^0}$. This isomorphism respects the $\tau_J$-action and
we have, by isomorphisms~\eqref{eq:homL2alphaL2beta}
and~\eqref{eq:homAA}:
$$
\Omega_{E/ \Ksharp_J} \otimes \hom(L^1_{\alpha,F}, L^1_{\beta,F}) 
\simeq
((\pi_J)_*\Omega_{E/ \Ksharp_J^0} )^t[\tau_J] 
\otimes_{B_J} \Hom(V_\rho, V_{\rho'}).
$$
Since $\otimes_{B_J}$ is exact for $B_J$-modules, the constant
sheaf $\Hom(V_\rho, V_{\rho'})$ factors out when we take global
sections, and we obtain~\eqref{eq:isom_K}.
\end{proof}
\begin{remark}
  \label{rem:retrictionCK}
  On the right hand side of formula~\eqref{eq:isom_K}, the restriction
  maps are given as follows. For another face $F'=F_{\Delta',J'} \in
  \F_{\alpha\beta}$ such that $F\subset \ovl{F'}$, we have
  $\Ksharp_{J'} \subset \Ksharp_J$, hence a groups morphism
  $a:\tau_{J'} \to \tau_J$ and a quotient map $p:E/\Ksharp_{J'}^0
  \to E/\Ksharp_J^0$. The inverse image of forms by $p$ is
  compatible with the action of $\tau_J$, $\tau_{J'}$ (via $a$) and we
  obtain a dg-algebras morphism:
  $$
  b:\sect(E/\Ksharp_J^0; \Omega_{E/\Ksharp_J^0})^t[\tau_J]
  \to \sect(E/\Ksharp_{J'}^0;
  \Omega_{E/\Ksharp_{\smash{J'}}^0})^t[\tau_{J'}] \otimes_{B_{J'}} B_J,
  $$
  by $b(\omega \otimes t) = (p^*(\omega) \otimes 1) \otimes (1
  \otimes t)$, for a form $\omega$ and $t\in \tau_{J}$.  Tensorisation
  with $\Hom(V_\rho, V_{\rho'})$ gives the desired restriction map
  from $\fck{\alpha\beta}{F}$ to $\fck{\alpha\beta}{F'}$.
\end{remark}
Let us explain how to recover the product $\fck{\alpha\beta}{}
\times \fck{\beta\gamma}{} \to \fck{\alpha\gamma}{}$ in the right
hand side of~\eqref{eq:isom_K}.  First we consider
isomorphism~\eqref{hom_et_tens}.  For $u:N_1\to N_2$, $v:M_1\to M_2$,
let us denote by $u \cdot v: N_1\otimes_R M_1 \to N_2\otimes_R M_2$
the image of $u\otimes v$ by~\eqref{hom_et_tens}.  For left and right
$R$-modules, $M_3$, $N_3$, $u':N_2\to N_3$, $v':M_2\to M_3$, we see
that $(u'\cdot v') \circ (u \cdot v) = (u'\circ u) \cdot (v'\circ v)$.

Hence the product in the right hand side of~\eqref{eq:homL2alphaL2beta}
is given by the composition in $\hom(A_J,A_J)$ and in
$\Hom(V_\rho, V_{\rho'})$. When we replace $\hom(A_J,A_J)$ by another
algebra, say $R_1$, as in~\eqref{eq:isom_K}, we use the following lemma
(with $R_2= \Hom(V, V)$, $V$ being the sum of irreducible representations
of $W$ and $\phi_2$ the action of $W$ on $V$).
\begin{lemma}
  \label{lem:produittordu}
  Let $W$ be a finite group, $R=\C[W]$ and $R_i$, $i=1,2$, algebras.
  Let $\phi_1:R^{op} \to R_1$, $\phi_2:R \to R_2$ be algebras
  morphisms. We consider the right $R\otimes R^{op}$-module structure
  on $R_1$ given by $a\cdot (w\otimes w') = \phi_1(w) a \phi_1(w')$,
  for $w,w'\in W$, $a\in R_1$, and the left structure on $R_2$,
  $(w\otimes w') \cdot a' = \phi_2(w) a' \phi_2(w')$.
  
  Then the formula $(a\otimes a', b\otimes b') \mapsto ab \otimes
  a'b'$, for $a,b \in R_1$, $a',b' \in R_2$, gives a well-defined
  product on $(R_1 \otimes_{R\otimes R^{op}} R_2)$.
\end{lemma}
\begin{proof}
  By symmetry, it is sufficient to prove, for $a,b\in R_1$, $a',b'\in
  R_2$, $w\, w'\in W$: $( \phi_1(w) a \phi_1(w') b) \otimes a'b' = ab
  \otimes (\phi_2(w) a' \phi_2(w') b')$. The tensor product is over
  $R\otimes R^{op}$, so that we are reduced to
  \begin{equation}
    \label{eq:prodtensabw}
    (a \phi_1(w') b) \otimes
  a'b' = ab \otimes ( a' \phi_2(w') b').
  \end{equation}
  Let us consider the subgroup $W'\subset W$ generated by $w'$, and
  $R' = \C[W']$. Then it is sufficient to prove
  that~\eqref{eq:prodtensabw} holds with a tensor product over
  $R'\otimes R'^{op}$. Hence, replacing $W$ by $W'$, we may assume
  from the beginning that $W = \langle w' \rangle$ is commutative. We
  decompose $R_i$ under the action of $w'$: $R_i = \oplus_{\lambda,\mu
    \in \C} R_i^{\lambda\mu}$, where $\forall a \in R_i^{\lambda\mu}$,
  $\phi_i(w') a = \lambda a$, $a \phi_i(w') = \mu a$.  By additivity,
  we may assume that $a,b,a',b'$ are elements of some
  $R_i^{\lambda\mu}$.  But $R_1^{\lambda\mu} \otimes_{R'\otimes
    R'^{op}} R_2^{\lambda'\mu'}$ is $0$ unless $\lambda = \lambda'$
  and $\mu = \mu'$. Similarly, for $a \in R_1^{\lambda\mu}$, $b\in
  R_1^{\mu'\nu}$, the product $ab$ is $0$ unless $\mu\mu' =1$.  Hence
  we may assume $a\in R_1^{\lambda\mu}$, $b \in R_1^{\mu^{-1}\nu}$,
  $a'\in R_2^{\lambda\mu}$, $b' \in R_2^{\mu^{-1}\nu}$.  In this case
  formula~\eqref{eq:prodtensabw} is obvious.
\end{proof}

\subsubsection{Formality of $\fck{}{}$}
Using~\eqref{eq:isom_K} we will deduce the formality of $\fck{}{}$ from
the formality of the de Rham algebras of $E/\Ksharp_J^0$ obtained
in lemma~\ref{lem:formalisersousgroupes}.

Let us denote by $\liek$, $\lieKsharp_J$, the Lie algebras of $K$,
$\Ksharp_J$, for $J\subset \{1,\ldots,l\}$.  Let us set for short
$H_J^0 = (S^\cdot(\lieKsharp_J))^{\Ksharp_J^0}$, viewed as a
dg-algebra with differential $0$; it is isomorphic to the
$\Ksharp_J^0$-equivariant cohomology algebra of the point,
$H_{\smash{\Ksharp_J^0}}^\cdot(\{pt\})$.

By lemma~\ref{lem:formalisersousgroupes}, the choice of a connection on
$\sect(E;\Omega_E)$ (in the sense of~\eqref{eq:defconnection})
gives functorial quasi-isomorphisms:
\begin{equation*}
  \sect(E/\Ksharp_J^0; \Omega_{E/\Ksharp_J^0})
\xfrom{f_J}  W(\liek)_{\lieKsharp_J -b}
\xto{g_J} H_J^0.
\end{equation*}
The normaliser $N_K(\Ksharp_J^0)$ acts on the above dg-algebras and
$\Ksharp_J^0$ acts trivially, so that $\tau_J = \Ksharp_J /
\Ksharp_J^0$ also acts. The morphisms $f_J$ and $g_J$ are
$\tau_J$-equivariant and yield quasi-isomorphisms:
\begin{equation}
  \label{eq:formaliteKJdiese}
  \sect(E/\Ksharp_J^0; \Omega_{E/\Ksharp_J^0})^t[\tau_J]
\xfrom{f'_J}  (W(\liek)_{\lieKsharp_J -b})^t[\tau_J]
\xto{g'_J} (H_J^0)^t[\tau_J].
\end{equation}

{\bf Definition of $\fcqu{}{}$}.  For $(\alpha,\beta) \in A^2$, we define
$\fckun{\alpha\beta}{}$ similarly as $\fck{\alpha\beta}{}$,
replacing the forms over $E/\Ksharp_J^0$ in
expression~\eqref{eq:isom_K} by a quasi-isomorphic dg-algebra. The
stalks at a face $F=F_{\Delta,J} \in \F_{\alpha\beta}$ are given by:
$$
\fckun{\alpha\beta}{F} = ( W(\liek)_{\lieKsharp_J -b}
)^t[\tau_J] \otimes_{B_J} \Hom(V_\rho, V_{\rho'}),
$$
with restriction morphisms defined as in
remark~\ref{rem:retrictionCK}, and product as in
lemma~\ref{lem:produittordu}.  As for $\fck{}{}$, other faces are reduced
to this case. In view of the isomorphism~\eqref{eq:isom_K} and the
quasi-isomorphism $f'_J$ of~\eqref{eq:formaliteKJdiese}, we have a
quasi-isomorphism of sheaves $\fckun{\alpha\beta}{} \to
\fck{\alpha\beta}{}$.  Setting
$$
\fcqu{\alpha\beta}{} = \fcsun{\alpha\beta}{} \otimes
\fckun{\alpha\beta}{} , \qquad \fcqu{}{} =\oplus_{(\alpha,\beta) \in
  A^2} \fcqu{\alpha\beta}{},
$$
with the product defined as the product of $\fctr{}{}$, we deduce a
quasi-isomorphism of sheaves of dg-algebras $\fcqu{}{} \to \fctr{}{}$.

{\bf Definition of $\fcci{}{}$}.
  We define $\fckde{\alpha\beta}{}$ similarly
as $\fckun{\alpha\beta}{}$, with stalks at a face $F=F_{\Delta,J}
\in \F_{\alpha\beta}$:
\begin{equation}
  \label{eq:deffCK2}
\fckde{\alpha\beta}{F} = (H_J^0)^t[\tau_J]
\otimes_{B_J} \Hom(V_\rho, V_{\rho'}).
\end{equation}
Then the $g'_J$ induce a quasi-isomorphism $\fckun{\alpha\beta}{} \to
\fckde{\alpha\beta}{}$ and, setting
$$
\fcci{\alpha\beta}{} = \fcsun{\alpha\beta}{} \otimes
\fckde{\alpha\beta}{} , \qquad \fcci{}{} =\oplus_{(\alpha,\beta) \in
  A^2} \fcci{\alpha\beta}{},
$$
we obtain a quasi-isomorphism of sheaves of dg-algebras $\fcqu{}{} \to
\fcci{}{}$.

\subsubsection{End of proof}
Since the differential in $\fcci{}{}$ is $0$, $\fcci{}{}$ coincides
with its cohomology algebra $\H$.  Finally we have built a sequence of
quasi-isomorphic sheaves of dg-algebras $\fC{}{} \to \fcun{}{} \from
\fcde{}{} \to\cdots \to \H$, as required.

To conclude the proof of proposition~\ref{prop:formalite} we remark
that in the above sequence of quasi-isomorphisms, each of the
intermediate sheaves, say $\fD$, is defined as a sum $\fD =
\oplus_{(\alpha,\beta) \in A^2} \fD^{\alpha\beta}$, so that, for
$\alpha\in A$, $\fD_\alpha = \oplus_{\alpha_1} \fD^{\alpha_1\alpha}$
has a natural structure of $\fD$-module.  It follows that, for a
quasi-isomorphism $\fD \to \fD'$ in the above sequence, the
equivalence of categories between $\D_{\fD}$ and $\D_{\fD'}$ sends the
$\fD$-module $\fD_\alpha$ to a $\fD'$-module isomorphic to
$\fD'_\alpha$.  This shows that $N_\alpha$ corresponds to $\H_\alpha$.

\section{proof of theorem~\ref{thm}}
\label{findemo}
Now we prove theorem~\ref{thm} using the equivalence of
proposition~\ref{prop:formalite} between $\D^b_{G}(X_Z)$ and $\D_\H
\langle \H_\alpha, \alpha \in A \rangle$. In view of this equivalence
and the equality $\H = \oplus_{\alpha\in A} \H_\alpha$, the algebra
$\E$ of the theorem is isomorphic to $\Ext^\cdot_{\D_{\H}}(\H,\H)$ and
$\E_\O^\rho \simeq \Ext^\cdot_{\D_{\H}}(\H,\H_\alpha)$, for $\alpha =
(\O,\rho)$.  We have to prove that $\D_\H \langle \H_\alpha \rangle$
is equivalent to $\D_\E \langle \E_\O^\rho \rangle$. We recall the
following construction:

\smallskip 1) Let $\fAintro$ be a sheaf of dg-algebras on a finite set
$I$ (hence $\M_\fAintro$ has enough $K$-projectives). Let
$M^1,\ldots,M^r$ be $\fAintro$-modules, $P^i \to M^i$ a $K$-projective
resolution of $M^i$ and $P=\oplus_i P^i$. The composition of morphisms
induces a structure of dg-algebra on $R = \Hom^\cdot(P,P)$ (see
section~\ref{eq_der_cat} for the definition of
$\Hom^\cdot(\cdot,\cdot)$). We have a functor $F:\M_\fAintro \to
\M_R$, $M\mapsto \Hom^\cdot(P,M)$, which sends quasi-isomorphisms to
quasi-isomorphisms because $P$ is $K$-projective. Hence it induces
$F:\D_\fAintro \to \D_R$. We set $N^i = F(M^i)$.  Then $F$ restricts
to an equivalence of categories between $\D_\fAintro \langle M^i
\rangle$ and $\D_R \langle N^i \rangle$ (for example this is a very
special case of~\cite{K94}, theorem 4.3).

\smallskip 2) We apply this to $\D_\H \langle \H_\alpha \rangle$.
Since $\H_\alpha$ is a direct summand of $\H$, $(\H_\alpha)_{U_F}$ is
$K$-projective, for any $F\in \F$. Moreover, since $\F$ satisfies (i)
of assumptions~\ref{assumptionstratification}, we see, as in the proof
of proposition~\ref{prop:eqL'M}, that the \v{C}ech resolution (where
we fix a total order on $\F$):
\begin{equation}
  \label{eq:resolHalpha}
P_\alpha \quad = \quad \cdots\to \bigoplus_{F_1<\cdots<F_k\in\F}
(\H_\alpha)_{U_{F_1}\cap\ldots\cap U_{F_k}} \to \cdots \to
\bigoplus_{F\in\F} (\H_\alpha)_{U_{F}} \to 0.
\end{equation}
is a $K$-projective resolution of $\H_\alpha$.  We set $P =
\oplus_{\alpha} P_\alpha$ and $R= \Hom^\cdot(P,P)$.
Following~\cite{L95} we use the fact that the differential of $\H$ is
$0$ and define $\H'$ to be the sheaf of non-graded algebras underlying
$\H$. Replacing $\H$ by $\H'$ above, we define similarly $\H'_\alpha$,
$P'_\alpha$, $P'$, $R'$ (we note that $P'_\alpha$ still is a
$K$-projective resolution of $\H'_\alpha$).  The algebras $R$ and $R'$
are canonically isomorphic as differential algebras, but they do not
have the same graduation. If we write $P$ as a double complex
$P=\oplus_{i,j} P^{ij}$ and set $Q_{ij}^{kl} =
\Hom_\C(P^{ij},P^{kl})$, then $R^d = R \cap (\oplus_{k+l = i+j+d}
Q_{ij}^{kl})$ and $R'^d = R\cap (\oplus_{k = i+d} Q_{ij}^{kl})$.
\begin{claim}
  \label{concentre_en0}
  the dg-algebra $R'$ is concentrated in degree $0$.
\end{claim}
We will prove this below. Let us see why it implies the theorem.  We
set $R_0 = \tau_{\leq 0} R' = \cdots \oplus R'^{-2} \oplus R'^{-1}
\oplus \ker d_0$. This is a differential sub-algebra of $R'$ (or of
$R$ as well) and the claim implies that we have quasi-isomorphisms of
differential algebras:
$$
R' \xfrom{u} R_0 \xto{v} H(R'), \qquad R \xfrom{u'} R_0 \xto{v'}
H(R).
$$
In view of the decompositions of $R$ and $R'$ above, we have $R'^d
= \oplus_n (R^n\cap R'^d)$. Hence we may endow $R_0$ with the
graduation induced by the embedding $u'$. Then $v'$ is a graded
morphism: indeed $v'$ is the composition of the projection from $R_0$
to $R'^0 \cap R_0$ and the projection from $\ker d_{R}$ to $H(R)$;
both morphisms are graded and so is $v'$. Now we just remark that
$H(R) =\E$ by definition. Moreover, the functor from $\D_\H$ to $\D_R$
sends $P_\alpha \simeq \H_\alpha$ to $\Hom^\cdot(P,P_\alpha)$.  Since
this last object is a summand of $R$, we see that, in the equivalence
from $\D_R$ to $\D_{H(R)}$, it is sent to its cohomology,
$\Ext^\cdot_{\D_\H}(P,P_\alpha) \simeq \Ext^\cdot_{\D_\H}(\H,\H_\alpha) =
\E_\O^\rho$.  Summing up we obtain an equivalence between $\D_\H
\langle \H_\alpha \rangle$ and $\D_\E \langle \E_\O^\rho \rangle$, as
desired.

\begin{proof}[Proof of claim~\ref{concentre_en0}]
  We have to prove that $\Hom_{\D_{\H'}}(\H',\H'[p])$ is $0$, for
  $p\not=0$. We use the \v{C}ech resolution~\eqref{eq:resolHalpha},
  with $\H'$ instead of $\H_\alpha$. Since $\Hom^\cdot(\H'_{U_F},\H')
  \simeq \sect(U_F; \H')$, we obtain
  $$
  \RHom(\H',\H') \quad \simeq\quad 0\to \bigoplus_{F\in\F}
  \sect(U_F;\H') \to \bigoplus_{F_1<F_2\in\F} \sect(U_{F_1}\cap
  U_{F_2};\H') \to\cdots.
  $$
  Since the open sets $U_{F_1}\cap\ldots\cap U_{F_r}$ are
  fundamental open sets, the functor of sections over them is exact. Hence
  the above resolution computes $H^\cdot(\F;\H')$ and the claim
  follows from the next lemma.
\end{proof}
\begin{lemma}
  For any $G$-stable open subset $V$ of $X_Z$, and $U= \pxf(V)$, we
  have $H^i(U;\H') = 0$ for $i>0$, where $\H'$ is the non-graded sheaf
  underlying $\H$.
\end{lemma}
\begin{proof}
  By proposition~\ref{prop:formalite}, we have an isomorphism $\H'
  \simeq \fcci{}{}$, and, by definition, $\fcci{}{} =
  \oplus_{(\alpha,\beta)\in A^2} \fcci{\alpha\beta}{}$.  Hence it
  sufficient to prove the vanishing of $H^i(U;\fcci{\alpha\beta}{})$.
  We fix $(\alpha,\beta)\in A^2$ for the remainder of the proof and
  set for short $\fD= \fcci{\alpha\beta}{}$.
  
  \smallskip
 
  1) Let us set $V_\alpha = X_Z \setminus \bigcup_{v\in \Delta'_\alpha}
  D_v$, $V_\beta = X_Z \setminus \bigcup_{v\in \Delta'_\beta} D_v$,
  $V_{\alpha\beta} = V_\alpha \cap V_\beta$, $U_{\alpha\beta} =
  \pxf(V_{\alpha\beta})$. We let $j:U_{\alpha\beta} \to \F$ be the
  inclusion.  We first prove that $\fD = j_*
  (\fD|_{U_{\alpha\beta}})$.  For a face $F\in \F$, repeated
  applications of hypothesis (iii) of
  assumptions~\ref{assumptionstratification}, show that there exists a
  face $F'$ such that $U_F \cap U_{\alpha\beta} = U_{F'}$.  Our claim
  is equivalent to the fact that for any $F$, the restriction map
  $\fD_F \to \fD_{F'}$ is an isomorphism. We have already seen, by
  definition of $\Delta'_\alpha$, $\Delta'_\beta$ and
  lemma~\ref{lem:monodromy}, that $L_\alpha \simeq
  (L_\alpha)_{V_\alpha}$, $L_\beta \simeq \rsect_{V_\beta}(L_\beta)$.
  It follows that $\Rhom(L_\alpha, L_\beta) \simeq
  \rsect_{V_{\alpha\beta}} \Rhom(L_\alpha, L_\beta)$.
  Using~\eqref{eq:germecohomo}, this implies the desired result:
  \begin{multline*}
    \fD_F \simeq
    \Ext^\cdot_{\D_G(X_Z)}(L_\alpha|_{V_F},L_\beta|_{V_F}) \simeq
    H^\cdot_G(V_F; \Rhom(L_\alpha,L_\beta)) \\
    \simeq H^\cdot_G(V_F \cap V_{\alpha\beta};
    \Rhom(L_\alpha,L_\beta)) \simeq \fD_{F'}.
  \end{multline*}
  Now we remark that the functor $j_*$ is exact, because, for any
  sheaf $B$ on $\F$, any face $F$, we have $(j_*B)_F \simeq B_{F'}$,
  for $F'$ satisfying $U_F \cap U_{\alpha\beta} = U_{F'}$ as above.
  Hence we have $H^i(U; \fD) \simeq H^i(U\cap U_{\alpha\beta};\fD)$.
  
  \smallskip

  2) This means that we may assume $V\subset V_{\alpha\beta}$.  We
  prove the result by induction on the number of $G$-orbits in $V$.
  If $V$ consists of the open orbit of $X_Z$, then $U$ is the
  fundamental open set $U_F$, with $F= F_{\emptyset,\emptyset}$. Hence
  $\sect(U;\cdot) = (\cdot)_F$ is exact and we are done.
  
  Now let us assume that $V=W\sqcup \O_\Delta$, where $W \subsetneq V$
  is a $G$-stable open subset of $X_Z$. By induction the result is true
  for $W$.  Let $F_{\Delta,J_\Delta}$ be the closed face of
  $\O_\Delta$. Since $U_{F_{\Delta,J_\Delta}}$ contains $\O_\Delta$,
  we have $\pxf(V)=\pxf(W) \cup U_{F_{\Delta,J_\Delta}}$. Setting $U'
  = \pxf(W) \cap U_{F_{\Delta,J_\Delta}}$ and using the Mayer-Vietoris
  sequence, it is sufficient to prove:
  \begin{equation}
    \label{eq:annulationcohom}
      \forall i>0,\:\: H^i(U'; \fD)=0, \quad \text{and} \quad
      H^0(U_{F_{\Delta,J_\Delta}}; \fD) \to H^0(U';\fD ) \:\text{is
        surjective.}
  \end{equation}

  \smallskip

  3) We compute $H^\cdot(U'; \fD)$ with the help of a \v{C}ech
  covering.  Remember that $U_{F_{\Delta',J'}} \subset
  U_{F_{\Delta,J}}$ if and only if $\Delta'\subset \Delta$ and $J
  \subset J'$.  For $\Delta' \subset \Delta$, we have $\O_\Delta \cap
  U_{F_{\Delta',J'}}= \emptyset$ if and only if $\Delta' \not=
  \Delta$; but this implicitly assume that $F_{\Delta',J'}$ is a face,
  i.e. $\O_{\Delta'} \cap C_{J'}$ is non-empty.  For any
  $\Delta'\subset \Delta$, we have indeed $F_{\Delta',J_\Delta} \not=
  \emptyset$: this is easily seen if $X_Z = X$, in which case $\F$ is
  the set of faces of $[0,1]^l$. This implies the general case because
  $C_{X_Z,J_\Delta} = \pxzx^{-1}\pxzx(C_{X,J_\Delta})$ (where $\pxzx$
  is the map $X_Z\to X$).  It follows that
  \begin{equation}
    \label{eq:recouvrement}
  \textstyle
  U' = \bigcup_{\Delta'\subsetneq \Delta} U_{F_{\Delta',J_\Delta}}.
  \end{equation}
  By~\eqref{intersectionUF}, this covering is stable by taking
  intersections. Since it consists of fundamental open sets, on which
  the functor of sections is exact, it can be used to compute
  $H^\cdot(U';\cdot)$. For this, we have to know $\fD(U_{F})$. By
  definition, for any face $F= F_{\Delta_1,J_1} \in \F$:
  $$
  \fD(U_F) = \fD_F = \fcsun{\alpha\beta}{F} \otimes
  \fckde{\alpha\beta}{F} .
  $$
  Let us describe more precisely the components of the tensor
  product. Recall that $F \in U_{\alpha\beta}$; hence either $F
  \not\in \F_{\alpha\beta} \cup \F'_{\alpha\beta}$ or $F\in
  \F_{\alpha\beta}$ (see~\eqref{eq:deffabfprimeab}).  In the first
  case we have $\fckde{\alpha\beta}{F} =
  \fcsun{\alpha\beta}{F} = 0$. In the second case, by
  definition~\eqref{eq:defCS1} we have
  $$
  \fcsun{\alpha\beta}{F_{\Delta_1,J_1}} = \C[X_v;\: v\in
  \Delta_1],
  $$
  and, by~\eqref{eq:deffCK2}, $\fckde{\alpha\beta}{F}$
  only depends on $J_1$, say $\fckde{\alpha\beta}{F} =
  M_{J_1}$.  These descriptions of $\fcsun{\alpha\beta}{F}$
  and $\fckde{\alpha\beta}{F}$ assume that $F_{\Delta_1,J_1}$
  is a face. We have seen that this is the case for
  $F_{\Delta',J_\Delta}$ with $\Delta'\subset \Delta$.

  \smallskip
  
  4) Since, in the covering~\eqref{eq:recouvrement}, all faces have
  the same ``$J$-index'', we obtain $H^i(U';\fD) = M(J_\Delta)\otimes
  H^i(C^\cdot)$, where $C^\cdot =C^\cdot(S_\Delta)$ is the following
  complex.  We let $S_\Delta$ be the set of subsets $\Delta'$ of
  $\Delta$ such that $(\Delta_\alpha \cup \Delta_\beta) \subset
  \Delta' \subsetneq \Delta$ and we consider any total order on
  $S_\Delta$:
  $$
  C^\cdot(S_\Delta) \:\: =\:\: 
  0\to \bigoplus_{\makebox[0cm]{$\scriptstyle
      \Delta'_1 \in S_\Delta$}} \C[X_v;\: v\in \Delta'_1 ] \to
  \bigoplus_{\makebox[0cm]{$\scriptstyle \Delta'_1 < \Delta'_2 \in
      S_\Delta$}} \C[X_v;\: v\in \Delta'_1\cap \Delta'_2] \to\cdots
  .
  $$
  We have a bijection between $S_\Delta$ and the set, $S_\Phi$, of
  strict subsets of $\Phi = \Delta \setminus (\Delta_\alpha \cup
  \Delta_\beta)$. Setting $C_1 = \C[X_v ;\: v \in \Delta_\alpha \cup
  \Delta_\beta]$ and $C^\cdot_2 = C^\cdot(S_\Phi)$, we obtain $C^\cdot
  = C_1 \otimes C^\cdot_2$.
  
  Hence $H^i(U';\fD) = M(J_\Delta)\otimes C_1 \otimes H^i(C^\cdot_2)$.
  Since $H^0(U_{F_{\Delta,J_\Delta}}; \fD) \simeq M(J_\Delta) \otimes
  \C[X_v;\: v\in \Delta]$, \eqref{eq:annulationcohom} will follow from
  \begin{equation}
    \label{eq:annulationcohombis}
      \forall i>0,\:\: H^i( C^\cdot_2)=0, \quad \text{and} \quad
       \C[X_v;\: v \in \Phi] \to H^0(C^\cdot_2 ) \:\text{is
        surjective.}
  \end{equation}

  \smallskip

  5) We may interpret $C^\cdot_2$ as another \v{C}ech complex: we
  consider the following sheaves on the quadrant $\q =\R_{\geq
    0}^{\Phi}$,
  $$
  M_1 = \oplus_{k\in \Phi} \C_{\{x_k =0\}}, \quad \text{$M_2$
    associated to $O \mapsto Sym(M_1(O))$},
  $$
  where $O\subset \q$ is open and $Sym(\cdot)$ denotes the
  symmetric algebra.  We also consider the open subsets of $\q$, for
  $\Phi' \subset \Phi$, $U_{\Phi'} = \{ x_k >0; \: k\in \Phi \setminus
  \Phi' \}$. The $U_{\Phi'}$, for $\Phi'\subsetneq \Phi$, give a
  covering of $\q \setminus \{0\}$, and
  $$
  \sect( U_{\Phi'}; M_2) \simeq \C[X_v;\: v\in \Phi' ].
  $$
  Hence $C_2^\cdot$ is isomorphic to the \v{C}ech complex
  $\mathcal{C}^\cdot( (U_{\Phi'})_{\Phi' \subsetneq \Phi}; M_2)$.  Now
  $M_2$ is constructible for the stratification of $\q$ by the strata
  $\q_{\Phi'} = \{x_k =0,\, k\in \Phi'; \: x_k>0,\, k\not\in \Phi'
  \}$, $\Phi' \subset \Phi$. Each open $U_{\Phi'}$ is contractible to
  a point by a homotopy preserving the closures of strata; hence
  $H^i(U_{\Phi'};M_2) =0$, for $i>0$.  It follows that $\forall i$,
  $H^i(C^\cdot_2) \simeq H^i( \q \setminus \{0\} ; M_2)$.
  
  For $\Phi_1 \subset \Phi$ we have of course $\otimes_{k\in \Phi_1}
  \C_{\{x_k =0\}} \simeq \C_{\{x_k =0; k\in\Phi_1\}}$. Hence $M_2$ is
  a sum of sheaves of the type $\C_{\{x_k=0; k\in\Phi_1\}}$. For each
  of them, say $N = \C_{\{x_k=0; k\in\Phi_1\}}$, we have $\forall
  i>0$, $H^i( \q \setminus \{0\} ; N) =0$ and the map $H^0( \q ; N)
  \to H^0( \q \setminus \{0\} ; N)$ is surjective. By additivity, both
  assertions are true with $M_2$ instead of $N$. We also have $H^0( \q
  ; M_2) \simeq \C[X_v;\: v \in \Phi]$. Hence we
  obtain~\eqref{eq:annulationcohombis} and the lemma is proved.
\end{proof}

\end{document}